\documentclass[article,onefignum,onetabnum]{siamart171218}

\usepackage{lipsum}
\usepackage{amsfonts}
\usepackage{graphicx}
\usepackage{epstopdf}
\usepackage{verbatim}
\usepackage{algorithmic}
\ifpdf
\DeclareGraphicsExtensions{.eps,.pdf,.png,.jpg}
\else
\DeclareGraphicsExtensions{.eps}
\fi

\newsiamremark{remark}{Remark}
\newsiamremark{assumption}{Assumption}
\newsiamremark{hypothesis}{Hypothesis}
\crefname{hypothesis}{Hypothesis}{Hypotheses}
\newsiamthm{claim}{Claim}

\headers{Second-order Guarantees Of Distributed Gradient Algorithms}{A. Daneshmand, G. Scutari, and V. Kungurtsev}

\title{Second-order Guarantees Of Distributed Gradient Algorithms\thanks{Part of this work has been presented at the 56th Annual Allerton Conference on Communication, Control, and Computing \cite{8636044}.
		\funding{The work of Daneshmand  and Scutari has been supported by the USA National Science Foundation under Grants  CIF 1564044, CIF 1719205, and CMMI 1832688; and the Army Research Office under Grant W911NF1810238.
			The work of Kungurtsev was supported by the OP VVV project
CZ.02.1.01/0.0/0.0/16 019/0000765 ``Research Center for Informatics''.}}}

\author{Amir Daneshmand\thanks{School of Industrial Engineering, Purdue University, West-Lafayette, IN, USA
		(\email{adaneshm@purdue.edu}, \url{http://web.ics.purdue.edu/\string~adaneshm/}).}
	\and Gesualdo Scutari\thanks{School of Industrial Engineering, Purdue University, West-Lafayette, IN, USA 
		(\email{gscutari@purdue.edu}, \url{https://engineering.purdue.edu/\string~gscutari/}).} \and  Vyacheslav Kungurtsev \thanks{Department of Computer Science, Czech Technical University in Prague, Department of Computer Science, Faculty of Electrical Engineering, Prague, Czech Republic   (\email{vyacheslav.kungurtsev@fel.cvut.cz}).}
}

\usepackage{amsopn}
\DeclareMathOperator{\diag}{diag}

\usepackage{soul} 

\newcommand{\norm}[1]{\left\lVert#1\right\rVert}
\newcommand{\trace}[1]{\text{tr}\left(#1\right)}

\newcommand\floor[1]{\lfloor#1\rfloor}

\usepackage{enumitem}

\usepackage{amssymb}
\usepackage{hyperref}
\usepackage{amsmath}
\usepackage{multirow}

\newcommand{\bt}{\boldsymbol{\theta}}
\newcommand{\bx}{\mathbf{x}}
\newcommand{\by}{\mathbf{y}}
\newcommand{\bT}{\boldsymbol{\Theta}}

\def \algname/{{\sc DOGT}}

\usepackage{mathrsfs}

\begin{document}
	
\maketitle

\begin{abstract}
	We consider distributed smooth nonconvex unconstrained optimization over networks, modeled as a   connected graph. We examine the behavior  of   distributed gradient-based  algorithms  near strict saddle points. Specifically, 
	we establish that (i) the renowned Distributed Gradient Descent (DGD) algorithm likely converges to a \emph{neighborhood} of a  Second-order Stationary (SoS)   solution;      
	and (ii)   the more recent  class of distributed algorithms based on gradient tracking--implementable also over digraphs--likely converges to   \emph{exact} SoS solutions,  
	thus avoiding   strict saddle-points. 
	Furthermore, new convergence rate results to first-order  critical points is established for the latter class of algorithms.
	\vspace{-0.2cm}

\end{abstract}

\begin{keywords}
	distributed optimization, decentralized methods, nonconvex optimization, second-order guarantees, distributed gradient descent, stable-manifold theorem, gradient tracking \vspace{-0.2cm}
\end{keywords}

\vspace{-0.2cm}

\section{Introduction}\label{Sec_intro}
We consider smooth unconstrained nonconvex optimization over networks, in the following form:\vspace{-0.2cm}
\begin{equation}
\min_{ {\bt}\in \mathbb{R}^m}\, F(\bt)\triangleq \sum_{i=1}^n f_i(\bt),
\tag{P}\label{eq:P}\vspace{-0.2cm}
\end{equation}
where $n$ is the number of agents in the network;  and  $f_i:\mathbb{R}^m\rightarrow\mathbb{R}$ is the cost function of agent $i$, assumed to be smooth and known only to agent $i$. 
Agents are connected through a communication network, modeled as a (possibly directed,  strongly) connected graph. No specific topology  is assumed for the graph (such as star or hierarchical structure). In this setting, agents seek  to cooperatively solve
Problem \eqref{eq:P}   by exchanging information with their immediate neighbors in the network.

Distributed \emph{nonconvex} optimization  in the form   \eqref{eq:P} has found a wide range of applications in several areas, including   network  information processing, machine learning,  communications, and multi-agent control;  see, e.g., \cite{Scutari_Ying_LectureNote}. For instance, this is the typical scenario of in-network data-intensive (e.g., sensor-network) applications wherein data are scattered across the  agents (e.g., sensors, clouds, robots), and     the sheer volume and spatial/temporal disparity of  data  render centralized processing and storage infeasible or inefficient. Communication networks modeled as \emph{directed} graphs capture  simplex communications between adjacent nodes.   {This is the case, e.g., in several wireless (sensor) networks wherein nodes  transmit at different power  and/or communication channels are not symmetric.}

\textbf{Main objective:} We call $\bt$ a critical point of $F$ if $\nabla F(\bt)=\mathbf{0}$; a critical point  $\bt$ is a \emph{strict saddle} of $F$ if  $\nabla^2 F(\bt)$ has at least one negative eigenvalue; and it is  a \emph{Second-order  Stationary} (SoS) solution if  $\nabla^2 F(\bt)$ is positive semidefinite. Critical points that are not minimizers are of little interest in the nonconvex setting.  It is thus desirable to consider    methods for  \eqref{eq:P} that are not attracted to such points. 
When $F$ has a favorable structure, stronger guarantees can   be claimed. For instance, a wide range of  salient objective functions arising from applications in  machine learning and signal processing   have been shown to enjoy  the so-called \emph{strict saddle} property: all the critical points of $F$ are either strict saddles or local minimizers. Examples include principal component analysis and fourth order tensor factorization \cite{pmlr-v40-Ge15}, low-rank matrix completion \cite{Ge_localMinsSpurious_NIPS16}, and some instances of neural networks \cite{NIPS2016_6112}, just to name a few. In all these cases, converging to SoS solutions--and thus circumventing strict saddles--guarantees finding  a local minimizer.  

This paper  studies for the first time  second-order guarantees of two renowned distributed gradient-based algorithms for Problem \eqref{eq:P}, namely: the   Distributed Gradient Descent (DGD)  \cite{Nedic2009,Nedic2010} and the family of distributed algorithms based on gradient-tracking \cite{dilorenzo2015distributed,NEXT,xu2015augmented}. The former is implementable on undirected graphs while  the latter is suitable also for directed graphs.    Convergence of these schemes applied to {\it convex} instances of \eqref{eq:P} is well understood; however,  less is known  in the {\it nonconvex} case,  let alone second-order guarantees; the relevant works 
are discussed next.\vspace{-0.1cm}

\subsection{Literature review} { Recent years have witnessed many studies proving  asymptotic solution-   and  {convergence}  rate-guarantees of a variety of algorithms for specific classes of nonconvex  optimization problems (e.g., satisfying suitable regularity conditions); a good overview can be found in   \cite{chi2019nonconvex}.   Since these analyses are heavily tailored to specific applications and it is unclear how to generalize them to a wider class of nonconvex functions, we omit further details and discuss next only  results of centralized and distributed algorithms for {\it general} nonconvex instances of \eqref{eq:P}. }\vspace{-0.2cm}
\subsubsection{Second-order guarantees of centralized optimization algorithms}
Second-order guarantees of  centralized  solution methods for general nonconvex optimization  \eqref{eq:P}  have been extensively studied in the literature. 

\textbf{Hessian-based methods:}  Algorithms based on \emph{second-order} information   have long been known to converge to SoS solutions of \eqref{eq:P}; they rely on computing the Hessian to distinguish between first- and second-order stationary points.   The classical cubic-regularization  \cite{Griewank_technical_report, Nesterov2006,  Cartis2011_part1, Cartis2011_part2, Agarwal:2017:FAL:3055399.3055464} and trust region   (e.g. \cite{doi:10.1137/0904038,Powell1984, Curtis2017,Dauphin:2014:IAS:2969033.2969154}) methods can  provably find {\it approximate} SoS solutions  in polynomial time (by approximate SoS we mean  $\boldsymbol{\theta}$ such that  $||\nabla F(\boldsymbol{\theta})||\leq \epsilon_g$ and $\lambda_{\min}(\nabla^2 F(\boldsymbol{\theta}))\geq  -\epsilon_h$, for small $\epsilon_g,\epsilon_h>0$); they however require access to the  full Hessian matrix.   A recent line of works \cite{Carmon_2018, Agarwal_et_al_2017, Carmon_Duchi_2016} show that the requirement of full Hessian access can be relaxed to Hessian-vector products in each iteration, hence solving simpler sub-problems per iteration, but at the cost of requiring more iterations to reach approximate SoS solutions.
\textbf{First-order methods:} For general nonconvex problems, Gradient Descent (GD) is   known to find a stationary point  in polynomial time \cite{Nesterov04}.      
In \cite{pmlr-v49-lee16}, it was proved that randomly initialized GD with a fixed step-size   converges to SoS solutions almost surely.  
The elegant analysis of \cite{pmlr-v49-lee16}, leveraging tools from the theory of dynamical systems (e.g., the Stable Manifold Theorem), has been later extended in a number of follow-up works establishing same kind of second-order guarantees of a variety of  first-order methods, including   the proximal point algorithm, block  coordinate descent,  mirror descent  \cite{Jordan_FOM_AvoidsSaddles};   the heavy-ball method and the Nesterov's accelerated method \cite{SWright_AGDavoidssaddles};   block coordinate descent and alternating minimization \cite{Li2019}; and   a primal-dual optimization procedure for solving linear equality constrained nonconvex optimization problems  \cite{Mingyi_SOAlgs_AvoidsSaddles}. 
These results are all asymptotic in nature and it is unclear whether polynomial convergence rates can be obtained for these methods. In \cite{du2017gradient} it was actually proven that, even with fairly natural random initialization schemes and for non-pathological functions, GD can be significantly slowed down by saddle points, taking exponential time to escape. Recent work has analyzed variations of GD that include stochastic perturbations. It has been shown that when perturbations are incorporated into GD at each step the resulting algorithm can escape strict saddle points in polynomial time  \cite{pmlr-v40-Ge15}; the same conclusion was earlier  established in  \cite{Pemantle90} for stochastic gradient methods, although without escape time guarantees.  It has also been shown that episodic perturbations suffice; in particular, \cite{pmlr-v70-jin17a} introduced an algorithm that occasionally adds a perturbation to GD, and proved that 
the number of iterations to escape saddle points  depends only poly-logarithmically on dimension (i.e., it is   nearly dimension-independent).  Fruitful follow-up results show that other   first-order    perturbed algorithms escape from strict saddle points efficiently \cite{Jin2018,Songtao2019}.

\subsubsection{Distributed algorithms for \eqref{eq:P} and guarantees}\label{sec:intro_distributed} Distributed algorithms for {\it convex} instances of \eqref{eq:P} have a long history; less results are available for  nonconvex objectives. Since the focus on this paper is on nonconvex problems, next, we mainly comment on distributed algorithms for minimizing    nonconvex objectives.  
	
\noindent  \textbf{$\bullet$ DGD and its variants:} DGD   (and its variants)  is  unquestionably among the first and most studied  decentralizations of the gradient descent algorithm for \eqref{eq:P} \cite{Nedic2009,Nedic2010}.    The instance of DGD considered in this paper reads:  given $\mathbf{x}_i^0\in \mathbb{R}^m$, $i\in[n]$,   \vspace{-0.4cm}\begin{equation}
\label{DGD_intro}
\mathbf{x}_i^{\nu+1}=\sum_{j=1}^n D_{ij}\,\mathbf{x}_j^\nu-\alpha \nabla f_i(\mathbf{x}_i^\nu),\quad i\in[n],\vspace{-0.3cm}
\end{equation}
where $\mathbf{x}_i^\nu$ is the agent $i$'s  estimate at iteration $\nu$ of the vector variable $\bt$;  $\{D_{ij}\}_{i,j}$ are suitably chosen  set of nonnegative weights (cf. Assumption \ref{DoublyStochastic_D}), matching the graph topology (i.e., $D_{ij}>0$ if there is a link between node $i$ and $j$, and $D_{ij}=0$  otherwise); and $\alpha>0$ is the step-size.  Roughly speaking, the update of each agent $i$ in \eqref{DGD_intro} is the linear combination of two components: i) the gradient $\nabla f_i$ evaluated at the agent's latest iterate (recall that agents do not have access to the entire gradient $\nabla F$); and ii) a convex combination of the current  iterates of the neighbors of   agent $i$ (including agent $i$ itself). The latter term (a.k.a. consensus step) is instrumental to asymptotically enforcing agreement  among the agents' local variables.
	
When each $f_i$ in \eqref{eq:P} is (strongly) {\it convex}, convergence  of   DGD   is well understood. With a diminishing step-size,   agents' iterates  converge to a consensual {\it exact} solution;   if a constant step-size is used,  convergence is generally faster but only   to a neighborhood of the solution, and exact consensus is not achieved.  
When \eqref{eq:P} is {\it nonconvex}, the available convergence guarantees are weaker. 
In  \cite{WYin_ncvxDGD_SIAM2018} it was shown that if a constant step-size is employed, every limit point $(\bx_1^{\infty},\ldots, \bx_n^{\infty})$ of the sequence generated by \eqref{DGD_intro} satisfies $\sum_{i=1}^n \nabla_{x_i} f_i(\bx_i^\infty)=\mathbf{0}$; the limit points of agents' iterates are not consensual; asymptotic consensus is achieved only using a diminishing step-size.
Since in general $f_i$ are all different, such limit points are {\it not} critical points of $F$.  Nothing is known about the connection of the critical  points of $\sum_{i=1}^n f_i(\bx_i)$  and those of $F$, {\it let alone its second-order guarantees}.   { A first contribution of this paper is to establish  second-order guarantees of  DGD  \eqref{DGD_intro}  applied to  \eqref{eq:P} over undirected graphs}. 

Several extensions/variants of the vanilla DGD followed the seminal works \cite{Nedic2009,Nedic2010}. The projected (stochastic) DGD for nonconvex constrained instances of \eqref{eq:P} was proposed in \cite{bianchi2013convergence}; with a diminishing step-size,   the algorithm converges to a stationary solution of the problem (almost surely, if noisy instances of the local gradients are used). The extension  of DGD to 
{\it digraphs} was   studied in \cite{Nedic2015}  for convex unconstrained optimization, and later  extended in \cite{tatarenko2017nonconvex} to nonconvex objectives. The algorithm, termed push-sum DGD, combines a local gradient step with the push-sum algorithm \cite{benezit2010weighted}. When    a diminishing step-size is employed, push-sum DGD 
converges to an exact  stationary solution of \eqref{eq:P};    and its   noisy perturbed version almost surely  converges   to local minimizers, provided that  $F$ does not have any saddle point \cite{tatarenko2017nonconvex}.  { To our knowledge, no other guarantees are known for DGD-like algorithms in the nonconvex setting. In particular, it is   unclear whether DGD \eqref{DGD_intro} escapes strict saddles of $F$. }

\smallskip 

\noindent\textbf{$\bullet$ Gradient tracking-based methods:}\label{DOGT_method_intro}
To cope with the speed-accuracy dilemma of DGD, \cite{dilorenzo2015distributed,NEXT} proposed a new class of   distributed gradient-based   methods that converge to an \emph{exact} consensual
solution of nonconvex (constrained) problems while using a \emph{fixed step-size}. The algorithmic framework, termed NEXT, 
introduces the idea of \emph{gradient tracking} to  correct  the
DGD direction  and cancel the steady state error in it while using a fixed step-size: each agent updates its own local variables along a surrogate direction that
tracks the   gradient  $\nabla F$ of the entire objective (the same idea was proposed independently in \cite{xu2015augmented} for convex unconstrained smooth problems).  
The generalization of NEXT to  digraphs--the SONATA algorithm--
was proposed in \cite{sun2016distributed,Scutari_Ying_LectureNote,SONATA_technical_report, YingDanScuConvergence18}, with \cite{SONATA_technical_report,YingDanScuConvergence18} proving   convergence of the agents' iterates to consensual stationary solutions of nonconvex problems at a sublinear rate. {\it No second-order guarantees have been established for these methods}.  Extensions of the SONATA family based on different   choices of the weight matrices   were later introduced in  \cite{xin2018linear,ANedich_GeoConvAlg_arxive_2018} for {\it convex} smooth unconstrained problems.  In this paper we consider the following family of distributed algorithms based on gradient tracking, which encompasses the majority of the above schemes  (see, e.g., \cite[Sec. 5]{SONATA_technical_report}), and refer to it as   
Distributed Optimization   with Gradient Tracking  (\algname/): 
\vspace{-0.25cm} 
\begin{align} 
\mathbf{x}_i^{\nu+1} &=\sum_{j=1}^n{R}_{ij}\bx_j^\nu\ -\alpha\,\mathbf{y}^\nu_i,  \label{eq:DOGT_x_update_intro} 
\\ 
\mathbf{y}^{\nu+1}_i &=\sum_{j=1}^n{C}_{ij}\by_j^\nu +\nabla f_i\big(\mathbf{x}^{\nu+1}_i\big)-\nabla f_i\big(	\mathbf{x}^\nu_i\big),\quad\text{(Gradient Tracking)}\label{eq:DOGT_track_update_intro}
\end{align}
where  $(R_{ij})_{i,j}$ and $(C_{ij})_{i,j}$ are suitably chosen nonnegative weights   compliant to the graph structure (cf. Assumption \ref{matrix_C_R}); and  $\mathbf{y}_i\in \mathbb{R}^m$ is   an auxiliary variable, controlled by agent $i$ via the update (\ref{eq:DOGT_track_update_intro}), which  aims at tracking locally   the gradient sum $\sum_i\nabla f_i(\mathbf{x}^\nu_i)$. Overall, the update (\ref{eq:DOGT_track_update_intro}) in conjunction with the consensus step in (\ref{eq:DOGT_x_update_intro}) 
is meant to  ``correct''  the local gradient direction $-\nabla f_i(\mathbf{x}^\nu_i)$ (as instead used in the DGD algorithm) and thus  nulls asymptotically the steady error $\nabla f_i(\mathbf{x}^\nu_i)-\nabla F(\mathbf{x}^\nu_i)$. This permits the use of a constant step-size $\alpha$ while still achieving {\it exact} consensus without penalizing the convergence rate. Another important difference between DOGT and DGD in (\ref{DGD_intro}) is that the former serves as a unified platform for distributed algorithms  applicable over both  undirected and directed graphs. Convergence of DOGT in the form \eqref{eq:DOGT_x_update_intro}-(\ref{eq:DOGT_track_update_intro}) when $F$ is nonconvex   remains an open problem, let alone second-order guarantees. A second contribution of this paper is to fill this gap and provide a first- and second-order convergence analysis of DOGT.


\noindent\textbf{$\bullet$ Primal-dual distributed algorithms:}	\label{PrimalDual_Algs_Subsec}
We conclude this literature review by commenting  on distributed algorithms for nonconvex \eqref{eq:P} using a primal-dual form \cite{Zhu:gm,hong2016decomposing,P-ProxPDA}. 
Because of their primal-dual nature, all these schemes  are implementable only  over {\it undirected} graphs.     
In \cite{Zhu:gm}  
a distributed approximate dual (sub)gradient algorithm, coupled with a consensus step is introduced. 
Assuming zero-duality gap, the algorithm is   proved to asymptotically find a pair of primal-dual solutions  of an auxiliary problem, which however    might not be   critical points of $F$; also, consensus is not guaranteed. No rate analysis is provided.   In \cite{hong2016decomposing}, a proximal primal-dual algorithm is proposed;  the algorithm,  termed Prox-PDA, employs either a constant or increasing penalty parameter (which plays the role of the step-size); a  sublinear convergence rate of a suitably defined primal-dual gap is proved. 
A perturbed  version of Prox-PDA, P-Prox-PDA,  was introduced  in \cite{P-ProxPDA}, which can also deal with nonsmooth convex, additive functions in the objective of     \eqref{eq:P}. P-Prox-PDA  converges to an $\epsilon$-critical point  (and thus also to {\it inexact} consensus), under a proper choice of the   penalty parameters  that depends on  $\epsilon$. A sublinear convergence rate is also proved. No second-order guarantees have been established for the above schemes. The only primal-dual algorithms we are aware of with provable convergence to SoS solutions is the one in \cite{Mingyi_SOAlgs_AvoidsSaddles}, proposed for    a linearly constrained nonconvex optimization problem. When linear constraints are used to enforce consensus, the primal-dual method \cite{Mingyi_SOAlgs_AvoidsSaddles} becomes distributed and applicable to Problem \eqref{eq:P}, but only for {\it undirected} graphs (DOGT is instead implementable also over digraphs). Second-order guarantees of such a scheme are established under slightly stronger assumptions than those required for DOGT (cf.  Remark \ref{DOGT_vs_ProxPDA_remark}, Sec. \ref{sec:SOS_main}). Finally, notice that, since \cite{Mingyi_SOAlgs_AvoidsSaddles} substantially differs from   DGD and DOGT--the former is a primal-dual scheme while the latter are primal methods--the convergence analysis put forth in  \cite{Mingyi_SOAlgs_AvoidsSaddles} is not applicable to DGD and DOGT. Since  DGD and DOGT  in their general form encompass two classic algorithms for distributed optimization,
the open problem of their second-order properties leaves a significant gap in the literature. 
	\vspace{-0.1cm}

\subsection{Major results} We establish for the first time second-order guarantees of DGD (\ref{DGD_intro}) and DOGT \eqref{eq:DOGT_x_update_intro}-\eqref{eq:DOGT_track_update_intro}. The main results are summarized next.\vspace{-0.2cm}
\subsubsection{DGD (\ref{DGD_intro})}\label{sec:DGD_intro_contribution} We prove that: \begin{description}
	\item[(i)] { For a sufficiently small step-size $\alpha$, agents' iterates $\{\bx^\nu\}$ generated by      (\ref{DGD_intro}) converge     to an $O(\alpha)$-critical point of $F$ for all   initializations--see  {Lemma \ref{L_criticalPoints_are_alpha_statPoints}}; neighborhood convergence to critical points is also established (cf.Theorem \ref{prop_DGD_local_convergence}).  This  complements the   convergence results in \cite{WYin_ncvxDGD_SIAM2018}; }
	\item[(ii)]  The average sequence $\{\overline{\bx}^\nu\triangleq (1/n) \sum_{i=1}^n \bx_i^\nu\}$  converges almost surely to a neighborhood  of a SoS solution  of \eqref{eq:P}, where the probability is taken over the initializations--see {Theorem \ref{SOG_of_DGD_thm}}.  
\end{description}
	
To prove \textbf{(ii)}, we employ a novel  analysis, which represents a major  technical contribution  of this work. In fact, existing techniques developed to established second-order guarantees of   the centralized GD are not readily applicable to  DGD--roughly speaking, this is due to the fact that   DGD (\ref{DGD_intro})  converges only to a neighborhood of critical points of $F$ [fixed points of (\ref{DGD_intro}) are not critical points of $F$]. We  elaborate next on this challenge and outline our analysis.  
		

The  elegant roadmap developed  in \cite{pmlr-v49-lee16,Jordan_FOM_AvoidsSaddles}   to establish second-order guarantees of the centralized GD builds on the  Stable Manifold theorem: roughly speaking,   fixed-points of the gradient map corresponding to strict saddles of the objective function  are ``unstable''    (more formally, the stable set\footnote{ Given   $\mathcal{X}\subseteq \mathbb{R}^m$,   $g:\mathbb{R}^m\to \mathbb{R}^m$, and the fixed-point iterate $\bx^{\nu+1}=g(\bx^\nu)$,  the stable set of $\mathcal{X}$ is    $\{\bx:\lim_\nu g^\nu(\bx)\in\mathcal{X}\}$, i.e., the set of initial  points such that  $\{\bx^{\nu}\}$ converges to a member of $\mathcal{X}$.} of strict saddles has zero measure),  implying almost sure convergence of GD iterates to SoS points  \cite[Corollary 2]{Jordan_FOM_AvoidsSaddles}. It is known that the DGD iterates (\ref{DGD_intro}) can be interpreted as instances of the GD applied to the following auxiliary function \cite{WYin_DGD_SIAM2016, WYin_ncvxDGD_SIAM2018}: denoting  $\bx\triangleq [\bx_1^\top,\ldots \bx_n^\top]^\top$, \vspace{-0.2cm}
\begin{equation}
L_\alpha (\mathbf{x})\triangleq \underbrace{\sum_{i=1}^n f_i(\bx_i)}_{ {\triangleq F_c(\bx)}}+\frac{1}{2\alpha}\sum_{i=1}^n\sum_{i=j}^n(e_{ij}-D_{ij})\bx_i^\top\bx_j,
\label{LyapunovFunc_intro}\vspace{-0.2cm}
\end{equation}
where $e_{ij}=1$ if there is an edge in the graph between agent $i$ and agent $j$; and $e_{ij}=0$ otherwise. 	Using \eqref{LyapunovFunc_intro}, (\ref{DGD_intro})  can be rewritten as: denoting  $\bx^\nu\triangleq [\bx_1^{\nu\,\top},\ldots \bx_n^{\nu\,\top}]^\top$,\vspace{-0.2cm}
\begin{equation}
\mathbf{x}^{\nu+1}=\mathbf{x}^{\nu}-\alpha\nabla L_\alpha(\mathbf{x}^\nu).\vspace{-0.2cm}
\label{DGD_update_GDupdate_intro}
\end{equation} 
One can then apply the above argument (cf. \cite[Corollary 2]{Jordan_FOM_AvoidsSaddles})  to \eqref{DGD_update_GDupdate_intro} and readily establish  the following result (see  Theorem \ref{L_criticalPoints_are_alpha_statPoints} for the formal statement) 
\begin{description}
	\item[Fact 1 (informal):] For sufficient small  $\alpha>0$, randomly initialized DGD \eqref{DGD_update_GDupdate_intro}  [and thus (\ref{DGD_intro})] converges almost surely to a second-order critical point of $L_\alpha$. 
\end{description}
Unfortunately, this result alone is not satisfactory, as  no connection is known between the   critical points   of $L_\alpha$ and those of $F$ (note that $L_\alpha:\mathbb{R}^{n\cdot m}\to \mathbb{R}$ whereas $F:\mathbb{R}^{m}\to \mathbb{R}$). To cope with this issue we prove the following two facts. 
\begin{description}
\item[Fact 2 (informal):] Every limit point $\overline{\bx}^\infty$ of the average sequence  $\overline{\bx}^\nu=1/n\sum_{i=1}^n \bx_i^\nu$ can be made arbitrarily close to a critical point of $F$ by using a sufficiently small $\alpha>0$ (Theorem \ref{prop_DGD_local_convergence});
\item[Fact 3 (informal):]  Whenever the limit point $\bar{\mathbf{x}}^\infty=1/n\,\sum_{i=1}^n\mathbf{x}^\infty_i$ belongs to a sufficiently small neighborhood  of a  strict saddle  of $F$,   ${\mathbf{x}}^\infty=[\bx_1^{\infty\top},\ldots,\bx_n^{\infty\top}]^\top$ must be a   {strict  saddle of  $L_\alpha$} (Proposition \ref{DGD_SSP_L_and_F_connect} and Corollary \ref{DGD_SSP_L_and_F_connect_2_AMIREM}).
\end{description}  
	The above three facts will then ensure that, for sufficiently small $\alpha>0$, with almost complete certainty, 
	$\{\bar{\mathbf{x}}^\nu\}$  will not get trapped in   a    neighborhood of a strict saddle of $F$--as $\mathbf{x}^\infty$ would be a strict saddle of $L_\alpha$--thus  landing  in a   neighborhood of a SoS solution of \eqref{eq:P}.



Facts 2 \& 3 above  are proved under a regularity condition on $F$ which recalls (albeit  slightly weaker than)  \cite{Gelfand_et_al_1991}. Roughly speaking,  the gradient flow over some   \emph{annulus} must be uniformly positive correlated with any  outward (from the origin) direction   (cf. Assumption \ref{Strip_assumption}). This condition is quite mild and is satisfied by functions arising, e.g.,  from several machine learning applications, including  distributed PCA, matrix sensing,  and  binary classification problems; see Sec. \ref{ProblemSetting_SubSe} for more details. Furthermore, this condition is also sufficient to prove convergence of DGD  without assuming the objective function to be globally $L$-smooth (but just locally L-smooth, $LC^1$ for short),  a requirement  that instead is common to   existing (first-order) convergence conditions of DGD. Notice that the loss functions arising from many of the aforementioned machine learning problems are not globally   $L$-smooth. 
	






\subsubsection{DOGT \eqref{eq:DOGT_x_update_intro}-\eqref{eq:DOGT_track_update_intro}} For DOGT, we establish the following three results.

	\begin{description}
		\item[(i)]  { When $F$ is nonconvex and the graph is either undirected or directed, it is proved that every limit point of the sequence generated by DOGT is  a critical point  of $F$. Furthermore, a  merit function,  measuring distance of the iterates from stationarity and consensus disagreement is introduced, and  proved to vanish at a sublinear rate--see  {Theorem \ref{Th_asympt_convergence}}. This extends    convergence results   \cite{ANedich_GeoConvAlg_arxive_2018, xin2018linear}, established only for convex functions.  To deal with nonconvexity, our   analysis builds on  a novel Lyapunov-like function  {[cf. \eqref{eq:def_L}]}, which properly combines  optimization error dynamics, consensus and tracking disagreements. While  these three terms alone do not ``sufficiently'' decrease along the iterates--as local optimization and consensus/tracking steps might act as competing forces--a suitable  combination of them,  as captured by the Lyapunov function, does monotonically decrease.}  
		\item[(ii)] When $F$  satisfies the Kurdyka-{\L}ojasiewicz (K\L) property \cite{MR0160856,Kurdyka1998} at any of its critical points,  convergence of the entire sequence to a critical point of $F$  is proved (cf.   {Theorem \ref{main_global_conv_thm}}), and a convergence rate is provided (cf.    {Theorem \ref{Rate_theorem}}). Although inspired by   \cite{Attouch2013}, establishing similar convergence results (but no rate analysis)  for  centralized first-order methods, our proof follows a different path building on the descent of the Lyapunov function introduced in (i), which  does not satisfy {\cite[conditions H1-H2]{Attouch2013}}); 
		see   Sec. \ref{Global_Conv_Subsec} for details. 


\item[(iii)] The sequence of iterates generated by DOGT   is shown to converge to  SoS solutions of \eqref{eq:P} almost surely, when initial points are randomly drawn from a suitably chosen   linear subspace--see  {Theorem \ref{Main_SOS_guarantee_long}}. This result is proved for undirected and directed networks. The proofs build on the stable manifold theorem, based upon the  interpretation of DOGT dynamics as   fixed-point iterates of a suitably defined map. 
The challenge in finding such a map is  ensuring that  the stable set of its undesirable fixed-points--those associated with the strict saddles of $F$--has measure zero in the  subspace where the initialization of DOGT takes place. Note that this subspace is not full dimensional.
			
%

\end{description}



{%
%

\smallskip

While our paper was under review after its initial arXiv posting \cite{Danarxiv2018} and its companion conference version \cite{8636044}, we became aware of a followup line of related works \cite{vlaski2019distributed_part1,vlaski2019distributed_part2}.
These schemes study second-order guarantees of  variations of the DGD algorithm (\ref{DGD_intro}). Specifically,  
\cite{vlaski2019distributed_part1,vlaski2019distributed_part2} studied the behavior of (the ATC version  of)  DGD wherein exact gradients are replaced by stochastic approximations; the algorithm is proved to  return approximate second-order stationary points in polynomial number of iterations. Finally, \cite{li2019landscape} proposed a variant of  DGD   to solve the distributed low-rank matrix factorization problem and they prove almost sure convergence to global minima of the problem. 


\vspace{-0.1cm}

\subsection{Paper organization}  
The rest of the paper is organized as follows.  The   main assumptions on the optimization problem and network    are introduced in Sec.~\ref{ProblemSetting_SubSe}. Sec.~\ref{sec:DGD}   studies    guarantees of  DGD    over undirected graphs, along the following steps: i) existing convergence results are discussed in Sec. \ref{sec:DGD_convergence_preliminary}; ii) Sec. \ref{sec:critical-DGD} studies convergence to a neighborhood  of a critical point of $F$; and  iii)    Sec. \ref{sec:SOS-DGD}  establishes second-order guarantees.
\algname/  algorithms are studied in   Sec.~\ref{sec:DOGT}  along the following steps: i) Sub-sequence convergence is proved in  {Sec.}~\ref{Subsequence_Conv_Subsec}; ii) Sec.~\ref{Global_Conv_Subsec} establishes global convergence under the {K\L}  property of $F$; and iii) Sec.~\ref{Guarantees_Sec} derives second-order guarantees  over undirected and directed graphs.   Finally, Sec.~\ref{numerical_experiments} presents some   numerical results.  \vspace{-0.1cm}

\subsection{Notation}\label{sec:notation}
The set of nonnegative integers is denoted by $\mathbb{N}_+$ and we use  $[n]$ as a   shorthand for   $\{1,2,\ldots,n\}$. All vectors are denoted by bold letters and assumed to be column vectors; given a vector $\bx$, 
$||\mathbf{x}||$ denotes the  $\ell_2$ norm of $\mathbf{x}$;  any other specific vector norm is subscripted accordingly.
$\mathbf{x}$ is called \emph{stochastic} if all its components are  nonnegative and sum to one; and $\mathbf{1}$ is the vector of all ones (we write $\mathbf{1}_m$ for the $m$--dimensional vector, if the dimension is not clear from the context).  Given  sets $\mathcal{X},\mathcal{Y}\subseteq \mathbb{R}^m$, we denote $\mathcal{X}\setminus \mathcal{Y}\triangleq \{x\in\mathcal{X}:x\notin\mathcal{Y}\}$, $\overline{\mathcal{X}}\triangleq \mathbb{R}^m\setminus\mathcal{X}$ (complement of $\mathcal{X}$),  and   $\mathbf{x}+\mathcal{X}=\{\mathbf{x}+\mathbf{z}:\mathbf{z}\in\mathcal{X}\}$.   {  $\mathcal{V}_{\mathbf{x}}$ and $\mathcal{B}(\mathbf{x},r)^d$ denote a neighborhood of $\mathbf{x}$ and the $d$-dimensional closed ball of radius $r>0$ centered at $\mathbf{x}$, respectively; when the ball is centered at $\mathbf{0}$, we will write $\mathcal{B}^d_r$. We further define an \emph{annulus} by $\mathcal{S}_{r,\epsilon}\triangleq \mathcal{B}^d_r\setminus\mathcal{B}^d_{r-\epsilon}$, with some $r>\epsilon>0$.}
The Euclidean projection of $\mathbf{x}\in\mathbb{R}^m$  onto the convex closed set  $\mathcal{X}\subseteq \mathbb{R}^m$ is  $\mathrm{proj}_\mathcal{X}(\mathbf{x})\triangleq \arg\min_{\mathbf{y}\in \mathcal{X}}||\mathbf{x}-\mathbf{y}||$. The sublevel set of a function $U$ at $u$ is denoted by $\mathcal{L}_{U}(u)\triangleq \{\mathbf{x}:U(\mathbf{x})\leq u\}$.

Matrices are denoted by capital bold letters;   $A_{ij}$ is the the    $(i,j)$-th element of $\mathbf{A}$; $\mathcal{M}_m(\mathbb{R})$ is the set of all $m\times m$ real   matrices;   $\mathbf{I}$ is the identity matrix (if the dimension is not clear from the context, we write $\mathbf{I}_m$ for the $m\times m$ identity matrix);   $\mathbf{A}\geq 0$ denotes a nonnegative matrix; and    $\mathbf{A}\geq \mathbf{B}$ stands for   $\mathbf{A}-\mathbf{B}\geq 0$.  
The spectrum of a square real matrix $\mathbf{M}$ is denoted by $\text{spec}(\mathbf{M})$ and  its spectral radius is $\text{spradii}(\mathbf{M})\triangleq\max\{|\lambda|:\lambda\in\text{spec}(\mathbf{M})\}$; the spectral norm is $||\mathbf{M}||\triangleq \max_{||\mathbf{x}||\neq 0}||\mathbf{Mx}||/||\mathbf{x}||$, and any other matrix norm is subscripted accordingly. 
Finally, the   minimum (resp. maximum) singular value  are denoted  by   $\sigma_\mathrm{min}(\mathbf{M})$ (resp. $\sigma_\mathrm{max}(\mathbf{M})$) and minimum (resp. maximum) eigenvalue by $\lambda_\mathrm{min}(\mathbf{M})$ (resp. $\lambda_\mathrm{max}(\mathbf{M})$).

 {  The sequence generated by   DGD (and DOGT) depends on the step-size $\alpha$ and the initialization $\bx^0$. When necessary, 
 we  write $\{\bx^\nu(\alpha,\bx^0)\}$ for $\{\bx^\nu\}$.  
 
  Throughout the paper, we assume that all the probability measures are absolutely continuous with respect to the Lebesgue measure.}

\section{Problem \& network  setting}
\label{ProblemSetting_SubSe}\vspace{-0.1cm}
In this section, we introduce the various assumptions on the functions $f_i$ and the graph, under which our results are derived. 

\begin{assumption}[{On Problem \ref{eq:P}}]\label{P_assumption} Given  Problem~\eqref{eq:P}, 
	\begin{enumerate}[leftmargin=1.5cm,label=(\roman*)]
		\item $f_i$ ($\forall i$) is $r+1$ times continuously differentiable for some  $r\geq 1$,   and $\nabla f_i$ is $L_{i}$-Lipschitz continuous. Denote $L_\mathrm{max}\triangleq \max_i~L_i$;
		\item $F$ is 
		coercive.
	\end{enumerate}
\end{assumption}
For some convergence results of DGD we need the following slightly stronger condition.  

\hypertarget{P_assumption_prime}{\emph{Assumption} 2.1'.}   Assumption \ref{P_assumption}-(i) is satisfied  and (ii) each $f_i$ is coercive.\smallskip

We also make the blanket  assumption   that each agent $i$ knows only its  own $f_i$ but not the rest of the objective function.	

Note that Assumption \ref{P_assumption}, particularly  the global Lipschitz gradient continuity   of $f_i$, is quite standard in the literature. Motivated by some applications of interest (see examples below), we will also prove convergence of DGD   under $LC^1$ only and the mild  condition \eqref{Strip_assumption_eq} below (cf. Assumption \ref{Strip_assumption}).    Although   strictly not necessary, coercivity in Assumptions \ref{P_assumption} \& \hyperlink{P_assumption_prime}{2.1'} simplifies  some of our derivations; our results can be extended under the weaker assumption that  \eqref{eq:P} has a   solution.

{Some of the convergence results of DGD and DOGT are established under the assumption that $F$ satisfies the Kurdyka-\L ojasiewicz ({K\L}) inequality \cite{Kurdyka1998, MR0160856}.
\begin{definition}[{K\L} property]\label{KL_def}
	Given a function $U:\mathbb{R}^N\rightarrow\mathbb{R}\cup\{+\infty\}$, we set   $[a<U<b]\triangleq \{\mathbf{z}\in \mathbb{R}^N\,:\, a<U( \mathbf{z})<b\}$, and 
	\begin{enumerate}[label=(\alph*)]
		\item The function $U$ has {K\L} property at $\acute{\mathbf{z}}\in\mathrm{dom}~\partial U$ if there exists $\eta\in(0,+\infty]$, a neighborhood $\mathcal{V}_{\acute{\mathbf{z}}}$, and a continuous concave function $\phi:[0,\eta)\rightarrow \mathbb{R}_+$ such that:
		\begin{enumerate}[label=(\roman*)]
			\item $\phi(0)=0$,
			\item $\phi$ is $\mathcal{C}^1$ on $(0,\eta)$,
			\item for all $s\in(0,\eta)$, $\phi'(s)>0$,
			\item for all $\mathbf{z}\in\mathcal{V}_{\acute{\mathbf{z}}}\cap[U(\acute{\mathbf{z}})<U<U(\acute{\mathbf{z}})+\eta]$, the {K\L} inequality holds:\vspace{-0.1cm}
			\begin{equation}
			\phi'\left(U(\mathbf{z})-U(\acute{\mathbf{z}})\right)\mathrm{dist}(0,\partial U(\mathbf{z}))\geq 1.
			\label{KL_ineq}\vspace{-0.1cm}
			\end{equation}
		\end{enumerate}
		\item A proper lower-semicontinuous function $U$ is called {K\L} if it satisfies the {K\L} inequality at every point in $\mathrm{dom}~\partial U$.
	\end{enumerate}
\end{definition} 

Many problems involve functions satisfying the {K\L}   inequality; real semi-algebraic functions provide a very rich class of functions satisfying the K\L, see   \cite{attouch2010proximal} for a thorough discussion.
}

Second-order guarantees of DGD are obtained under the following  two extra assumptions below; Assumption \ref{Lipschitz_Hessian} is quite standard and widely used in the literature to establish second-order guarantees  of centralized algorithms (e.g.,  \cite{pmlr-v40-Ge15, pmlr-v70-jin17a, Nesterov2006,  Cartis2011_part1, Cartis2011_part2, Agarwal:2017:FAL:3055399.3055464,Curtis2017})  as well as of  distributed algorithms  \cite{Mingyi_SOAlgs_AvoidsSaddles,vlaski2019distributed_part1,vlaski2019distributed_part2}. Assumption  \ref{Strip_assumption} is introduced for this paper and commented below. 

\begin{assumption}
\label{Lipschitz_Hessian}
Each $f_i:\mathbb{R}^m\rightarrow\mathbb{R}$ is twice differentiable and $\nabla^2 f_i$ is $L_{\nabla_i^2}$-Lipschitz continuous.
The Lipschitz constant of $\nabla^2 F$ and $\nabla^2 F_c$ are   $L_{\nabla^2}=\sum_{i=1}^n L_{\nabla^2_i}$ and $L_{\nabla^2_c}= \max_{i}~L_{\nabla^2_i}$, respectively,  {where $F_c$ is defined in eq. \eqref{LyapunovFunc_intro}.}
\end{assumption}



{

\begin{assumption}
\label{Strip_assumption}
(i) Each $f_i$ is $LC^1$; and (ii) there exist $0<\epsilon<R$ and $\delta>0$  such that \vspace{-0.2cm}
\begin{equation}
\label{Strip_assumption_eq}
\inf_{\bt\in\mathcal{S}_{R,\epsilon}}\left\langle\nabla f_i(\bt),\bt/\norm{\bt}\right\rangle\geq \delta,\quad\forall i\in[n].
\end{equation}
\end{assumption}
Roughly speaking, the condition above postulates that the gradient   $\nabla f_i(\bt)$ is positively correlated with any radial direction $\bt/\norm{\bt}$, for all $\bt$ in the annulus $\mathcal{S}_{R,\epsilon}$. A slightly  more restrictive form of the above assumption has appeared in   \cite[Assumption A3]{Gelfand_et_al_1991}. 
Many functions of practical interest satisfy this assumption; some examples arising from machine learning applications are listed below.   
\begin{description}
	\item[Distributed PCA \cite{Eckart1936}:] Given  matrices $\mathbf{M}_i\in\mathbb{R}^{m\times m}$, $i\in [n]$,   the distributed PCA problem is to find the leading eigenvector of   $\sum_{i=1}^n\mathbf{M}_i$ by solving 
	\begin{equation}
	\label{Dist_PCA}
	\min_{\bt\in\mathbb{R}^{m}}\quad \frac{1}{4}\Big\|\boldsymbol{\theta\theta}^\top-\sum_{i=1}^n\mathbf{M}_i\Big\|_F^2,
	\end{equation}
	which can be rewritten in the form  \eqref{eq:P};
	\item \textbf{Phase retrieval { \cite{chi2019nonconvex}}}: Let $\{(\mathbf{a}_i,y_i)\}_{i=1}^n$, with  $\mathbf{a}_i\in\mathbb{R}^{m}$ and $y_i\in\mathbb{R}$ such that  
	$y_i=\mathbf{a}_i^\top\mathbf{M}^\ast\mathbf{a}_i=(\mathbf{a}_i^\top\bt^\ast)^2$, and  $\mathbf{M}^\ast=\bt^{\ast}\bt^{\ast\top}\in\mathbb{R}^{m\times m}$. 
	The phase retrieval  problem reads 
	\begin{equation}\label{PR_problem}
	\min_{\bt\in\mathbb{R}^{m}}\quad \frac{1}{4}\sum_{i=1}^n\left(||\mathbf{a}_i^\top\bt||^2-y_i\right)^2+\frac{\lambda}{2}||\bt||^2,
	\end{equation}
	where $\lambda>0$ is a given parameter.
	\item \textbf{Matrix sensing {\cite{chi2019nonconvex}}: } Let $\{(\mathbf{A}_i,y_i)\}_{i=1}^n$,  with $\mathbf{A}_i\in\mathbb{R}^{m\times m}$ and $y_i\in\mathbb{R}$ such that  $y_i=\langle \mathbf{A}_i,\mathbf{M}^\ast\rangle$, and  $\mathbf{M}^\ast=\bT^{\ast}\bT^{\ast\top}\in\mathbb{R}^{m\times m}$,  $\bT^\ast\in\mathbb{R}^{m\times r}$. The matrix sensing problem reads  \vspace{-0.2cm}
	\begin{equation}\label{MS_problem}
	\min_{\bT\in\mathbb{R}^{m\times r}}\quad \frac{1}{4}\sum_{i=1}^n\Big(\Big\langle \mathbf{A}_i,\bT\bT^\top\Big\rangle -y_i\Big)^2+\frac{\lambda}{2}\norm{\bT}^2_F,
	\end{equation}
	where $\lambda>0$ is a given parameter.
	\item \textbf{Gaussian mixture model {\cite{li2019landscape}}}: 
	Let $\{\mathbf{z}_i\}_{i=1}^n$ be $n$ points drawn from a mixture of $q$ Gaussian distributions, i.e., $\mathbf{z}_i\sim \sum_{d=1}^q\mathcal{N}(\boldsymbol{\mu}_d^\ast,\boldsymbol{\Sigma})$, where   $\mathcal{N}(\boldsymbol{\mu}_d^\ast,\boldsymbol{\Sigma})$ is the Gaussian distribution with mean $\boldsymbol{\mu}_d^\ast\in\mathbb{R}^{m}$ and covariance $\boldsymbol{\Sigma}\in\mathbb{R}^{m\times m}$. The goal is to estimate the mean values $\boldsymbol{\mu}_1^\ast, \ldots, \boldsymbol{\mu}_q^\ast$ by solving the maximum likelihood problem \vspace{-0.2cm}
	\begin{equation}
	\min_{\{\bt_d\in \mathbb{R}^{m}\}_{d=1}^q}\, -\sum_{i=1}^n\log{\bigg(\sum_{d=1}^q\phi_m(\mathbf{z}_i-\bt_d)\bigg)}+\frac{\lambda}{2}\norm{\bt_d}^2,
	\end{equation}
	where $\phi_m(\bt)$ is the multivariate normal distribution with   $\mathbf{0}$ mean and  covariance ${\boldsymbol{\Sigma}}$;
	\item \textbf{Bilinear logistic regression} \cite{Dyrholm_JMLR_2007}: The description of the problem along with some numerical results can be found in  Sec. \ref{BLR_simulations};
	\item \textbf{Artificial neuron { \cite{auer1996exponentially,zhao2010convex}}}: Let $\{(\mathbf{s}_i,\xi_i)\}_{i=1}^n$ be $n$ samples, with $\mathbf{s}_i\in\mathbb{R}^m$, $\xi_i\in\mathbb{R}$, and measurement model $\xi_i=\sigma (\mathbf{s}_i^\top\bt^\ast)$, where $\bt^\ast$ is the optimal weights and $\sigma(\cdot)$ is a \emph{transfer} function; e.g., the logistic regression function $\sigma(\theta)=1/(1+\exp(-\theta))$. The goal is to estimate $\bt^\ast$ by solving
	\begin{equation}
	\min_{\bt\in\mathbb{R}^m}\quad  \sum_{i=1}^n  ~\frac{1}{2n}\left[\left(\xi_i-\sigma(\mathbf{s}_i^\top\bt)\right)^2+\frac{\lambda}{2}\norm{\bt}^2\right],
	\end{equation}
	where $\lambda>0$ is a given parameter. Further binary classification models    satisfying Assumption \ref{Strip_assumption}  include $f_i$ functions such as \cite{zhao2010convex}
	\vspace{-0.2cm}\begin{equation}
	\begin{aligned}
	f_i(\bt)=&1-\tanh{\xi_i\mathbf{s}_i^\top\bt}+\frac{\lambda}{2}\norm{\bt}^2,
	\\
	f_i(\bt)=&\left(1-\sigma(\xi_i\mathbf{s}_i^\top\bt)\right)^2+\frac{\lambda}{2}\norm{\bt}^2,
	\\
	f_i(\bt)=&-\ln{\sigma(\xi_i\mathbf{s}_i^\top\bt)}+\ln{\sigma(\xi_i\mathbf{s}_i^\top\bt+\mu)}+\frac{\lambda}{2}\norm{\bt}^2,
	\end{aligned}
	\end{equation}
	where $\lambda>0$ and $\mu>0$ are given parameters.
\end{description}
In all these examples, Assumption \ref{Strip_assumption} is satisfied for any sufficiently large $R$ and $R-\epsilon$; the proof can be found in   Appendix \ref{Strip_assumption_appendix}.
Note that many of the functions listed above are not    $L$-smooth  on their  entire domain, violating thus (part of) Assumption \ref{P_assumption}(i). Motivated by these examples, we will extend existing convergence results of DGD,   replacing Assumption \ref{P_assumption}(i) with  Assumption \ref{Strip_assumption}.

} 
			
\textbf{Network model:} 
The network   is modeled as a (possibly) directed graph $\mathcal{G}=(\mathcal{V},\mathcal{E})$, where the set of vertices $\mathcal{V}=[n]$ coincides with  the  set  of agents, and the set of edges $\mathcal{E}$ represents the agents' communication links:    $(i,j)\in\mathcal{E}$  if and only if there is link directed from agent $i$ to agent $j$.   The  {in-neighborhood}  of agent $i$ is defined as  $\mathcal{N}_i^{\rm in}=\{j|(j,i)\in\mathcal{E}\}\cup\{i\}$   and represents the set of agents that can send information to agent $i$ (including agent $i$ itself, for notational simplicity). The  {out-neighborhood} of agent $i$ is similarly defined $\mathcal{N}_i^{\rm out}=\{j|(i,j)\in\mathcal{E}\}\cup\{i\}$. When the graph is undirected, these two sets coincide and we use  
$\mathcal{N}_i$ to denote the neighborhood of agent $i$ (with a slight abuse of notation, we use the same symbol $\mathcal{G}$ to denote either directed or undirected graphs). Given a nonnegative matrix $\mathbf{A}\in\mathcal{M}_n(\mathbb{R})$, the directed graph  induced by   $\mathbf{A}$ is defined as $\mathcal{G}_{A}=(\mathcal{V}_A,\mathcal{E}_A)$, where $\mathcal{V}_A\triangleq [n]$ and $(j,i)\in \mathcal{E}_A$ if and only if   $A_{ij}>0$. The set of roots of all the directed spanning trees in $\mathcal{G}_{A}$ is denoted by $\mathcal{R}_A$.
We make the following blanket  standard  assumptions on $\mathcal{G}$.

\begin{assumption}[On the network] \label{Net-Assump}  The   graph  (resp. digraph)  $\mathcal{G}$ is   connected (resp. strongly connected). 
\end{assumption}

\section{The DGD algorithm}\label{sec:DGD} 
{ Consider Problem \eqref{eq:P} and assume that the network is modeled as an undirected graph $\mathcal{G}$. As described in Sec. \ref{Sec_intro}, the DGD algorithm is based on a   decentralization of GD as described  in  \eqref{DGD_intro}.
	It is convenient to rewrite the update \eqref{DGD_intro}  in the matrix/vector form: Using  the definition of \emph{aggregate} function $F_c(\mathbf{x})$ [cf. \eqref{LyapunovFunc_intro}]     and  
	$\bx^\nu\triangleq [\bx_1^{\nu\top},\ldots \bx_n^{\nu\top}]^\top$,
	we have 
	\begin{equation}
	\mathbf{x}^{\nu+1}=\mathbf{W}_D\,\mathbf{x}^{\nu}-\alpha\nabla F_c(\mathbf{x}^\nu),
	\label{DGD_update}
	\end{equation}
	given  $\bx^0\in\mathbb{R}^{mn}$, where 
	$\mathbf{W}_D\triangleq \mathbf{D}\otimes \mathbf{I}_m$, and $\mathbf{D}\in\mathcal{M}_n(\mathbb{R})$ satisfying 
	the following assumption.}
\begin{assumption}
	\label{DoublyStochastic_D}
	$\mathbf{D}\in\mathcal{M}_n(\mathbb{R})$ is nonnegative,  doubly-stochastic, and compliant to $\mathcal{G}$, i.e.,   $D_{ij}>0$ if and only if $(j,i)\in\mathcal{E}$, and $D_{ij}=0$  otherwise. 
\end{assumption} 

\subsection{Existing convergence results}\label{sec:DGD_convergence_preliminary} 
Convergence of  DGD  applied to the nonconvex problem \eqref{eq:P} has been established \cite{WYin_DGD_SIAM2016,WYin_ncvxDGD_SIAM2018}, and   summarized below.
\begin{theorem}[\cite{WYin_DGD_SIAM2016,WYin_ncvxDGD_SIAM2018}]
	\label{DGD_conv_cons_thm}  
	Let Assumptions  \hyperlink{P_assumption_prime}{2.1'}, \ref{Net-Assump} hold. Given arbitrary $\mathbf{x}^0\in\mathbb{R}^{mn}$ and $0<\alpha< \alpha_{\max}\triangleq \sigma_\mathrm{min}(\mathbf{I}+\mathbf{D})/L_c$, let  $\{\mathbf{x}^\nu\}$ be the sequence generated by the DGD algorithm \eqref{DGD_update} under Assumption \ref{DoublyStochastic_D}. Then $\{\mathbf{x}^\nu\}$ is bounded and   
	\begin{enumerate}[label=(\roman*)]
		\item[(i)] \emph{[almost consensus]:} for all  $i\in[n]$ and   $\nu\in \mathbb{N}_+$,\vspace{-0.2cm}
		\begin{equation*}
		\|\mathbf{x}_{i}^\nu-\bar{\mathbf{x}}^\nu\|\leq (\sigma_2)^\nu||\mathbf{x}_i^0|| +\frac{\alpha H}{1-\sigma_2},  \vspace{-0.2cm} 
		\end{equation*}  
		where $\sigma_2<1$ is the second largest singular value of $\mathbf{D}$, and   $H$ is   a universal upper-bound of   $\{||\nabla F_c(\mathbf{x}^\nu)||\}$;
		\item[(ii)] \emph{[stationarity]:}  every limit point $\mathbf{x}^\infty$ of  $\{\mathbf{x}^\nu\}$ is such that  $\mathbf{x}^\infty\in\mathrm{crit}~L_{\alpha}$. 		
	\end{enumerate}
	
	In addition, if	 $L_\alpha$ is a {K\L}  function, then $\{\mathbf{x}^\nu\}$ is globally convergent to some  $\mathbf{x}^\infty\in\mathrm{crit}~L_{\alpha}$.
\end{theorem} 
}

Although $L$-smoothness of $f_i$'s is a common assumption in the literature, above convergence results can also be established without this condition but  under Assumption \ref{Strip_assumption}--see Remark \ref{local_lip_extension_remark} and Appendix  \ref{DGD_conv_cons_thm_extended_proof_appendix} for details.

{Since  \eqref{DGD_update} is the gradient update applied to $L_\alpha$ (cf. \eqref{DGD_update_GDupdate_intro}), non-convergence of the DGD algorithm to strict saddle points of $L_\alpha$ can be established by applying 
\cite[Corollary 2]{Jordan_FOM_AvoidsSaddles} to \eqref{DGD_update_GDupdate_intro}; the statement is given in Theorem \ref{DGD_naive_SOG} below. The following extra assumption on the weight matrix $\mathbf{D}$ is needed.}

\begin{assumption}
\label{nonsingular_W}
The matrix $\mathbf{D}\in\mathcal{M}_n(\mathbb{R})$ is nonsingular.	
\end{assumption}
\vspace{-0.2cm}

{ 
\begin{theorem}
	\label{DGD_naive_SOG}
	Consider Problem \eqref{eq:P}, under Assumptions   \hyperlink{P_assumption_prime}{2.1'}, \ref{Net-Assump}, and further assume that  each $f_i$ is a {K\L}  function.
	Let $\{\mathbf{x}^\nu\}$  be the sequence generated by the DGD algorithm with step-size    $0<\alpha<\frac{\sigma_\mathrm{min}(\mathbf{D})}{L_c}$ and weight matrix $\mathbf{D}$ satisfying  Assumptions \ref{DoublyStochastic_D} and \ref{nonsingular_W}. Then, the stable set of strict saddles has measure zero. Therefore,   $\{\mathbf{x}^\nu\}$  convergences almost surely  to a SoS solution of $L_\alpha$, where the probability is taken over the random initialization $\mathbf{x}^0\in\mathbb{R}^{mn}$.
\end{theorem}
}

As anticipated in Sec. \ref{sec:DGD_intro_contribution}, the above second-order guarantees    are not  satisfactory as they do not provide any information on the behavior of   DGD   near critical points of $F$, including the strict saddles  of $F$. In the following, we fill this gap. We first show that the DGD algorithm convergences to  
neighborhood of the critical points of $F$,  whose size is controlled by the step-size $\alpha>0$ (cf. Section \ref{sec:critical-DGD}). Then, we prove that, for sufficiently small $\alpha>0$,   such critical points are almost surely SoS solutions of \eqref{eq:P}, where the randomization is taken on the initial point (cf. Section \ref{sec:SOS-DGD}).

\subsection{DGD  converges to a neighborhood of  critical points of $F$} \label{sec:critical-DGD}
Let us begin with introducing the definition of $\epsilon$-critical points of $F$.

\begin{definition}
	\label{eps_critical_points_def}
	A point $\boldsymbol{\theta}\in\mathbb{R}^m$  such that 
	$||\nabla F(\boldsymbol{\theta})||\leq \varepsilon$,   with $\varepsilon> 0$, is called  $\varepsilon$-critical point of $F$. The set of $\varepsilon$-critical points of $F$ is denoted by $\mathrm{crit}_\varepsilon F$.
\end{definition} 

In this section, we prove  that when the step-size is sufficiently small and DGD is initialized in a compact set, the  iterates   $\{\bx_i^\nu\}$, $i\in[n]$, converge to an arbitrarily small neighborhood of critical points of $F$--the result is formally stated in  Theorem \ref{prop_DGD_local_convergence}.  Roughly speaking, this is proved chaining  the following  intermediate results: \begin{description}
	\item[i) Lemma \ref{L_criticalPoints_are_alpha_statPoints}:] Every limit point of DGD is an $\mathcal{O}(\alpha)$-critical point of $F$;
	\item[ii) Lemma \ref{Y_set_found}:] Every sequence generated by DGD for given $\alpha>0$ and    initialization in a compact set, is enclosed in some compact set, for all  $\alpha\downarrow 0$; and \item[iii) Lemma \ref{approximate_critical_points_lemma}:] Any $\epsilon$-critical point of $F$  {achievable} by DGD 
	is arbitrarily close to a critical point of $F$, when $\epsilon$ is sufficiently small. 
\end{description}  
Lemma \ref{L_criticalPoints_are_alpha_statPoints}  implies that, for any given $\epsilon>0$,  one can find arbitrarily   small $\alpha>0$ so that every limit point of each $\{\bx_i^\nu\}$  (whose existence is guaranteed by Lemma \ref{Y_set_found}) is an $\epsilon$-critical point of $F$. Finally, Lemma  \ref{approximate_critical_points_lemma} guarantees that every such $\epsilon$-critical point can be made  arbitrarily close to a critical point of $F$ as $\epsilon \downarrow 0$. The proof of the above three lemmata follows.

%
%

\begin{lemma}
	\label{L_criticalPoints_are_alpha_statPoints}
	Let Assumptions  \hyperlink{P_assumption_prime}{2.1'} and \ref{Net-Assump} hold. Given arbitrary $\mathbf{x}^0\in\mathbb{R}^{mn}$ and $0<\alpha<\sigma_\mathrm{min}(\mathbf{I}+\mathbf{D})/L_c$, 
	every limit point  
	$\mathbf{x}^\infty=[\mathbf{x}_1^{\infty\top},\ldots,\mathbf{x}_n^{\infty\top}]^\top$
	of   $\{\mathbf{x}^\nu\}$ generated  by the DGD algorithm satisfies $\bar{\mathbf{x}}^\infty\in \mathrm{crit}_{K'\alpha}F$, with $\bar{\mathbf{x}}^\infty\triangleq (1/n)\,\sum_{i=1}^n\mathbf{x}^\infty_i$ and $K'= n\sqrt{n}L_cH/$ $(1-\sigma_2)$, where $H$ and $\sigma_2$ are defined in Theorem \ref{DGD_conv_cons_thm}.
\end{lemma}
\begin{proof}
	By Theorem \ref{DGD_conv_cons_thm}(ii),     $({\mathbf{1}\otimes\mathbf{I}})^\top \nabla L_\alpha(\mathbf{x}^\infty)=\mathbf{0}$, which using (\ref{LyapunovFunc_intro}) and the column stochasticity of $\mathbf{D}$ yields  
	$({\mathbf{1}\otimes\mathbf{I}})^\top$ $\nabla F_c(\mathbf{x}^\infty)= \mathbf{0}$. Hence,
	\begin{equation}
	\label{nablaD_epsCriticalbound}
	\begin{aligned}
	\norm{\nabla F(\bar{\mathbf{x}}^\infty)}= & \norm{({\mathbf{1}\otimes\mathbf{I}})^\top \left(\nabla F_c(\mathbf{1}\otimes\bar{\mathbf{x}}^\infty)- \nabla F_c(\mathbf{x}^\infty)\right)}
	\\
	\leq & L_c\sqrt{n}\norm{\mathbf{x}^\infty-\mathbf{1}\otimes\bar{\mathbf{x}}^\infty}\overset{(a)}{\leq} \alpha\cdot \frac{n\sqrt{n}L_cH}{1-\sigma_2}, 
	\end{aligned}\vspace{-0.3cm}
	\end{equation}
	where 
	in  (a) we used  Theorem \ref{DGD_conv_cons_thm}(i).
\end{proof}

To proceed,   we limit DGD initialization to $\mathbf{x}_i^0\in\mathcal{X}_i$,  $i\in[n],$ where $\mathcal{X}_i^0\subseteq \mathbb{R}^{m}$ is some compact set with positive Lebesgue measure.

\begin{lemma}
	\label{Y_set_found}
	Consider Problem \eqref{eq:P}, under Assumptions  \hyperlink{P_assumption_prime}{2.1'}, \ref{Strip_assumption} and \ref{Net-Assump}.
	Let $\{\mathbf{x}^\nu(\alpha,\bx^0)\}$ be any sequence generated by DGD under Assumption \ref{DoublyStochastic_D}, with step-size $\alpha$ and initialization  $\mathbf{x}^0$. Then, there exists a bounded set $\mathcal{Y}$ such that  $\{\mathbf{x}^\nu(\alpha,\bx^0)\}\subseteq \mathcal{Y}$, for all $0<{\alpha}\leq\alpha_{\max}= \sigma_\mathrm{min}(\mathbf{I}+\mathbf{D})/L_c$ and $\mathbf{x}_i^0\in\mathcal{X}_i\subseteq\mathcal{B}^m_R,  i\in[n]$, where $R$ is defined in Assumption \ref{Strip_assumption}. 
\end{lemma}

\begin{proof}
We proceed by induction. For the sake of  notation,  throughout the proof, we will use for $\bx^\nu(\alpha,\bx^0)$ 
the shorthand   $\bx^\nu$. 
Define $h\triangleq \max_{i\in[n],\boldsymbol{\theta}\in\mathcal{B}^m_R}||\nabla f_i(\boldsymbol{\theta})||$. By assumption, there holds $\mathbf{x}_i^0\in\mathcal{B}^m_R$, for all $i$. Suppose $\mathbf{x}_i^\nu\in\mathcal{B}^m_R$, for all $i$. If $\mathbf{x}_i^{\nu}\in\mathcal{B}^m_{R-\epsilon}$ and $\alpha\leq \epsilon D_{ii}/h$, then $\mathbf{x}_i^\nu-\frac{\alpha}{D_{ii}}\nabla f_i(\mathbf{x}_i^\nu)\in\mathcal{B}^m_R$, since \vspace{-0.2cm}
\begin{equation}
\big\|\mathbf{x}_i^\nu-\frac{\alpha}{D_{ii}}\nabla f_i(\mathbf{x}_i^\nu)\big\|\leq  \norm{\mathbf{x}_i^\nu}+\frac{\alpha}{D_{ii}} \norm{\nabla f_i(\mathbf{x}_i^\nu)}
\leq  R-\epsilon+\frac{\alpha h}{D_{ii}}.
\end{equation}

If $\mathbf{x}_i^{\nu}\in\mathcal{S}_{R,\epsilon}$ and $\alpha\leq 2D_{ii}\delta(R-\epsilon)/h^2$, then $\mathbf{x}_i^\nu-\frac{\alpha}{D_{ii}}\nabla f_i(\mathbf{x}_i^\nu)\in\mathcal{B}^m_R$, since
\begin{equation}
\begin{aligned}
  \big\|\mathbf{x}_i^\nu-\frac{\alpha}{D_{ii}}\nabla f_i(\mathbf{x}_i^\nu)\big\|^2 & = 
       \norm{\mathbf{x}_i^\nu}^2-\frac{2\alpha||\mathbf{x}_i^\nu||}{D_{ii}}\left\langle \frac{\mathbf{x}_i^\nu}{||\mathbf{x}_i^\nu||},\nabla f_i(\mathbf{x}_i^\nu)\right\rangle+\frac{\alpha^2}{D^2_{ii}} \norm{\nabla f_i(\mathbf{x}_i^\nu)}^2\\ & \leq R^2-\frac{2\alpha \delta(R-\epsilon)}{D_{ii}}+\frac{\alpha^2 h^2}{D^2_{ii}}.
\end{aligned}
\end{equation}

By agents' updates   $\mathbf{x}_i^{\nu+1}=\sum_{j\neq i}D_{ij}\mathbf{x}_j^\nu+D_{ii}(\mathbf{x}_i^\nu-\frac{\alpha}{D_{ii}}\nabla f_i(\mathbf{x}_i^\nu))$ and convexity of the norm, we conclude that  if  $\mathbf{x}_i^{\nu}\in\mathcal{B}^m_{R}$,  for all $i$, and  $0<\alpha\leq \alpha_b\triangleq \min_i\min\{\epsilon D_{ii}/h,2D_{ii}\delta(R-\epsilon)/h^2\}$, then $\mathbf{x}_i^{\nu+1}\in\mathcal{B}^m_R$. This proves that, for $\alpha\in(0,\alpha_b]$, any sequence $\{\bx^{\nu}_i\}$ initialized in $\mathcal{B}^m_R$ lies in $\mathcal{B}^m_R$, for all $i$.

We prove now the same result for  $\alpha\in[\alpha_b,\sigma_\mathrm{min}(\mathbf{I}+\mathbf{D})/L_c]$. Note that since each $f_i$ is coercive (cf. Assumption \hyperlink{P_assumption_prime}{2.1'(ii)}), any sublevel set of $L_\alpha$ is compact. Also, since $\{L_\alpha(\mathbf{x}^\nu)\}$ is non-increasing for all $\alpha\in(0,\sigma_\mathrm{min}(\mathbf{I}+\mathbf{D})/L_c]$  (cf. \cite[lemma 2]{WYin_DGD_SIAM2016}), then   $\{\mathbf{x}^\nu\}\subseteq\mathcal{L}_{L_\alpha}(F_c(\mathbf{x}^0)+\frac{1}{2\alpha}||\mathbf{x}^0||^2_{\mathbf{I-W}})$, and furthermore,  
\begin{equation}
\begin{aligned}
& \mathcal{L}_{L_\alpha}\left(F_c(\mathbf{x}^0)+\frac{1}{2\alpha}||\mathbf{x}^0||^2_{\mathbf{I-W}}\right)\subseteq  \mathcal{L}_{L_\alpha}\left(F_c(\mathbf{x}^0)+\frac{1}{2\alpha_b}||\mathbf{x}^0||^2_{\mathbf{I-W}}\right)
\\
\subseteq & \mathcal{L}_{F_c}\left(F_c(\mathbf{x}^0)+\frac{1}{2\alpha_b}||\mathbf{x}^0||^2_{\mathbf{I-W}}\right) \subseteq \mathcal{L}_{F_c} \left(\max_{\mathbf{x}_i^0\in\mathcal{B}^m_R, i\in [n]} \left\{ F_c(\mathbf{x}^0)+\frac{1}{2\alpha_b}||\mathbf{x}^0||^2_{\mathbf{I-W}}\right\}\right).
\end{aligned}
\end{equation}
Since $||\mathbf{x}^0||^2_{\mathbf{I-W}}\leq 2||\mathbf{x}^0||^2$, it follows   
\begin{equation}\label{L_0_Y_inclusion}
 \mathcal{L}_{L_\alpha}\left(F_c(\mathbf{x}^0)+\frac{1}{2\alpha}||\mathbf{x}^0||^2_{\mathbf{I-W}}\right)\subseteq \underbrace{\mathcal{L}_{F_c}\left(\max_{\mathbf{x}_i^0\in\mathcal{B}^m_R,   i\in [n]} \left\{\sum_{i=1}^n f_i(\mathbf{x}_i^0)\right\}+\frac{R^2}{\alpha_b}\right)}_{\triangleq \bar{\mathcal{L}}}.
\end{equation}
The statement of the lemma hods with   $\mathcal{Y}=\bar{\mathcal{L}}\cup\prod_{i=1}^n\mathcal{B}^m_{R}$.  
\end{proof}

The following lemma shows that any $\epsilon$-critical point of $F$ achievable by DGD (i.e., any point in $\mathrm{crit}_\varepsilon F\cap\bar{\mathcal{Y}}$) can be made  arbitrarily close to a critical point of $F$, when $\epsilon>0$ (and thus $\alpha>0$) is sufficiently small. 




\begin{lemma}
\label{approximate_critical_points_lemma}
Suppose $F:\mathbb{R}^{m}\rightarrow\mathbb{R}$ is continuously differentiable. For any given compact set  $\bar{\mathcal{Y}}\subseteq \mathbb{R}^m$, there holds 
\vspace{-0.3cm}
\begin{equation}\label{eq:vanishing_crit_f_error_statement}
\lim_{\varepsilon\rightarrow 0}~\max_{\mathbf{q}\in\mathrm{crit}_\varepsilon F\cap\bar{\mathcal{Y}}}~\mathrm{dist}(\mathbf{q},\mathrm{crit}~F)=0.\vspace{-0.1cm}
\end{equation}
\end{lemma}

\begin{proof}
We prove the lemma  by contradiction. Suppose \vspace{-0.2cm}
\begin{equation}\label{eq:vanishing_crit_f_error}
\limsup_{\varepsilon\rightarrow 0}~\max_{\mathbf{q}\in\mathrm{crit}_\varepsilon F\cap\bar{\mathcal{Y}}}~\mathrm{dist}(\mathbf{q},\mathrm{crit}~F)=\gamma>0.\vspace{-0.2cm}
\end{equation}
Then,  one can construct   $\{\mathbf{q}^\nu\}$ with   $\mathbf{q}^\nu \in\mathrm{crit}_{1/\nu} F\cap\bar{\mathcal{Y}}$ such that $\mathrm{dist}(\mathbf{q}^\nu,\mathrm{crit} F)\geq \gamma$ for all  $\nu\in\mathbb{N}$.  Since  $\nabla F$  is continuous, 
$\mathrm{crit_{1}} F$ is closed and $\mathrm{crit_{1}} F\cap\bar{\mathcal{Y}}$ is compact. 
Note that $\{\mathbf{q}^\nu\}\subseteq\mathrm{crit_{1}} F\cap\bar{\mathcal{Y}}$, which ensures    
$\{\mathbf{q}^\nu\}$ is bounded. Let $\{\mathbf{q}^{t_\nu}\}$ be a convergent subsequence of $\{\mathbf{q}^\nu\}$;  its    limit point $\mathbf{q}^\infty$ satisfies $\mathrm{dist}(\mathbf{q}^\infty,\mathrm{crit}~F)\geq \gamma$. By construction,  for any $\acute{\nu}\in\mathbb{N}$, $\{\mathbf{q}^{t_\nu}\}$ eventually settles  in $\mathrm{crit}_{1/\acute{\nu}}~F\cap\bar{\mathcal{Y}}$, thus $\mathbf{q}^\infty\in\mathrm{crit}_{1/\acute{\nu}}~F\cap\bar{\mathcal{Y}}$. This  means that $||\nabla F(\mathbf{q}^{\acute{\nu}})||\leq 1/\acute{\nu}$, for all $\acute{\nu}\in\mathbb{N}$, implying $||\nabla F(\mathbf{q}^\infty)||=0$. Hence  $\mathrm{dist}(\mathbf{q}^\infty,\mathrm{crit}~F)=0$, which contradicts \eqref{eq:vanishing_crit_f_error}.
	\end{proof}
	
We can now combine Lemmas \ref{L_criticalPoints_are_alpha_statPoints}-\ref{approximate_critical_points_lemma}   with Theorem \ref{DGD_conv_cons_thm}(i) and state the   main result of this section.  




 \begin{theorem}
	\label{prop_DGD_local_convergence}	
	Let Assumptions  \hyperlink{P_assumption_prime}{2.1'}, \ref{Strip_assumption} and \ref{Net-Assump} hold.  Let  $\epsilon>0$. There exists $\bar{\alpha}>0$ (which depends on $\epsilon$)
	such that with any initialization $\mathbf{x}_i^0\in\mathcal{X}_i^0\subseteq \mathcal{B}^m_R$ ($R>0$ is defined in Assumption \ref{Strip_assumption}), $i\in [n]$, and any step-size $0<\alpha\leq \bar{\alpha}$, all the limit points 
	$\mathbf{x}^\infty(\alpha,\bx^0)=[\mathbf{x}_1^{\infty}(\alpha,\bx^0)^\top,\ldots,\mathbf{x}_n^{\infty}(\alpha,\bx^0)^\top]^\top$
	of the sequence $\{\mathbf{x}^\nu(\alpha,\bx^0)\}$, generated by DGD   satisfies
	\begin{equation}\label{eq:crit_neig}
	\mathrm{dist}\big(\bar{\mathbf{x}}^\infty(\alpha,\bx^0),\mathrm{crit\,F}\big)<\epsilon \quad \text{and} \quad \norm{\mathbf{x}^\infty(\alpha,\bx^0)-\mathbf{1}\otimes \bar{\mathbf{x}}^\infty(\alpha,\bx^0)}<\epsilon,\vspace{-0.1cm}
	\end{equation}
	where  $\bar{\mathbf{x}}^\infty(\alpha,\bx^0)\triangleq (1/n)\,\sum_{i=1}^n\mathbf{x}^\infty_i(\alpha,\bx^0)$.
\end{theorem}
\begin{proof}
Combining Lemmata \ref{L_criticalPoints_are_alpha_statPoints}-\ref{approximate_critical_points_lemma}  proves that there exists some $\alpha_1>0$ such that  $\mathrm{dist}(\bar{\mathbf{x}}^\infty(\alpha,\bx^0), \mathrm{crit\,F})<\epsilon$, for all $\alpha\leq \alpha_1$. In addition, Theorem \ref{DGD_conv_cons_thm}(i), with $H=\sup_{\bx\in\mathcal{Y}} F_c(\bx)$, implies that there exists some $\alpha_2>0$   such that $||\mathbf{x}^\infty(\alpha,\bx^0)-\mathbf{1}\otimes \bar{\mathbf{x}}^\infty(\alpha,\bx^0)||<\epsilon$,   for all $\alpha\leq \alpha_2$. Hence, choosing $\bar{\alpha}= \min\{\alpha_1,\alpha_2\}$ proves  \eqref{eq:crit_neig}.
\end{proof}

	\subsection{DGD likely converges to a neighborhood of SoS solutions of $F$} \label{sec:SOS-DGD} 
	We study now second-order guarantees of DGD. 
	Our path to prove almost sure convergence   to a neighborhood of SoS solutions of \eqref{eq:P} will pass through  the non-convergence  of DGD to strict saddles of $L_\alpha$ (cf.   Theorem \ref{L_criticalPoints_are_alpha_statPoints}).  Roughly speaking, our idea is to show that whenever $\bar{\mathbf{x}}^\infty=1/n\,\sum_{i=1}^n\mathbf{x}^\infty_i$ belongs to a sufficiently small neighborhood  of a  strict saddle  of $F$ inside the region (\ref{eq:crit_neig}),   
	$\mathbf{x}^\infty=[\mathbf{x}_1^{\infty\top},\ldots,\mathbf{x}_n^{\infty\top}]^\top$
	 must be a  \emph{strict  saddle of  $L_\alpha$}.
	The escaping properties of DGD from strict saddles of  $L_\alpha$   will then ensure that it is unlikely that $\{\bar{\mathbf{x}}^\nu= 1/n\,\sum_{i=1}^n\mathbf{x}^\nu_i\}$  gets trapped in   a    neighborhood of a strict saddle of $F$, thus ending in a   neighborhood of a SoS solution of \eqref{eq:P}.
	Proposition \ref{DGD_SSP_L_and_F_connect} makes this argument formal; in particular,  conditions (i)-(iii) identify the   neighborhood of a strict saddle of $F$ with the mentioned escaping properties.\vspace{-0.1cm}

\begin{proposition}
\label{DGD_SSP_L_and_F_connect}
Consider the setting of Lemma \ref{Y_set_found} and further assume that   Assumption \ref{Lipschitz_Hessian} hold. Let   $\bar{\mathcal{Y}}$ be the image of the compact set $\mathcal{Y}$ (defined in Lemma \ref{Y_set_found}) through the linear operator $(\mathbf{1}_n\otimes \mathbf{I}_m)^\top$. Suppose that the  limit point 
$\mathbf{x}^\infty=[\mathbf{x}_1^{\infty\top},\ldots,\mathbf{x}_n^{\infty\top}]^\top$
of $\{\mathbf{x}^\nu\}$, along with  $\bar{\mathbf{x}}^\infty= 1/n\sum_{i=1}^n\mathbf{x}_{i}^\infty$, satisfy  
\begin{enumerate}[label=(\roman*)]
	\item $\mathrm{dist}(\bar{\mathbf{x}}^\infty,\mathrm{crit\,F})<\dfrac{\delta}{2L_{\nabla^2}}$,
	\item $\norm{\mathbf{x}^\infty-\mathbf{1}\otimes \bar{\mathbf{x}}^\infty}<\dfrac{\delta}{2nL_{\nabla^2_c}}$,
	\item   There exists   $\boldsymbol{\theta}^\ast\in\mathrm{proj}_{\mathrm{crit\, F}}(\bar{\mathbf{x}}^\infty)\cap \Theta^\ast_{ss}$,
\end{enumerate}
\vspace{.25cm} 
for some  $\delta$ such that $\delta\leq -\lambda_\mathrm{min}\left(\nabla^2 F(\boldsymbol{\theta}^\ast)\right),\forall \boldsymbol{\theta}^\ast\in\Theta_{ss}^\ast\cap\bar{\mathcal{Y}}$.  Then, ${\mathbf{x}^\infty}$ is a strict saddle point of $L_\alpha$.
\end{proposition}
	\begin{proof}
		Given $\boldsymbol{\theta}\in\mathbb{R}^{m}$, let $\boldsymbol{\upsilon}(\boldsymbol{\theta})$ denote  the unitary eigenvector of  $\nabla^2 F(\boldsymbol{\theta})$ associated with the smallest eigenvalue, and define  $\tilde{\boldsymbol{\upsilon}}(\boldsymbol{\theta})\triangleq \mathbf{1}\otimes \boldsymbol{\upsilon}(\boldsymbol{\theta})$. Then,  we have
		\begin{equation}\label{LHessian_bound0}
		\hspace{-0.4cm}\begin{aligned}
		&\tilde{\boldsymbol{\upsilon}}(\boldsymbol{\theta})^\top\nabla^2 L_\alpha (\mathbf{x}^\infty)\tilde{\boldsymbol{\upsilon}}(\boldsymbol{\theta})\overset{(a)}{=} \tilde{\boldsymbol{\upsilon}}(\boldsymbol{\theta})^\top\nabla^2 F_c (\mathbf{x}^\infty)\tilde{\boldsymbol{\upsilon}}(\boldsymbol{\theta}) 	
		\\
		& \qquad {\leq}\, {\boldsymbol{\upsilon}}(\boldsymbol{\theta})^\top\nabla^2 F (\boldsymbol{\theta}){\boldsymbol{\upsilon}}(\boldsymbol{\theta})
		\\
		&\qquad \quad +||\nabla^2 F (\bar{\mathbf{x}}^\infty)-\nabla^2 F (\boldsymbol{\theta})||\norm{\boldsymbol{\upsilon}(\boldsymbol{\theta})}^2
		+||\nabla^2 F_c (\mathbf{x}^\infty)-\nabla^2 F_c (\mathbf{1}\otimes\bar{\mathbf{x}}^\infty)||\norm{\tilde{\boldsymbol{\upsilon}}(\boldsymbol{\theta})}^2
		\\ 
		& \qquad \overset{(b)}{\leq}\, {\boldsymbol{\upsilon}}(\boldsymbol{\theta})^\top\nabla^2 F (\boldsymbol{\theta}){\boldsymbol{\upsilon}}(\boldsymbol{\theta})
		+L_{\nabla^2}\norm{\bar{\mathbf{x}}^\infty-\boldsymbol{\theta}}+n\,L_{\nabla^2_c} \norm{\mathbf{x}^\infty-\mathbf{1}\otimes \bar{\mathbf{x}}^\infty}
		\end{aligned}
		\end{equation}
		where   (a) follows from $\tilde{\boldsymbol{\upsilon}}(\boldsymbol{\theta})\in\mathrm{null}(\mathbf{W}_D-\mathbf{I})$;  
		and  (b) is due to Assumption \ref{Lipschitz_Hessian}. 
		Let us now evaluate \eqref{LHessian_bound0} at some $\boldsymbol{\theta}^\ast$ as defined in condition (iii) of the proposition; using    ${\boldsymbol{\upsilon}}(\boldsymbol{\theta}^\ast)^\top\nabla^2 F (\boldsymbol{\theta}^\ast){\boldsymbol{\upsilon}}(\boldsymbol{\theta}^\ast)\leq -\delta$ and conditions (i) and (ii), yields  
		$\tilde{\boldsymbol{\upsilon}}(\boldsymbol{\theta}^\ast)^\top\nabla^2 L_\alpha (\mathbf{x}^\infty)$ $\tilde{\boldsymbol{\upsilon}}(\boldsymbol{\theta}^\ast)
		<0.$ 
		By the Rayleigh-Ritz theorem, 
		it must be $\lambda_\mathrm{min}(\nabla^2 L_\alpha (\mathbf{x}^\infty))<0$. This, together with $\mathbf{x}^\infty\in\mathrm{crit}~L_{\alpha}$ (cf. Theorem \ref{DGD_conv_cons_thm}(ii)),  proves the proposition. 
	\end{proof}
	
	Invoking now Theorem \ref{prop_DGD_local_convergence}, we infer that there exists a sufficiently small $\alpha>0$ such that conditions (i) and (ii) of Proposition \ref{DGD_SSP_L_and_F_connect} are always satisfied, implying that ${\mathbf{x}^\infty}$ is a strict saddle   of $L_\alpha$ if there exists a strict saddle of $F$ ``close'' to $\bar{\mathbf{x}}^\infty$ [in the sense of (iii)].
	This is formally stated next.
	

\begin{corollary}
	\label{DGD_SSP_L_and_F_connect_2_AMIREM}
	Consider the setting of Theorem \ref{prop_DGD_local_convergence} and  Proposition \ref{DGD_SSP_L_and_F_connect}.  
	There exists a sufficiently small $\alpha>0$ such that, if $\mathrm{proj}_{\mathrm{crit\, F}}(\bar{\mathbf{x}}^\infty)\cap \Theta^\ast_{ss}\neq \emptyset$, then  ${\mathbf{x}^\infty}$ is a strict saddle   of $L_\alpha$. 
\end{corollary}

To state our final result, let us introduce the following merit function: given 
$\mathbf{x}=[\mathbf{x}_1^{\top},\ldots,\mathbf{x}_n^{\top}]^\top$
let\vspace{-0.3cm}
\begin{equation*}
M(\mathbf{x})\triangleq \max\Big(\mathrm{dist}(\bar{\mathbf{x}},\mathcal{X}_{SOS}),\norm{\mathbf{x}-\mathbf{1}\otimes \bar{\mathbf{x}}}\Big),\vspace{-0.1cm}
\end{equation*}
where  $\mathcal{X}_{SOS}$ denotes the set of SoS solutions of \eqref{eq:P}, and $\bar{\mathbf{x}}=  1/n\sum_{i=1}^n \mathbf{x}_i$.   $M(\mathbf{x})$ capture the distance of the    average $\bar{\mathbf{x}}$ from the set of SoS solutions of \eqref{eq:P} and well as the consensus disagreement of the agents' local variables $\bar{\mathbf{x}}_i$. \vspace{-0.1cm}

\begin{theorem}
\label{SOG_of_DGD_thm}
Consider Problem \eqref{eq:P} under   Assumptions \hyperlink{P_assumption_prime}{2.1'},  \ref{Lipschitz_Hessian}, \ref{Strip_assumption},  and \ref{Net-Assump};    further assume that   {each $f_i$ is a {K\L}  function}. For every $\epsilon>0$, there exists sufficiently small  $0<\bar{\alpha}<\frac{\sigma_\mathrm{min}(\mathbf{D})}{L_c}$ such that \vspace{-0.2cm}
\begin{equation*}
\mathbb{P}_{\mathbf{x}^0}\big(M(\mathbf{x}^\infty)\leq \epsilon\big)=1,	 
\end{equation*}
where 
$\mathbf{x}^\infty=[\mathbf{x}_1^{\infty\top},\ldots,\mathbf{x}_n^{\infty\top}]^\top$
is the limit point of the   sequence  $\{\mathbf{x}^\nu\}$    generated by the DGD algorithm \eqref{DGD_update} with  $\alpha\in (0,\bar{\alpha}]$,  the  weight matrix $\mathbf{D}$ satisfying Assumptions \ref{DoublyStochastic_D} and  \ref{nonsingular_W}, and initialization $\mathbf{x}^0\in \prod_{i=1}^n \mathcal{X}^0_i\subseteq \prod_{i=1}^n \mathcal{B}^m_R$;    $R$ is defined in Assumption \ref{Strip_assumption} and each $\mathcal{X}_i^0$ has positive Lebesgue measure; 
and 	the probability is taken over the  initialization $\mathbf{x}^0\in \prod_{i=1}^n \mathcal{X}^0_i$. Furthermore, any $\boldsymbol{\theta}^\ast\in\mathrm{proj}_{\mathrm{crit~F}}(\bar{\mathbf{x}}^\infty)$ is almost surely a SoS solution of $F$ where $\bar{\mathbf{x}}^\infty=(1/n)\sum_{i=1}^n \mathbf{x}^\infty_i$.
\end{theorem}
\begin{proof}
For sufficiently small $\alpha<\bar\alpha_1$, if $\mathrm{proj}_{\mathrm{crit~F}}(\bar{\mathbf{x}}^\infty)$ contains a strict saddle point of $F$, then $\mathbf{x}^\infty$  is also a strict saddle point of $L_\alpha$ (by Corollary \ref{DGD_SSP_L_and_F_connect_2_AMIREM}). Let also $\bar\alpha_2$ be a sufficiently small step-size such that  every limit point $\mathbf{x}^\infty$ satisfies $\mathrm{dist}(\bar{\mathbf{x}}^\infty,\mathrm{crit}~F)\leq \epsilon$ and $\norm{\mathbf{x}^\infty-\mathbf{1}\otimes \bar{\mathbf{x}}^\infty}\leq\epsilon$ (by Theorem~\ref{prop_DGD_local_convergence}).  Now consider DGD update \eqref{DGD_update} with $\alpha<\min\{\bar\alpha_1,\bar\alpha_2\}$ and $\mathbf{x}^0$ being drawn randomly from the set of probability one measure $\prod_{i=1}^n \mathcal{X}^0_i$ for which the algorithm converges to a SoS solution  of $L_\alpha$ (by Theorem \ref{DGD_naive_SOG}\footnote{Note that the conclusion of Theorem \ref{DGD_naive_SOG} is valid also when  the set of initial points is restricted to $\prod_{i=1}^n\mathcal{X}_i^0$, as $\prod_{i=1}^n\mathcal{X}_i^0$ has positive measure (the Cartesian product of sets with positive measure has positive measure-- cf. \cite[Sec. 35]{halmos1976measure}).}). Finally, by the above properties of $\alpha$, it holds that $M(\mathbf{x}^\infty)\leq\epsilon$ and $\mathrm{proj}_{\mathrm{crit~F}}(\bar{\mathbf{x}}^\infty)$ must  contain only  SoS solutions of $F$. Therefore,    there exists a $\boldsymbol{\theta}^\ast\in \mathrm{crit}~F$ such that $\boldsymbol{\theta}^\ast\in\mathcal{X}_{SoS}$ and  $\|\bar{\mathbf{x}}^\infty-\boldsymbol{\theta}^\ast\|\le \epsilon$. 
\end{proof}


\begin{remark}\label{local_lip_extension_remark}
All (first- and second-order) convergence results of DGD established in this section remain valid when   $\nabla f_i$'s are not globally Lipschitz continuous [Assumption \ref{P_assumption}(i)] but Assumption \ref{Strip_assumption} holds. Specifically, Theorems \ref{DGD_conv_cons_thm}, \ref{DGD_naive_SOG}, \ref{prop_DGD_local_convergence} and Lemmata \ref{L_criticalPoints_are_alpha_statPoints}-\ref{Y_set_found} hold if one replaces   Assumption \ref{P_assumption}(i)  with  Assumption \ref{Strip_assumption}   and  the global Lipschitz constant $L_c$ with the  Lipschitz constant of $\nabla F_c$ \emph{restricted} to the compact set $\tilde{\mathcal{Y}}$, defined in   Appendix \ref{DGD_conv_cons_thm_extended_proof_appendix}, where we refer to for the technical details.   
\end{remark}

	\section{DOGT Algorithms} \label{sec:DOGT} 
	
	
	{
		The family of DOGT algorithms is introduced in  Sec. \ref{DOGT_method_intro}. We begin here rewriting \eqref{eq:DOGT_x_update_intro}-\eqref{eq:DOGT_track_update_intro} in matrix/vector form.  Denoting $\bx^\nu\triangleq [\bx_1^{\nu\top},\ldots \bx_n^{\nu\top}]^\top$ and $\by^\nu\triangleq [\by_1^{\nu\top},\ldots \by_n^{\nu\top}]^\top$, we have  
		\begin{equation}\label{NEXT_generic}
		\left\lbrace\begin{aligned}
		\mathbf{x}^{\nu+1} &=\mathbf{W}_R\, \mathbf{x}^{\nu}-\alpha\,\mathbf{y}^\nu,
		\\
		\mathbf{y}^{\nu+1} &=\mathbf{W}_C\, \mathbf{y}^{\nu}+\nabla F_c\big(\mathbf{x}^{\nu+1}\big)-\nabla F_c\big(	\mathbf{x}^\nu\big),
		\end{aligned}\right.
		\end{equation}
		where    $\mathbf{W}_R\triangleq \mathbf{R} \otimes \mathbf{I}_m$ and  $\mathbf{W}_C\triangleq \mathbf{C} \otimes \mathbf{I}_m$ with $\mathbf{R}\triangleq (R_{ij})_{i,j=1}^n$ and $\mathbf{C}\triangleq (C_{ij})_{i,j=1}^n$ being some \emph{column-stochastic} and \emph{row-stochastic} matrices (respectively) compliant to the graph $\mathcal{G}$ (cf. Assumption \ref{matrix_C_R} below).  	
		The initialization of (\ref{NEXT_generic}) is set to   $\mathbf{x}^{0}\in \mathbb{R}^{mn}$ and $\mathbf{y}^0\in\nabla F_c(\mathbf{x}^0)+\text{span}\left(\mathbf{W}_C-\mathbf{I}\right)$. Note that the latter condition   is instrumental to preserve the \emph{total-sum} of the ${y}$-variables, namely $\sum_i\mathbf{y}_i^\nu=\sum_if_i(\bx_i^\nu)$ (which holds due to the column-stochasticity of matrix $\mathbf{C}$--cf. Assumption \ref{matrix_C_R}). This property is imperative for the ${y}$-variables to track the sum-gradient. Notice that the condition used in the literature \cite{NEXT,Harnessing_Na_Li, SONATA_technical_report, nedich2016achieving,xin2018linear}--$\mathbf{y}^0=\nabla F_c(\mathbf{x}^0)$--is a special case of the proposed initialization.  On the practical side, 
		this initialization  can  be enforced in a distributed way, with minimal coordination. For instance,       agents first choose independently      a vector  $\mathbf{y}_i^{-1}\in \mathbb{R}^m$; then they run one step of consensus on the $y$-variables using the values ${y}_i^{-1}$'s and  weights matrix $\mathbf{C}$,  and set $\mathbf{y}_i^0 = \nabla f_i(\bx_i^0) +\sum_{j\in \mathcal{N}_i^{in}} C_{ij} \mathbf{y}_j^{-1}- \mathbf{y}_i^{-1}$, resulting in $\mathbf{y}^0\in \nabla F_c(\mathbf{x}^0)+ \text{span}(\mathbf{W}_C-\mathbf{I})$.   
		
	}

	Different choices for $\mathbf{R}$ and $\mathbf{C}$ are possible, resulting in different existing algorithms. For instance, if $\mathbf{R}=\mathbf{C}\in\mathcal{M}_n(\mathbb{R})$ are  doubly-stochastic matrices compliant to the graph $\mathcal{G}$, \eqref{NEXT_generic}  reduces to the NEXT algorithm \cite{dilorenzo2015distributed,NEXT}  (or the one in \cite{xu2015augmented}, when \eqref{eq:P} is convex). If $\mathbf{R}$ and $\mathbf{C}$ are allowed to be time-varying (suitably chosen) 
	(\ref{NEXT_generic}) reduces to the SONATA algorithm applicable to (possibly time-varying) digraphs   \cite{sun2016distributed,Scutari_Ying_LectureNote,SONATA_technical_report, YingDanScuConvergence18} [or the one later proposed in \cite{nedich2016achieving} for strongly   convex instances of \eqref{eq:P}]. Finally, if  $\mathbf{R}$ and $\mathbf{C}$ are chosen according to Assumption~\ref{matrix_C_R} below, the scheme  (\ref{NEXT_generic}) becomes the   algorithm proposed independently in     \cite{ANedich_GeoConvAlg_arxive_2018} and \cite{xin2018linear}, for strongly convex objectives  in  \eqref{eq:P}, and implementable over fixed digraphs.

	\begin{assumption}(On the matrices $\mathbf{R}$ and $\mathbf{C}$)\label{matrix_C_R} The weight matrices $\mathbf{R},\mathbf{C}\in \mathcal{M}_n(\mathbb{R})$ satisfy the following:
		\begin{itemize}
			\item[(i)] $\mathbf{R}$ is nonnegative row-stochastic and ${R}_{ii}>0$, for all $i\in [n]$;
			\item[(ii)]  $\mathbf{C}$ is nonnegative column-stochastic and ${C}_{ii}>0$, for all $i\in [n]$;
			\item[(iii)] The graphs   $\mathcal{G}_R$ and $\mathcal{G}_{C^\top}$ each contain    at least one spanning tree; and   $\mathcal{R}_R\cap\mathcal{R}_{{C^\top}}\neq \emptyset$.
		\end{itemize}
	\end{assumption}
	It is not difficult to check that   matrices  $\mathbf{R}$ and $\mathbf{C}$ above exist if and only if   the digraph $\mathcal{G}$ is strongly connected; however, $\mathcal{G}_R$ and $\mathcal{G}_{C^\top}$ need not be so. Several choices for such matrices are discussed in \cite{ANedich_GeoConvAlg_arxive_2018,xin2018linear}. Here, we only point out the following property of $\mathbf{R}$ and $\mathbf{C}$, as a consequence of Assumption \ref{matrix_C_R}, which will be used in our analysis. The result is a consequence of \cite[Lemma 1]{xin2019distributed}. 
	
	
	{
		\begin{lemma}  \label{R_C_norms_plus_realAnalyticity} 
			Given $\mathbf{R}$ and $\mathbf{C}$ satisfying Assumption \ref{matrix_C_R} with  stochastic left eigenvector $\mathbf{r}$  (resp. right eigenvector $\mathbf{c}$ ) of $\mathbf{R}$ (resp. $\mathbf{C}$) associated with the eigenvalue one, then there exist matrix norms\vspace{-0.1cm}
			\begin{align}
			||\mathbf{X}||_R&\triangleq ||\diag(\sqrt{\mathbf{r}})\mathbf{X}\diag(\sqrt{\mathbf{r}})^{-1}||_2,\label{R_matNorm}
			\\
			||\mathbf{X}||_C&\triangleq ||\diag(\sqrt{\mathbf{c}})^{-1}\mathbf{X}\diag(\sqrt{\mathbf{c}})||_2,\label{C_matNorm}
			\end{align}
			such that $\rho_R\triangleq \|\mathbf{R}-\mathbf{1}\mathbf{r}^\top\|_R<1$ and $\rho_C\triangleq \|\mathbf{C}-\mathbf{c}\mathbf{1}^\top\|_C<1$. Furthermore,  $\mathbf{r}^\top\mathbf{c}>0$.
		\end{lemma}
		
		Using Lemma \ref{R_C_norms_plus_realAnalyticity}, it is not difficult to check that the following properties hold:\vspace{-0.1cm}
		\begin{align}
		&\rho_R=\sigma_2\left(\diag(\sqrt{\mathbf{r}})\mathbf{R}\diag(\sqrt{\mathbf{r}})^{-1}\right),
		\label{explicit_rho_R}
		\\
		&\rho_C=\sigma_2\left(\diag(\sqrt{\mathbf{c}})^{-1}\mathbf{C}\diag(\sqrt{\mathbf{c}})\right),\label{explicit_rho_C}
		\\
		&||\mathbf{R}||_R=||\mathbf{1r}^T||_R=||\mathbf{I}-\mathbf{1r}^T||_R=1,
		\label{explicit_Rnorm_Eqs}
		\\
		&||\mathbf{C}||_C=||\mathbf{c1}^T||_C=||\mathbf{I}-\mathbf{c1}^T||_R=1.
		\label{explicit_Cnorm_Eqs}
		\end{align}

		The vector norms associated with above matrix norms are\vspace{-0.1cm}
		\begin{align}
		||\mathbf{x}||_R&=||\diag(\sqrt{\mathbf{r}})\mathbf{x}||_2,\label{R_vecNorm}
		\\
		||\mathbf{x}||_C&=||\diag(\sqrt{\mathbf{c}})^{-1}\mathbf{x}||_2;\label{C_vecNorm}
		\end{align}
		and $||\cdot ||_a\leq K_{a,b}||\cdot||_b$ holds for $a,b\in\{R,C,2\}$ with\vspace{-0.1cm}
		\begin{equation}\label{eq:def_K_1}
		\begin{aligned}
		&K_{R,2}=\sqrt{r_{\max}},\quad K_{2,R}=1/\sqrt{r_{\min}},
		\\
		&K_{C,2}=1/\sqrt{c_{\min}},\quad K_{2,C}=\sqrt{c_{\max}},
		\\
		&K_{R,C}=\sqrt{r_{\max}c_{\max}},\quad K_{C,R}=1/\sqrt{c_{\min}r_{\min}},
		\end{aligned}
		\end{equation}
		where $r_{\min}$ (resp. $c_{\min}$) and $r_{\max}$ (resp. $c_{\max}$) are minimum and maximum elements of $\mathbf{r}$ (resp. $\mathbf{c}$).
		
	}

	Convergence of   DOGT algorithms  in the form  \eqref{NEXT_generic} (with  $\mathbf{R}$ and $\mathbf{C}$ satisfying Assumption~\ref{matrix_C_R}) has not been studied in the literature when   $F$ is nonconvex. In next subsection we fill this gap and provide a full characterization of the convergence behavior of DOGT including its second-order guarantees.
	
	
	\subsection{First-order   convergence \& rate analysis}
	\label{Subsequence_Conv_Subsec}  
	In this section, we study asymptotic convergence to first-order stationary solutions; we assume  $m=1$ (scalar optimization variables); while this simplifies the notation, all the conclusions hold for the general case $m>1$. 
	As in \cite{ANedich_GeoConvAlg_arxive_2018}, define the weighted sums  \vspace{-0.1cm} \begin{equation}\label{eq:def_w_vectors}
	\bar{{x}}^{\nu}\triangleq \mathbf{r}^\top\mathbf{x}^\nu,\quad   \bar{{y}}^{\nu}\triangleq \mathbf{1}^\top\mathbf{y}^\nu,\quad \text{and}\quad \bar{{g}}^{\nu}\triangleq \mathbf{1}^\top\nabla F_c(\mathbf{x}^\nu),\vspace{-0.1cm}
	\end{equation}
	where we recall that  $\mathbf{r}$  is  the Perron vector  associated with   $\mathbf{R}$ (cf.~Lemma~\ref{R_C_norms_plus_realAnalyticity}). { Note that $\nabla F_c$ is $L_c$-Lipschitz continuous with  $L_c\triangleq  L_\mathrm{max}$.
	}
	
	Using \eqref{NEXT_generic}, it is not difficult to check that the following holds 
	\vspace{-0.1cm}
	\begin{equation}
	\label{x_bar_eq}
	\begin{aligned}
	\bar{{x}}^{\nu+1}
	=\bar{{x}}^{\nu}-\zeta\alpha\bar{y}^\nu-\alpha\mathbf{r}^\top\left(\mathbf{y}^\nu-\mathbf{c}\bar{y}^\nu \right)\quad \text{and}\quad \bar{y}^\nu=\bar{g}^\nu,
	\end{aligned}\vspace{-0.1cm}
	\end{equation}
	where $\mathbf{c}$  is  the Perron vector  associated with   $\mathbf{C}$, and 
	$\zeta\triangleq \mathbf{r}^\top\mathbf{c}>0$  (cf.~Lemma \ref{R_C_norms_plus_realAnalyticity}). 
	
	\subsubsection{Descent on $F$} Using the descent lemma along with    \eqref{x_bar_eq} yields
	\begin{equation*}
	\begin{aligned}
	F(\bar{{x}}^{\nu+1})
	& = F\left(\bar{{x}}^{\nu}-\zeta\alpha\bar{y}^\nu-\alpha\mathbf{r}^\top\left(\mathbf{y}^\nu-\mathbf{c}\bar{y}^\nu \right)\right)
	\\
	& \leq   F(\bar{{x}}^{\nu})-\zeta\alpha\left\langle\nabla F({\bar{x}}^{\nu}),\bar{{y}}^{\nu}\right\rangle-\alpha\left\langle\nabla F({\bar{x}}^{\nu}),\mathbf{r}^\top\left(\mathbf{y}^\nu-\mathbf{c}\bar{y}^\nu \right)\right\rangle
	\\
	&\quad +\frac{L}{2}\norm{\zeta\alpha\bar{y}^\nu+\alpha\mathbf{r}^\top\left(\mathbf{y}^\nu-\mathbf{c}\bar{y}^\nu \right)}^2.
	\end{aligned}
	\end{equation*}
	Adding/subtracting suitably chosen  terms we obtain 
	\begin{equation}
	\begin{aligned}
	F(\bar{{x}}^{\nu+1})
	\leq & F(\bar{{x}}^{\nu})-\zeta\alpha\left\langle\nabla F({\bar{x}}^{\nu})-\bar{y}^\nu,\bar{{y}}^{\nu}\right\rangle-\zeta\alpha|\bar{y}^\nu|^2
	\\
	&-\alpha\left\langle\nabla F({\bar{x}}^{\nu})-\bar{y}^\nu,\mathbf{r}^\top\left(\mathbf{y}^\nu-\mathbf{c}\bar{y}^\nu \right)\right\rangle-\alpha \left\langle \bar{y}^\nu,\mathbf{r}^\top\left(\mathbf{y}^\nu-\mathbf{c}\bar{y}^\nu \right)\right\rangle
	\\
	&+L\zeta^2\alpha^2|\bar{y}^\nu|^2+L\alpha^2\norm{\mathbf{y}^\nu-\mathbf{c}\bar{y}^\nu}^2
	\\
	\leq & F(\bar{{x}}^{\nu})+\frac{\zeta\alpha}{2\epsilon_1}|\nabla F({\bar{x}}^{\nu})-\bar{y}^\nu|^2+\frac{\zeta\alpha\epsilon_1}{2}|\bar{y}^{\nu}|^2
	-\zeta\alpha|\bar{y}^\nu|^2
	\\
	&+\frac{\alpha}{2}|\nabla F({\bar{x}}^{\nu})-\bar{y}^\nu|^2+\frac{\alpha}{2}\norm{\mathbf{y}^\nu-\mathbf{c}\bar{y}^\nu}^2
	+\frac{\alpha \epsilon_2}{2}|\bar{y}^{\nu}|^2+\frac{\alpha }{2\epsilon_2}\norm{\mathbf{y}^\nu-\mathbf{c}\bar{y}^\nu}^2
	\\
	&+L\zeta^2\alpha^2|\bar{y}^\nu|^2+L\alpha^2\norm{\mathbf{y}^\nu-\mathbf{c}\bar{y}^\nu}^2
	\\
	=& F(\bar{{x}}^{\nu})+\left(\frac{\zeta\alpha\epsilon_1}{2}-\zeta\alpha+\frac{\alpha \epsilon_2}{2}+L\zeta^2\alpha^2\right)|\bar{y}^{\nu}|^2
	\\
	&+\left(\frac{\zeta\alpha}{2\epsilon_1}+\frac{\alpha}{2}\right)|\nabla F({\bar{x}}^{\nu})-\bar{y}^\nu|^2
	+\left(\frac{\alpha}{2}+\frac{\alpha }{2\epsilon_2}+L\alpha^2\right)\norm{\mathbf{y}^\nu-\mathbf{c}\bar{y}^\nu}^2,
	\end{aligned}
	\label{descent1}
	\end{equation}
	where $\epsilon_1$ and $\epsilon_2$ are some arbitrary positive quantities (to be chosen). { 
		By $\bar{y}^\nu=\bar{g}^\nu$ [cf.~\eqref{x_bar_eq}], it holds that  \vspace{-0.1cm}
		\begin{equation}
		|\nabla F({\bar{x}}^{\nu})-\bar{y}^\nu|=  \left|\sum_{i=1}^n \nabla f_i(\bar{{x}}^\nu)-\sum_{i=1}^n \nabla f_i({x}_i^\nu)\right|
		\leq   L_c\sqrt{n} \norm{\mathbf{x}^\nu-\mathbf{1}\bar{x}^\nu}.
		\label{grad_error_bound}
		\end{equation}
	}
	
	Combining   \eqref{descent1} and \eqref{grad_error_bound} yields
	{
		\begin{equation}
		\begin{aligned}
		& F(\bar{{x}}^{\nu+1}) 
		\\
		\leq &\,
		F(\bar{{x}}^{\nu})+\left(\frac{\zeta\alpha\epsilon_1}{2}-\zeta\alpha+\frac{\alpha \epsilon_2}{2}+L\zeta^2\alpha^2\right)|\bar{y}^{\nu}|^2
		\\
		&  +nL_c^2K_{2,R}^2\left(\frac{\zeta\alpha}{2\epsilon_1}+\frac{\alpha}{2}\right)\norm{\mathbf{x}^\nu-\mathbf{1}\bar{x}^\nu}_R^2
		+K_{2,C}^2\left(\frac{\alpha}{2}+\frac{\alpha }{2\epsilon_2}+L\alpha^2\right)\norm{\mathbf{y}^\nu-\mathbf{c}\bar{y}^\nu}_C^2,
		\end{aligned}
		\label{descent2}
		\end{equation}
		where $K_{2,R}=1/\sqrt{r_{\min}}$ and $K_{2,C}=\sqrt{c_{\max}}$  [cf. \eqref{eq:def_K}].}
	
	
	\subsubsection{Bounding the consensus and gradient tracking errors} 
	Let us bound  the consensus error $\|\mathbf{x}^\nu-\mathbf{1}\bar{x}^\nu\|_R$. Using $\|\mathbf{z}+\mathbf{w}||^2_R\leq (1+\epsilon)\,\|\mathbf{x}\|^2_R+(1+{1}/{\epsilon})\,\|\mathbf{y}\|^2_R$, for arbitrary $\mathbf{z}, \mathbf{w}\in \mathbb{R}^m$ and   $\epsilon>0$, along with Lemma~\ref{R_C_norms_plus_realAnalyticity}, yields
	{
		\begin{equation}
		\begin{aligned}
		& \norm{\mathbf{x}^{\nu+1}-\mathbf{1}\bar{x}^{\nu+1}}^2_R
		=   \norm{\left(\mathbf{R}-\mathbf{1}\mathbf{r}^\top\right)\left(\mathbf{x}^\nu-\mathbf{1}\bar{x}^\nu\right)-\alpha\left(\mathbf{I}-\mathbf{1}\mathbf{r}^\top\right)\left(\mathbf{y}^{\nu}-\mathbf{1}\bar{y}^\nu\right)}^2_R
		\\
		&\leq   (1+\epsilon_x)\norm{\left(\mathbf{R}-\mathbf{1}\mathbf{r}^\top\right)\left(\mathbf{x}^\nu-\mathbf{1}\bar{x}^\nu\right)}^2_R+\alpha^2 \left(1+\frac{1}{\epsilon_x}\right)\norm{\left(\mathbf{I}-\mathbf{1}\mathbf{r}^\top\right)\left(\mathbf{y}^{\nu}-\mathbf{1}\bar{y}^\nu\right)}^2_R
		\\
		&\leq   \rho_R^2 (1+\epsilon_x) \norm{\mathbf{x}^\nu-\mathbf{1}\bar{x}^\nu}^2_R+\alpha^2 \left(1+\frac{1}{\epsilon_x}\right)\|\mathbf{I}-\mathbf{1}\mathbf{r}^\top \|_R^2\norm{\mathbf{y}^{\nu}-\mathbf{1}\bar{y}^\nu}^2_R
		\\
		&\overset{\eqref{explicit_Rnorm_Eqs}}{\leq}   \rho_R^2  (1+\epsilon_x) \norm{\mathbf{x}^\nu-\mathbf{1}\bar{x}^\nu}^2_R+2\alpha^2  \left(1+\frac{1}{\epsilon_x}\right)\,\norm{\mathbf{y}^{\nu}-\mathbf{c}\bar{y}^\nu}^2_R
		\\
		&	\qquad  +2\alpha^2  (1+\frac{1}{\epsilon_x})\norm{(\mathbf{1}-\mathbf{c})\bar{y}^\nu}^2_R
		\\
		& \leq   \rho^2_R (1+\epsilon_x) \norm{\mathbf{x}^\nu-\mathbf{1}\bar{x}^\nu}^2_R+\alpha^2 K_2\norm{\mathbf{y}^{\nu}-\mathbf{c}\bar{y}^\nu}^2_C+\alpha^2 K_3|\bar{y}^\nu|^2_2,
		\end{aligned}
		\label{x_bound}
		\end{equation}
		where $\epsilon_x>0$ is arbitrary and  we defined 
		\begin{equation}\label{K_consts_1}
		K_2\triangleq 2K_{R,C}^2  \left(1+\frac{1}{\epsilon_x}\right),\quad K_3\triangleq 
		2n(1+\frac{1}{\epsilon_x}).
		\end{equation}} 
	
	Similarly, the tracking error can be bounded as  
	{
		\begin{equation}\hspace{-.2cm}
		\label{y_bound}
		\begin{aligned}
		&\norm{\mathbf{y}^{\nu+1}-\mathbf{c}\bar{y}^{\nu+1}}^2_C
		=  \norm{\left(\mathbf{C}-\mathbf{c}\mathbf{1}^\top\right)\mathbf{y}^{\nu}+\left(\mathbf{I}-\mathbf{c}\mathbf{1}^\top\right)\left(\nabla F_c(\mathbf{x}^{\nu+1})-\nabla F_c(\mathbf{x}^{\nu})\right)}^2_C
		\\
		& \leq   (1+\epsilon_y)\norm{\left(\mathbf{C}-\mathbf{c}\mathbf{1}^\top\right)\left(\mathbf{y}^{\nu}-\mathbf{c}\bar{y}^\nu\right)}^2_C
		\\
		&\quad +(1+\frac{1}{\epsilon_y})\norm{\left(\mathbf{I}-\mathbf{c}\mathbf{1}^\top\right)\left(\nabla F_c(\mathbf{x}^{\nu+1})-\nabla F_c(\mathbf{x}^{\nu})\right)}^2_C
		\\
		& \overset{\eqref{explicit_Cnorm_Eqs}}{\leq} \rho_C^2(1+\epsilon_y)\norm{\mathbf{y}^{\nu}-\mathbf{c}\bar{y}^\nu}^2_C+K_{C,2}^2L_{c}^2\left(1+\frac{1}{\epsilon_y}\right)\,\norm{\mathbf{x}^{\nu+1}-\mathbf{x}^{\nu}}^2
		\\
		& \overset{(a)}{=} \rho_C^2(1+\epsilon_y)\norm{\mathbf{y}^{\nu}-\mathbf{c}\bar{y}^\nu}^2_C
		\\
		&+3K_{C,2}^2L_{c}^2\left(1+\frac{1}{\epsilon_y}\right)\,\left[\norm{(\mathbf{R}-\mathbf{I})(\mathbf{x}^\nu-\mathbf{1}\bar{x}^\nu)}^2+\alpha^2\norm{\mathbf{y}^{\nu}-\mathbf{c}\bar{y}^\nu}^2 +\alpha^2|\bar{y}^\nu|^2\norm{\mathbf{c}}^2\right]
		\\
		&  \overset{\eqref{explicit_Rnorm_Eqs}}{\leq} \rho_C^2(1+\epsilon_y)\norm{\mathbf{y}^{\nu}-\mathbf{c}\bar{y}^\nu}^2_C
		\\
		&+3K_{C,2}^2L_{c}^2\left(1+\frac{1}{\epsilon_y}\right)\,\left[K_{2,R}^2\norm{\mathbf{x}^\nu-\mathbf{1}\bar{x}^\nu}^2+K_{2,C}^2\alpha^2\norm{\mathbf{y}^{\nu}-\mathbf{c}\bar{y}^\nu}_C^2 +\alpha^2|\bar{y}^\nu|^2\right]
		\\
		&  = \rho_C^2(1+\epsilon_y)\norm{\mathbf{y}^{\nu}-\mathbf{c}\bar{y}^\nu}^2_C
		\\
		&+3K_{C,2}^2L_{c}^2\left(1+\frac{1}{\epsilon_y}\right)\,\left[K_{2,R}^2\norm{\mathbf{x}^\nu-\mathbf{1}\bar{x}^\nu}^2+K_{2,C}^2\alpha^2\norm{\mathbf{y}^{\nu}-\mathbf{c}\bar{y}^\nu}_C^2 +\alpha^2|\bar{y}^\nu|^2\right]
		\\
		& \leq   \left(\rho^2_C+\frac{\alpha^2 K_4}{\epsilon_y}\right)(1+\epsilon_y) \norm{\mathbf{y}^{\nu}-\mathbf{c}\bar{y}^\nu}^2_C+\alpha^2 K_5|\bar{y}^\nu|^2_2+K_6\left(1+\frac{1}{\epsilon_y}\right)\norm{\mathbf{x}^\nu-\mathbf{1}\bar{x}^\nu}^2_R,
		\end{aligned}
		\end{equation}
		where in (a) we used $\bx^{\nu+1}-\bx^\nu=(\mathbf{R}-\mathbf{I})(\mathbf{x}^\nu-\mathbf{1}\bar{x}^\nu)-\alpha(\mathbf{y}^{\nu}-\mathbf{c}\bar{y}^\nu)-\alpha\mathbf{c}\bar{y}^\nu$ and the Jensen's inequality; and in the last inequality we defined 
		\begin{equation}\label{K_consts_2}
		K_4=3K_{C,2}^2K_{2,C}^2L_c^2,\quad K_5=3K_{C,2}^2L_c^2,\quad K_6=3K_{C,2}^2K_{2,R}^2L_c^2.	
		\end{equation}
		
	}


%
%

\subsubsection{Lyapunov function}
{
Let us introduce now the candidate Lyapunov function: denoting $\mathbf{J}_R\triangleq \mathbf{1}\mathbf{r}^\top$ and $\mathbf{J}_C\triangleq \mathbf{c}\mathbf{1}^\top$, define
\begin{equation}\label{eq:def_L}
L(\mathbf{x},\mathbf{y})\triangleq F_c(\mathbf{J}_R\mathbf{x})+\norm{(\mathbf{I}-\mathbf{J}_R)\mathbf{x}}^2_R+\varkappa\norm{(\mathbf{I}-\mathbf{J}_C)\mathbf{y}}^2_C,
\end{equation} 	
where  $\varkappa>0$ is a positive constant (to be properly chosen). Combining  \eqref{descent2}-\eqref{y_bound} and using $\bar{y}^\nu=\bar{g}^\nu=\sum_{i=1}^n\nabla f_i(x_i^\nu)$ [cf. \eqref{x_bar_eq}] leads to the following descent property for $L$:
\begin{equation}
L(\mathbf{x}^{\nu+1},\mathbf{y}^{\nu+1})\leq L(\mathbf{x}^{\nu},\mathbf{y}^{\nu})-d(\mathbf{x}^{\nu},\mathbf{y}^{\nu})^2,
\label{L_descent} 
\end{equation}
where
\begin{equation}	\label{d_def_func}
d(\bx,\by)\triangleq \sqrt{(1-\tilde{\rho}_R)\norm{(\mathbf{I}-\mathbf{J}_R)\mathbf{x}}^2_R+\varkappa(1-\tilde{\rho}_C)\norm{(\mathbf{I}-\mathbf{J}_C)\by}^2_C+\Gamma~\Big|\sum_{i=1}^n\nabla f_i(x_i)\Big|^2}
\end{equation}	
and 
	\begin{equation}
	\begin{aligned}
	\tilde{\rho}_R\triangleq &\rho^2_R (1+\epsilon_x)+\frac{\alpha nL_c^2K_{2,R}^2}{2}\left(1+\frac{\zeta}{\epsilon_1}\right)+\varkappa K_6\left(1+\frac{1}{\epsilon_y}\right),
	\\
	\tilde{\rho}_C\triangleq &\rho^2_C(1+\epsilon_y)+ \frac{\alpha K_{2,C}^2}{2\varkappa}\left(1+\frac{1}{\epsilon_2}\right)+\alpha^2\left(\frac{LK_{2,C}^2+K_2}{\varkappa}+K_4\left(1+\frac{1}{\epsilon_y}\right)\right),
	\\
	\Gamma\triangleq &\left(\zeta-\frac{\epsilon_1\zeta}{2}-\frac{\epsilon_2 }{2}\right)\alpha-\left(L\zeta^2+K_3+K_5\varkappa\right)\alpha^2.
	\end{aligned}
	\end{equation}

Note that the function $d(\bullet,\bullet)$ is a valid measure of optimality/consensus for DOGT: i) it is continuous; and ii)    $d(\bx,\by)=0$  implies $x_i=x_j=x^\ast$, for all $i,j\in[n]$ and some $x^\ast$ such that $\sum_{i=1}^n\nabla f_i(x^\ast)=0$, meaning that all  $x_i$ are consensual and equal to a critical point of $F$. }

{
To ensure   $\tilde{\rho}_R<1$, $\tilde{\rho}_C<1$, and $\Gamma>0$ in $d(\bx,\by)$, we choose the free parameters $\epsilon_x$, $\epsilon_y$, $\epsilon_1$, $\epsilon_2$, and $\varkappa$ as follows: 
\begin{equation}
\begin{aligned}
& 0<\epsilon_x<\frac{1-\rho_R^2}{2\rho_R^2},  \quad   & 0<\epsilon_y<\frac{1-\rho_C^2}{\rho_C^2},\quad &
\\
& \epsilon_1=\epsilon_2=\epsilon,\quad & 0<\epsilon<\frac{2\zeta}{1+\zeta},\qquad  & 0<\varkappa\leq \frac{\rho^2_R\epsilon_x}{K_6(1+1/\epsilon_y)},
\end{aligned}
\label{parameters_assymptotic_conv}
\end{equation}
and finally,   $\alpha>0$ must satisfy 
	\begin{equation}
	\begin{aligned}
	\alpha & <\frac{2}{nL_c^2K_{2,R}^2\left(1+\frac{\zeta}{\epsilon}\right)}\left(1-\rho^2_R (1+2\epsilon_x)\right),
	\\
	\alpha & <\frac{1-\rho^2_C (1+\epsilon_y)}{\frac{1}{2\varkappa}K_{2,C}^2\left(1+\frac{1}{\epsilon}+2L\right)+\frac{K_2}{\varkappa}+K_4\left(1+\frac{1}{\epsilon_y}\right)},
	\\
	\alpha & <\frac{\zeta-\frac{\epsilon}{2}\left(\zeta+1\right)}{L\zeta^2+K_3+K_5\varkappa}.
	\end{aligned}
	\label{alpha_conds_assymptotic_conv}
	\end{equation}
	 Substituting    \eqref{eq:def_K_1}, \eqref{K_consts_1}, and \eqref{K_consts_2}     in (\ref{alpha_conds_assymptotic_conv}) and setting for simplicity 	
	\begin{equation}\label{alpha_conds_assymptotic_conv_2}
	\epsilon_x=\frac{1-\rho_R^2}{4\rho_R^2},\quad \epsilon_y=\frac{1-\rho_C^2}{2\rho_C^2},\quad\epsilon=\frac{\xi}{1+\xi},\quad \varkappa= \frac{c_{\min}r_{\min}}{24L_c^2}(1-\rho_R^2)(1-\rho_C^2),
	\end{equation}
	we obtain   the following sufficient conditions for (\ref{alpha_conds_assymptotic_conv}):
	\begin{equation}
	\begin{aligned}
	\alpha & \leq \tilde{\alpha}_1\triangleq  \frac{r_{\min}(1-\rho_R^2)}{3nL_c^2},
	\\
	\alpha & \leq \tilde{\alpha}_2\triangleq\frac{(1-\rho_R^2)^2(1-\rho_C^2)^2r_{\min}^2c_{\min}^2}{1152L_c^2(2+L)},
	\\
	\alpha & \leq \tilde{\alpha}_3\triangleq \frac{r_{\min}c_{\min}(1-\rho_R^2)}{2(L+16n)}.
	\end{aligned}
	\label{alpha_conds_assymptotic_conv_3}
	\end{equation}
	A further simplification, leads to the   following final  more restrictive  condition on $\alpha$:
	\begin{equation}
	\begin{aligned}
	0<\alpha & \leq \frac{(1-\rho_R^2)^2(1-\rho_C^2)^2r_{\min}^2c_{\min}^2}{1152L_c^2(L+16n)}.
	\end{aligned}
	\label{alpha_conds_assymptotic_conv_final}
	\end{equation}

}

The descent property (\ref{L_descent}) readily implies the following convergence result for $\{L(\mathbf{x}^\nu,\mathbf{y}^\nu)\}$ and $\{d(\bx^\nu,\by^\nu)\}$. 
\begin{lemma}
	\label{sequence_props_a}  
	Under Assumptions \ref{P_assumption},  \ref{Net-Assump}, and \ref{matrix_C_R}, and the above choice of parameter, there hold:
	\begin{enumerate}[label=(\roman*)]
		\item The sequence $\{L(\mathbf{x}^\nu,\mathbf{y}^\nu)\}$ converges; 
		\item $\sum_{\nu=0}^\infty d(\bx^\nu,\by^\nu)^2<\infty$, and  thus $\lim_{\nu\to \infty}d(\bx^\nu,\by^\nu)= 0$. 
	\end{enumerate}
\end{lemma}

{ We conclude this subsection by lower bounding  $d(\bx^\nu,\by^\nu)$  by the magnitude of the gradient of the Lyapunov function $L$. This will allow us to transfer the convergence properties of   $\{d(\bx^\nu,\by^\nu)\}$ to $\{||\nabla L(\bx^\nu,\by^\nu)||\}$. The  lemma below will also be useful to  establish global convergence of DOGT under the {K\L} property  (cf. Sec.  \ref{sec:KL-main-results}).}

{
	\begin{lemma}
		\label{nablaL_bound_lemma}  Let  $\nabla L(\mathbf{x}^\nu,\mathbf{y}^\nu)\triangleq (\nabla_\mathbf{x} L(\mathbf{x}^\nu,\mathbf{y}^\nu),\nabla_\mathbf{y} L(\mathbf{x}^\nu,\mathbf{y}^\nu))$, where $\nabla_\mathbf{x} L$ (resp. $\nabla_\mathbf{y} L$) are the gradient of $L$ with respect to the first (resp. second) argument.  In the setting above, there holds\vspace{-0.2cm}
		\begin{equation}
		\norm{\nabla L(\mathbf{x}^\nu,\mathbf{y}^\nu)}\leq Md(\bx^\nu,\by^\nu),\quad  \nu\geq 0,
		\label{nablaL_bound_original}
		\end{equation}
		with \vspace{-0.3cm}
		\begin{equation}
		M=\sqrt{2}\max\left(\frac{(2r_{\max}+L_c\sqrt{n})^2}{r_{\min}(1-\tilde{\rho}_R)},\frac{2\varkappa c_{\max}}{c_{\min}^{2}(1-\tilde{\rho}_C)}, \frac{1}{\Gamma}\right)^{\frac{1}{2}}.
		\end{equation}
	\end{lemma}
	\begin{proof}
		Recall that $\mathbf{J}_R=\mathbf{1r}^\top$ and $\mathbf{J}_C=\mathbf{c1}^\top$. 
		By definition \eqref{eq:def_L} and Lemma \ref{R_C_norms_plus_realAnalyticity}, we can write\vspace{-0.1cm}
		\begin{equation}
		\begin{aligned}
		\nabla_\mathbf{x} L(\mathbf{x}^\nu,\mathbf{y}^\nu) &=\mathbf{J}_R^\top\nabla F_c(\mathbf{J}_R\mathbf{x}^\nu)+2(\mathbf{I}-\mathbf{J}_R)^\top\diag({\mathbf{r}})(\mathbf{I}-\mathbf{J}_R)\mathbf{x}^\nu
		\\
		& \overset{(a)}{=} \mathbf{r}\,\bar{y}^\nu+\mathbf{J}_R^\top\left(\nabla F_c(\mathbf{J}_R\mathbf{x}^\nu)-\nabla F_c(\mathbf{x}^\nu)\right)
		\\
		&\quad +2(\mathbf{I}-\mathbf{J}_R)^\top\diag({\mathbf{r}})(\mathbf{x}^\nu-\mathbf{1}\bar{x}^\nu),\smallskip
		\\ 
		\nabla_\mathbf{y} L(\mathbf{x}^\nu,\mathbf{y}^\nu) &=2\varkappa(\mathbf{I}-\mathbf{J}_C)^\top\diag({\mathbf{c}})^{-1}(\mathbf{I}-\mathbf{J}_C)\mathbf{y}^\nu\\
		& = 2\varkappa(\mathbf{I}-\mathbf{J}_C)^\top\diag({\mathbf{c}})^{-1}(\mathbf{y}^\nu-\mathbf{c}\bar{y}^\nu),
		\end{aligned}\label{nablaL_expression}\vspace{-0.1cm}
		\end{equation}
		where (a) is due to $\bar{y}^\nu=\bar{g}^\nu$ (cf. \eqref{x_bar_eq}). Thus there holds
		\begin{equation}
		\begin{aligned}
		||\nabla_\mathbf{x} L(\mathbf{x}^\nu,\mathbf{y}^\nu)|| &\leq ||\mathbf{r}||~ |\bar{y}^\nu|+||\mathbf{J}_R^\top\left(\nabla F_c(\mathbf{J}_R\mathbf{x}^\nu)-\nabla F_c(\mathbf{x}^\nu)\right)||
		\\
		&\quad +2||(\mathbf{I}-\mathbf{J}_R)^\top\diag({\mathbf{r}})(\mathbf{x}^\nu-\mathbf{1}\bar{x}^\nu)||
		\\
		& \overset{(b)}{\leq} |\bar{y}^\nu|+K_{2,R}\left(2r_{\max}+L_c\sqrt{n}\right)||\mathbf{x}^\nu-\mathbf{1}\bar{x}^\nu||_R,
		\\
		||\nabla_\mathbf{y} L(\mathbf{x}^\nu,\mathbf{y}^\nu)|| &\overset{(c)}{\leq}  2\varkappa K_{2,C}c_{\min}^{-1}||\mathbf{y}^\nu-\mathbf{c}\bar{y}^\nu||_C,
		\end{aligned}
		\label{nablaL_expression2}
		\end{equation}
		where (b) holds due to $||\diag(\mathbf{r)}||_R=||\diag(\mathbf{r})||_2=r_{\max}$, $||\mathbf{r}||\leq 1$, $||\mathbf{J}_R||_2\leq \sqrt{n}$ and \eqref{explicit_Rnorm_Eqs};  (c) is due to $||\diag(\mathbf{c})^{-1}||_C=||\diag(\mathbf{c})^{-1}||_2=c_{\min}^{-1}$ and  \eqref{explicit_Cnorm_Eqs}.
		Eq.~\eqref{nablaL_bound_original} follows readily from \eqref{nablaL_expression2}. 
\end{proof}}

\subsubsection{Main result}
We can now state the main convergence result of DOGT to critical points of $F$.    
{ \begin{theorem}\label{Th_asympt_convergence}
		\label{sequence_props} Consider Problem \eqref{eq:P}, and suppose that Assumptions \ref{P_assumption} and \ref{Net-Assump}  are satisfied. 
		Let     $\{(\mathbf{x}^\nu,\mathbf{y}^\nu)\}$ be the sequence generated by the DOGT Algorithm   \eqref{NEXT_generic}, with $\mathbf{R}$ and $\mathbf{C}$ satisfying Assumption \ref{matrix_C_R},   and $\alpha$ chosen according to \eqref{alpha_conds_assymptotic_conv_final} [or \eqref{alpha_conds_assymptotic_conv_2}]; let $\{\bar{x}^\nu\}$ and $\{\bar{y}^\nu\}$  be defined in \eqref{eq:def_w_vectors}; and let $\{d(\bx^\nu,\by^\nu)\}$ be defined in  \eqref{d_def_func}. Given   $\epsilon>0$, let    $T_\epsilon=\min\{\nu\in \mathbb{N}_+\,:\, d(\bx^\nu,\by^\nu)\leq \epsilon\}$.    Then, there hold 
		\begin{enumerate}[label=(\roman*)] 
			\item[(i)] \emph{[consensus]:} $\lim_{\nu\to \infty} \norm{\mathbf{x}^\nu-\mathbf{1}\bar{x}^\nu}=0$ and $\lim_{\nu\to \infty}  \bar{y}^\nu =0$;
			\item[(ii)] \emph{[stationarity]:} let $\mathbf{x}^\infty$ be a limit point of $\{\mathbf{x}^\nu\}$; then,   $\mathbf{x}^\infty=\theta^\infty \,\mathbf{1}$, for some   $\theta^\infty\in\mathrm{crit}~F$; 
			\item[(iii)] \emph{[sublinear rate]:}  $T_\epsilon=o(1/\epsilon^2)$. 
		\end{enumerate}
	\end{theorem}
	\begin{proof}
		(i) follows readily from  Lemma~\ref{sequence_props_a}(ii).
		
		We prove   (ii). Let  $(\mathbf{x}^\infty, \mathbf{y}^\infty)$ be a limit point of $\{(\mathbf{x}^\nu, \mathbf{y}^\nu)\}$. By (i), it must be $(\mathbf{I}-\mathbf{J}_R){\mathbf{x}}^\infty=\mathbf{0}$, implying    $ {\mathbf{x}}^\infty=\mathbf{1}\theta^\infty$, for some $\theta^\infty\in\mathbb{R}$.  Also, $\lim_{\nu\to \infty} \mathbf{1}^\top\nabla  F_c(\mathbf{x}^\nu)= \lim_{\nu\to \infty}  \bar{g}^\nu=\lim_{\nu\to \infty} \bar{y}^\nu=0$, which together  with   the continuity   of  $\nabla F_c$, yields $0=\mathbf{1}^\top\nabla F_c(\mathbf{1}\theta^\infty)=\nabla F(\theta^\infty)$. Therefore, $\theta^\infty\in\mathrm{crit}~F$.
		
		We prove now (iii). Using \eqref{L_descent} and the definition of  $T_\epsilon$, we can write \vspace{-0.3cm}
		\begin{equation}\label{T_eps_bound}
		\frac{T_\epsilon}{2}\epsilon^2	\leq \sum_{t=\floor{\frac{T_\epsilon}{2}}+1}^{T_\epsilon} d(\bx^t,\by^t)^2\leq l^{\floor{\frac{T_\epsilon}{2}}+1}-l^{T_\epsilon+1},\vspace{-0.2cm}
		\end{equation}
		where we  used the  shorthand $l^\nu\triangleq L(\mathbf{x}^\nu,\mathbf{y}^\nu)$. Consider the   following two cases: (1)   $T_\epsilon\rightarrow \infty$ as $\epsilon\rightarrow 0$, then   $l^{\floor{\frac{T_\epsilon}{2}}+1}-l^{T_\epsilon+1}\rightarrow 0$ (recall that $\{l^\nu\}$ converges, cf. Lemma \ref{sequence_props_a}(i)); and (2) $T_\epsilon<\infty$ as $\epsilon\rightarrow 0$, then  $\{l^\nu\}$ converges in a finite number of iterations. Therefore, by \eqref{T_eps_bound}, we have   $T_\epsilon=o(1/\epsilon^2)$.  		\end{proof}}

Note that, as a direct consequence of Lemma~\ref{nablaL_bound_lemma}, one can infer the following further property of the limit points  $(\mathbf{x}^\infty, \mathbf{y}^\infty)$ of the sequence  $\{(\mathbf{x}^\nu, \mathbf{y}^\nu)\}$: any such a $(\mathbf{x}^\infty, \mathbf{y}^\infty)$ is a critical point of $L$ [defined in (\ref{eq:def_L})]. 





	
\subsection{Convergence under the {K\L}  property} \label{Global_Conv_Subsec}
{
We now strengthen the subsequence convergence result in Theorem~\ref{Th_asympt_convergence}, under the additional assumption that    $F$ is a {K\L} function \cite{MR0160856,Kurdyka1998}: We prove that   the entire  sequence $\{\mathbf{x}^\nu\}$  converges to a critical point of $F$ (cf. Theorem \ref{main_global_conv_thm}), and establish   asymptotic  convergence rates (cf. Theorem \ref{Rate_theorem}).   We extend the  analysis developed in   \cite{Attouch2009, Attouch2013} for centralized first-order  methods to our distributed setting and complement it with a rate analysis.   The major difference with  \cite{Attouch2013} is that the     \emph{sufficient decent} condition postulated in \cite{Attouch2013} is neither satisfied by the objective value sequence $\{F(\mathbf{x}^\nu)\}$ (as requested in \cite{Attouch2013}), due to  consensus and gradient tracking errors, nor  by the Lyapunov function sequence $\{L(\mathbf{x}^\nu,\mathbf{y}^\nu)\}$, which instead satisfies \eqref{L_descent}. 
 A key step to cope with this issue is to establish necessary connections between   $\nabla L(\bx,\by)$ and $d(\bx,\by)$ (defined in \eqref{eq:def_L} and \eqref{d_def_func}, respectively)--see  Lemma \ref{nablaL_bound_original} and Proposition \ref{pre_conv_result}. }

\subsubsection{Convergence analysis} \label{sec:KL-main-results}
We begin proving the following abstract intermediate results similar to   \cite{Attouch2013} but extended to our distributed setting, which is at the core of the subsequent analysis; we still assume $m=1$ without loss of generality. 

\begin{proposition}
	\label{pre_conv_result}
	In the  setting of Theorem~ \ref{Th_asympt_convergence}, let $L$ defined in (\ref{eq:def_L}) is {K\L}  at some $\acute{\mathbf{z}}\triangleq (\acute{\mathbf{x}},\acute{\mathbf{y}})$. Denote by $\mathcal{V}_{\acute{\mathbf{z}}}$, $\eta$, and $\phi:[0,\eta)\rightarrow \mathbb{R}_+$ the objects appearing in Definition \ref{KL_def}. Let $\rho>0$ be  such that $\mathcal{B}(\acute{\mathbf{z}},\rho)^{2mn}\subseteq \mathcal{V}_{\acute{\mathbf{z}}}$.
	Consider the sequence  $\{\mathbf{z}^\nu\triangleq (\mathbf{x}^\nu,\mathbf{y}^\nu)\}$   generated by the DOGT Algorithm   \eqref{NEXT_generic},  with initialization  $\mathbf{z}^0\triangleq(\mathbf{x}^0,\mathbf{y}^0)$; and  define $\acute{l}\triangleq L(\acute{\mathbf{z}})$ and $l^\nu\triangleq L(\mathbf{z}^\nu)$. Suppose that 
	\begin{equation}
	\acute{l}<l^\nu<\acute{l}+\eta,\quad \forall \nu\geq 0,
	\label{nonImprovingIterates_cond}\vspace{-0.2cm}
	\end{equation}
	and\vspace{-0.2cm}
	\begin{equation}
	K \, M\, \phi(l^0-\acute{l})+\norm{\mathbf{z}^{0}-\acute{\mathbf{z}}}<\rho,
	\label{proximity_cond}
	\end{equation}
	where \vspace{-0.2cm}\begin{equation} K= \sqrt{3}(1+L_c)\max\left(\frac{4nK_{||}^2}{1-\tilde{\rho}_R},\frac{K_{||}^2}{\varkappa(1-\tilde{\rho}_C)}\left(\alpha+\frac{2\sqrt{n}}{1+L_c}\right)^2, \alpha^2/\Gamma\right)^{1/2},
	\label{eq:def_K}
	\end{equation} 
	and $M>0$ is defined in \eqref{nablaL_bound_original} (cf.~Lemma \ref{nablaL_bound_lemma}). 
	
	Then,  $\{ \mathbf{z}^\nu\}$ satisfies:
	\begin{enumerate}[label=(\roman*)]
		\item $\mathbf{z}^\nu\in\mathcal{B}(\acute{\mathbf{z}},\rho)^{2mn}$, for all $\nu\geq 0$;
		\item $\sum_{t=k}^{\nu}\norm{\mathbf{z}^{t+1}-\mathbf{z}^{t}} \leq KM\left(\phi(l^k-\acute{l})-\phi(l^{\nu+1}-\acute{l})\right)$ for all $\nu,k\geq 0$ and $\nu\geq k$;
		\item $l^\nu\rightarrow\acute{l},\quad\mathrm{as}~\nu\rightarrow \infty$.
	\end{enumerate} 
\end{proposition} 

\begin{proof}
{ Throughout the proof, we will use the following shorthand  $d^\nu\triangleq d(\bx^\nu,\by^\nu)$. Let  $d^\nu>0$, for all integer $\nu\geq 0$; otherwise,  $\{\mathbf{x}^{\nu}\}$ converges in a finite number of steps, and its limit point is $\mathbf{x}^\infty=\mathbf{1}\theta^\infty$, for some $\theta^\infty\in\mathrm{crit}~F$. 
}

	We first bound  the ``length" $\sum_{t=k}^{\nu}\norm{\mathbf{z}^{t+1}-\mathbf{z}^{t}}$. By \eqref{NEXT_generic}, there holds
	\begin{equation*}
	\begin{aligned}
	\mathbf{x}^{\nu+1}-\mathbf{x}^{\nu}=&\left(\mathbf{R}-\mathbf{I}\right)\left(\mathbf{x}^\nu-\mathbf{1}\bar{x}^\nu\right)-\alpha\left(\mathbf{y}^\nu-\mathbf{c}\bar{y}^\nu\right)-\alpha\mathbf{c}\bar{y}^\nu,
	\\
	\mathbf{y}^{\nu+1}-\mathbf{y}^{\nu}=&\left(\mathbf{C}-\mathbf{I}\right)\left(\mathbf{y}^\nu-\mathbf{c}\bar{y}^\nu\right) +\nabla F_c(\mathbf{x}^{\nu+1})-\nabla F_c(\mathbf{x}^{\nu}).
	\end{aligned}
	\end{equation*}
	Using  $
	||\mathbf{A}||_2\leq \sqrt{n}||\mathbf{A}||_\infty$ and  $
	||\mathbf{A}||_2\leq \sqrt{n}||\mathbf{A}||_1$,  with $\mathbf{A}\in\mathcal{M}_n(\mathbb{R})$; and $||\mathbf{R}-\mathbf{I}||_\infty\leq 2 $ and $||\mathbf{C}-\mathbf{I}||_1\leq 2 $,
	we get\vspace{-0.1cm}
	\begin{equation*}
	\begin{aligned}
	\sum_{t=k}^{\nu}\norm{\mathbf{x}^{t+1}-\mathbf{x}^{t}}
	\leq & \sum_{t=k}^{\nu}2\sqrt{n}\norm{\mathbf{x}^t-\mathbf{1}\bar{x}^t}+\alpha\norm{\mathbf{y}^t-\mathbf{c}\bar{y}^t}+\alpha |\bar{y}^t|,\\
	\sum_{t=k}^{\nu}\norm{\mathbf{y}^{t+1}-\mathbf{y}^{t}}
	\leq & \sum_{t=k}^{\nu}2\sqrt{n}\norm{\mathbf{y}^t-\mathbf{c}\bar{y}^t}+L_c \sum_{t=k}^{\nu}\norm{\mathbf{x}^{t+1}-\mathbf{x}^{t}},
	\end{aligned}
	\end{equation*}
	where $L_c$ is the Lipschitz constant of $\nabla F_c$. The above inequalities imply\vspace{-0.2cm}
	\begin{equation}
	\begin{aligned}
	&\sum_{t=k}^{\nu}\norm{\mathbf{z}^{t+1}-\mathbf{z}^{t}}
	\\
	\leq  & \sum_{t=k}^{\nu}2(1+L_c)\sqrt{n}K_{||}\norm{\mathbf{x}^t-\mathbf{1}\bar{x}^t}_R+K_{||}\left(\alpha(1+L_c)+2\sqrt{n}\right)\norm{\mathbf{y}^t-\mathbf{c}\bar{y}^\nu}_C
	\\
	&+\alpha (1+L_c)|\bar{y}^t|
	\leq   K\sum_{t=k}^{\nu}d^t,
	\end{aligned}
	\label{length_bound}
	\end{equation}
	where 
	$K$  is defined in \eqref{eq:def_K}. 
	
	We prove now the proposition, starting from statement (ii). Multiplying both sides of  \eqref{L_descent} by $\phi'(l^\nu-\acute{l})$ and using   $\phi'(l^\nu-\acute{l})>0$ [due to property (iii) in Definition~\ref{KL_def} and \eqref{nonImprovingIterates_cond}] and   the concavity of $\phi$, yield
	\begin{equation}
	(d^\nu)^2~ \phi'(l^\nu-\acute{l})\leq \phi'(l^\nu-\acute{l})~ \left(l^{\nu}-l^{\nu+1}\right)\leq \phi(l^\nu-\acute{l})-\phi(l^{\nu+1}-\acute{l}).
	\label{deltaPHI_bound1}
	\end{equation}
	
	For all
	$\mathbf{z}\in\mathcal{V}_{\acute{\mathbf{z}}}\cap[ \acute{l} <L<\acute{l}+\eta]$, the {K\L} inequality \eqref{KL_ineq} holds; hence, assuming
	$\mathbf{z}^{t}\in\mathcal{B}(\acute{\mathbf{z}},\rho)^{2mn}$  for all $t=0,\ldots,\nu$, yields
	\begin{equation}
	\phi'(l^t-\acute{l})||\nabla L(\mathbf{z}^t)||\geq 1,\quad t=0,\ldots,\nu,
	\label{KL_ineq_restated}
	\end{equation}
	which together with \eqref{deltaPHI_bound1} and \eqref{nablaL_bound_original} (cf. Lemma \ref{nablaL_bound_lemma}), gives
	\begin{equation*}
	M\left(\phi(l^t-\acute{l})-\phi(l^{t+1}-\acute{l})\right)\geq d^t,\quad t=0,\ldots,\nu,
	\end{equation*}
	and thus  \vspace{-0.3cm}
	\begin{equation}
	M\left(\phi(l^k-\acute{l})-\phi(l^{\nu+1}-\acute{l})\right)\geq \sum_{t=k}^\nu d^t.
	\label{sum_d_nu_ineq}
	\end{equation}
	
	Combining \eqref{sum_d_nu_ineq} with \eqref{length_bound}, we obtain \begin{equation}
	\sum_{t=k}^{\nu}\norm{\mathbf{z}^{t+1}-\mathbf{z}^{t}}  \leq  KM\left(\phi(l^k-\acute{l})-\phi(l^{\nu+1}-\acute{l})\right).
	\label{length_bound_2}
	\end{equation}
	{Ineq. \eqref{length_bound_2} proves (ii) if  $\mathbf{z}^{\nu}\in\mathcal{B}(\acute{\mathbf{z}},\rho)^{2mn}$  for all $\nu\geq 0$, which is shown next.}
	
	Now let us prove statement (i). Letting $k=0$ in \eqref{length_bound_2}, by \eqref{proximity_cond}, we obtain\vspace{-0.2cm}
	\begin{equation*}
	\norm{\mathbf{z}^{\nu+1}-\acute{\mathbf{z}}}\leq KM\left(\phi(l^0-\acute{l})-\phi(l^{\nu+1}-\acute{l})\right)+\norm{\mathbf{z}^{0}-\acute{\mathbf{z}}}<\rho.
	\end{equation*}
	Therefore,  $\mathbf{z}^\nu\in\mathcal{B}(\acute{\mathbf{z}},\rho)^{2mn}$, for all $\nu\geq 0$. 
	
	We finally prove statement (iii).   Inequalities \eqref{nablaL_bound_original} (cf. Lemma \ref{nablaL_bound_lemma}) and \eqref{KL_ineq_restated} imply \vspace{-0.2cm}
	\begin{equation}
	\phi'(l^\nu-\acute{l})~ d^\nu\geq 1/M,\quad \nu\geq 0.
	\label{KL_plus_nablaLineq}
	\end{equation}
	On the other hand, by Lemma \ref{sequence_props_a}-(i), as $\nu\rightarrow\infty$, we have $l^\nu\rightarrow p$, for some $p\geq \acute{l}$. In fact, $p =\acute{l}$,  otherwise $p-\acute{l}>0$, which would   contradict  \eqref{KL_plus_nablaLineq} (because $d^\nu\rightarrow 0$ as $\nu\rightarrow\infty$ and $\phi'(p-\acute{l})<\infty$).
\end{proof}

Roughly speaking, Proposition~\ref{pre_conv_result} states that,  if the algorithm is initialized in a suitably chosen neighborhood of   a point at which $L$ satisfies the {K\L} property, then  it will converge to that  point. 
Combining this property with the subsequence convergence proved in Theorem \ref{main_global_conv_thm} we can obtain  global convergence of the sequence to critical points of $F$, as stated next.  




{
	\begin{theorem}
		\label{main_global_conv_thm}
		Consider the setting of Theorem~\ref{sequence_props}, and furthermore assume that $F$ is real-analytic.  Any sequence $\{(\mathbf{x}^\nu,\mathbf{y}^\nu)\}$ generated by the DOGT Algorithm \eqref{NEXT_generic} converges to some $(\mathbf{x}^\infty,\mathbf{y}^\infty)\in\mathrm{crit}~L$. Furthermore, $\mathbf{x}^\infty=\mathbf{1}\otimes\mathbf{\theta}^\infty$, for some $\mathbf{\theta}^\infty\in\mathrm{crit}~F$.\vspace{-0.1cm}
	\end{theorem} 
	\begin{proof} 
		Let  $\mathbf{z}^\infty\triangleq (\mathbf{x}^\infty,\mathbf{y}^\infty)$ be a limit point of $\{\mathbf{z}^\nu\triangleq (\mathbf{x}^{\nu},\mathbf{y}^{\nu})\}$. Since $\{l^\nu\triangleq L(\mathbf{z}^\nu)\}$ is convergent (cf. Lemma \ref{sequence_props_a}) and $L$ is continuous, we deduce  $l^\nu \to l^\infty\triangleq L(\mathbf{z}^\infty)$. Since $F$ is real-analytic, $L$ is real analytic (due to Lemma \ref{R_C_norms_plus_realAnalyticity} and  the fact that summation/composition of functions preserve real-analytic property  \cite[Prop. 2.2.8]{krantz2002primer}) and thus  {K\L} at  at $\mathbf{z}^\infty$ \cite{MR0160856}. Set $ \acute{\mathbf{z}}=\mathbf{z}^\infty$ and $\acute{l}=l^\infty$;   denote by $\mathcal{V}_{\acute{\mathbf{z}}}$, $\eta$, and $\phi:[0,\eta)\rightarrow \mathbb{R}_+$ the objects appearing in Definition \ref{KL_def}; and let    $\rho>0$ be  such that $\mathcal{B}(\acute{\mathbf{z}},\rho)^{2mn}\subseteq \mathcal{V}_{\acute{\mathbf{z}}}$. By the continuity of $\phi$ and the properties above, we deduce that there exists an integer $\nu_0$ such that i) $l^\nu\in (\acute{l}, \acute{l}+\eta)$, for all $\nu \geq \nu_0$;  and ii)  
		$
		K ~ M~ \phi(l^{\nu_0}-\acute{l})+\norm{\mathbf{z}^{\nu_0}-\acute{\mathbf{z}}}<\rho,
		$ with $K$ and $M$ defined in (\ref{eq:def_K}) and  \eqref{nablaL_bound_original}, respectively. Global convergence of the sequence $\{\mathbf{z}^\nu\}$ follows by applying Proposition~\ref{pre_conv_result} to the sequence $\{\mathbf{z}^{\nu+\nu_0}\}$. 
		
		Finally, by Lemma \ref{sequence_props_a}(ii),   $d(\bx^\nu,\by^\nu)\rightarrow 0$ as $\nu\rightarrow\infty$. Invoking the continuity of $\nabla L$ and  Lemma \ref{nablaL_bound_lemma}, we have $\nabla L(\mathbf{x}^\infty,\mathbf{y}^\infty)=\mathbf{0}$, thus   $(\mathbf{x}^\infty,\mathbf{y}^\infty)\in\mathrm{crit}~L$. By Theorem \ref{Th_asympt_convergence}(ii),     $\mathbf{x}^\infty=\mathbf{1}\otimes \mathbf{\theta}^\infty $, with $\mathbf{\theta}^\infty\in\mathrm{crit}~F$.  
	\end{proof}
}

In the following theorem, we provide some  convergence rate estimates.
{\begin{theorem}
	\label{Rate_theorem}
	In the setting of  Theorem \ref{main_global_conv_thm}, let   $L$ be a {K\L} function with $\phi(s)=cs^{1-\theta}$, for some constant $c>0$ and $\theta\in[0,1)$. Let   $\{\mathbf{z}^\nu\triangleq (\mathbf{x}^\nu,\mathbf{y}^\nu)\}$ be a  sequence generated by  DOGT Algorithm \eqref{NEXT_generic}. Then, there hold:
	\begin{enumerate}[label=(\roman*)]
		\item If $\theta=0$, $\{\mathbf{z}^\nu\}$ converges to $ {\mathbf{z}}^\infty$ in a finite number of iterations;
		\item If $\theta\in(0,1/2]$, then $||\mathbf{z}^\nu-{\mathbf{z}}^\infty||\leq C\tau^\nu$, for all $\nu\geq \bar{\nu}$ for some $\tau\in[0,1)$, $\bar{\nu}\in\mathbb{N}_+$, $C>0$;
		\item If $\theta\in(1/2,1)$, then $||\mathbf{z}^\nu- {\mathbf{z}}^\infty||\leq C\nu^{-\frac{1-\theta}{2\theta-1}}$, for all $\nu\geq \bar{\nu}$ for some $\bar{\nu}\in\mathbb{N}_+$, $C>0$.
	\end{enumerate} 
\end{theorem}}
\begin{proof}
	{ For sake of simplicity of notation, denote $d^\nu\triangleq d(\bx^\nu,\by^\nu)$ and define $D^\nu\triangleq \sum_{t=\nu}^\infty d^t$.} By   \eqref{length_bound},
	we have\vspace{-0.2cm}
	\begin{equation}
	\norm{\mathbf{z}^{\nu+1}-\mathbf{z}^\infty}\leq \sum_{t=\nu}^\infty\norm{\mathbf{z}^{t+1}-\mathbf{z}^{t}} {\leq}K D^\nu.\vspace{-0.2cm}
	\end{equation}
	It is then sufficient to establish the convergence rates for the sequence  $\{D^\nu\}$. 
	
	By  {K\L}  inequality \eqref{KL_ineq} and \eqref{nablaL_bound_original}, we have\vspace{-0.1cm}
	\begin{equation}
	Md^\nu\phi'(l^\nu-l^\infty)\geq 1 
	\implies  \tilde{M} (d^\nu)^{(1-\theta)/\theta}\geq   (l^\nu-l^\infty)^{1-\theta},\quad \forall\nu\geq \bar{\nu}
	\label{rate_bound1}
	\end{equation}
	for sufficiently large $\bar{\nu}$, where $\tilde{M}= \left(Mc(1-\theta)\right)^{(1-\theta)/\theta}$, $l^\nu\triangleq L(\mathbf{z}^\nu)$, and $l^\infty\triangleq L({\mathbf{z}}^\infty)$. In addition, by    \eqref{sum_d_nu_ineq} (setting $\acute{l}=l^\infty$), we have
	$D^\nu\leq M\phi(l^\nu-l^\infty)=Mc(l^\nu-l^\infty)^{1-\theta}$,
	which together with \eqref{rate_bound1}, yields\vspace{-0.2cm}
	\begin{equation}
	D^\nu\leq \tilde{M}Mc(d^\nu)^{(1-\theta)/\theta}=\tilde{M}Mc(D^{\nu}-D^{\nu+1})^{(1-\theta)/\theta},\quad \forall\nu\geq \bar{\nu}.
	\label{rate_bound3}
	\end{equation}
	The convergence rate estimates as stated in the theorem can be derived from   \eqref{rate_bound3}, using the same line of analysis introduced in \cite{Attouch2009}. The remaining part of the proof is provided in Appendix \ref{Rate_theorem_proof_appendix} for completeness.
\end{proof}

\subsection{Second-order guarantees}
\label{Guarantees_Sec}
We prove that  the   DOGT algorithm almost surely converges to   SoS solutions of \eqref{eq:P}, under a  suitably chosen  initialization and some additional conditions on the weight matrices $\mathbf{R}$ and $\mathbf{C}$. Following a path first established in \cite{pmlr-v49-lee16}  and further developed in \cite{Jordan_FOM_AvoidsSaddles}, the key to our argument for the non-convergence to strict saddle points of $F$ lies in formulating the DOGT algorithm as a dynamical system while leveraging an   instantiation of the stable manifold theorem, as given in \cite[Theorem 2]{Jordan_FOM_AvoidsSaddles}. 
The nontrivial task is finding a self-map representing DOGT so that the stable set of the strict saddles of $F$ is zero measure with respect to the domain of the mapping; note that the domain of the map--which is the set of initialization points--is not full dimensional and is the same as the support of the probability measure.


Our analysis is organized in the following three steps: 1) Sec.~\ref{sec:prelim_manifold}  introduces the preparatory   background; 
2) Sec.~\ref{DOGT_dynamicalSys} tailors the results of Step 1 to the DOGT algorithm; and 3) finally,   Sec. \ref{sec:SOS_main} states our main results about convergence of the DOGT algorithm to SoS solutions of \eqref{eq:P}. 

\subsubsection{The stable manifold theorem and unstable fixed-points}\label{sec:prelim_manifold}
{Let $g:\mathcal{S}\to \mathcal{S}$ be  a   mapping  from $\mathcal{S}$ to itself, 
	where $\mathcal{S}$ is a manifold without boundary.} 
Consider the dynamical system $\mathbf{u}^{\nu+1}=g(\mathbf{u}^\nu)$, with $\mathbf{u}^0\in \mathcal{S}$; we denote by $g^\nu$ the  {$\nu$-fold} composition of $g$. Our focus is on the analysis of the trajectories of the dynamical system 
around the fixed points of $g$; in particular we are interested in the  set of unstable fixed points of $g$. We begin introducing the following definition. 

\begin{definition}[Chapter 3 of \cite{Absil:2007:OAM:1557548}]
	\label{Differential_def}
	The differential of the mapping $g:\mathcal{S}\to \mathcal{S}$, denoted as $\mathrm{D}g(\mathbf{u})$, is a linear operator from $\mathcal{T}(\mathbf{u})\rightarrow\mathcal{T}(g(\mathbf{u}))$, where $\mathcal{T}(\mathbf{u})$ is the tangent space of $\mathcal{S}$ at   $\mathbf{u}\in \mathcal{S}$.   Given a curve $\gamma$ in $\mathcal{S}$ with $\gamma(0)=\mathbf{u}$ and $\frac{d\gamma}{d t}(0)=\mathbf{v}\in\mathcal{T}(\mathbf{u})$, the linear operator is defined as $\mathrm{D}g(\mathbf{u})\mathbf{v}=\frac{d(g\circ\gamma)}{d t}(0)\in \mathcal{T}(g(\mathbf{u}))$. The determinant of the linear operator $\det(\mathrm{D}g(\mathbf{u}))$ is the determinant of the matrix representing $\mathrm{D}g(\mathbf{u})$ with respect to a  standard basis.\footnote{This determinant   may not be uniquely defined, in the sense of being completely invariant to the basis used for the geometry. In this work,
		we are interested in properties of the determinant that are independent of scaling, and thus the potentially arbitrary choice of a standard basis does not affect our conclusions.}
\end{definition} 

We can now introduce the definition of the  set of unstable fixed points of $g$. 

\begin{definition}[Unstable fixed points]
	\label{UFP_def}
	The set of unstable fixed points of $g$ is defined as\vspace{-0.3cm}
	\begin{equation}\label{def_unstable_fixpoint}
	\mathcal{A}_g=\Big\{\mathbf{u}: g(\mathbf{u})=\mathbf{u},~\mathrm{spradii}\big(\mathrm{D} g(\mathbf{u})\big)>1\Big\}.
	\end{equation}
\end{definition}  
{The theorem below, which is based on the stable manifold theorem \cite[Theorem III.7]{shub1987global}, provides  tools to let us connect $\mathcal{A}_g$ with the set of limit points which    $\{\mathbf{u}^\nu\}$ can escape from.}\vspace{-0.1cm}

\begin{theorem}[{\cite[Theorem 2]{Jordan_FOM_AvoidsSaddles}}]
	\label{StableSetofUnstable_fixed_point_is_zeroM}
	Let $g:\mathcal{S}\rightarrow\mathcal{S}$ be a $\mathcal{C}^1$ mapping and\vspace{-0.2cm}
	\begin{equation*}
\det\left(\mathrm{D}g(\mathbf{u})\right)\neq 0,\quad \forall\mathbf{u}\in\mathcal{S}.\vspace{-0.2cm}
	\end{equation*}
	Then, the set of initial points that converge to an   unstable fixed point  (termed stable set of $\mathcal{A}_g$) is zero measure in $\mathcal{S}$.
	Therefore,  \vspace{-0.2cm}
	\begin{equation*}
	\mathbb{P}_{\mathbf{u}^0} \left( \lim_{\nu\to \infty}g^\nu(\mathbf{u}^0) \in \mathcal{A}_g \right)=0,\vspace{-0.2cm}
	\end{equation*} 
	where the probability is taken over the starting point $\mathbf{u}^0\in \mathcal{S}$.
\end{theorem}

\subsubsection{DOGT as a dynamical system}
\label{DOGT_dynamicalSys} 
Theorem \ref{StableSetofUnstable_fixed_point_is_zeroM} sets the path to the analysis of the convergence of the DOGT algorithm to SoS solutions of $F$: it is sufficient to describe the DOGT algorithm  by a proper mapping     $g:\mathcal{S}\to \mathcal{S}$ satisfying the assumptions in the theorem  and  such that the non-convergence of $g^\nu(\mathbf{u}^0)$, $\mathbf{u}^0\in \mathcal{S}$, to $\mathcal{A}_g$ implies the non-convergence of the DOGT algorithm to   strict saddles   of $F$.

We begin rewriting the DOGT in an equivalent and more convenient form. Define  $\mathbf{h}^\nu\triangleq \mathbf{y}^{\nu}-\nabla F_c(\mathbf{x}^{\nu})$;    \eqref{NEXT_generic} can be rewritten as  
\begin{equation}
\left\lbrace
\begin{aligned}
& \mathbf{x}^{\nu+1}=\mathbf{W}_R\mathbf{x}^{\nu}-\alpha\left(\mathbf{h}^\nu+\nabla F_c(\mathbf{x}^{\nu})\right);
\\
& \mathbf{h}^{\nu+1} =\mathbf{W}_C\mathbf{h}^{\nu}+\left(\mathbf{W}_C-\mathbf{I}\right)\nabla F_c(\mathbf{x}^{\nu}),
\end{aligned}
\right.
\label{NEXT_generic4}
\end{equation}
with arbitrary $\mathbf{x}^0\in\mathbb{R}^{nm}$ and $\mathbf{h}^0\in\text{span}(\mathbf{W}_C-\mathbf{I})$. By Theorem \ref{Th_asympt_convergence}, every limit point $(\mathbf{x}^\infty,\mathbf{h}^\infty)$ of $\{(\mathbf{x}^\nu,\mathbf{h}^\nu)\}$ has the form $\mathbf{x}^\infty=\mathbf{1}_n\otimes\boldsymbol{\theta}^\infty$ and $\mathbf{h}^\infty=-\nabla F_c(\mathbf{1}_n\otimes \boldsymbol{\theta}^\infty)$, for some $\boldsymbol{\theta}^\infty\in \mathrm{crit}~F$. We are interested in the non-convergence of (\ref{NEXT_generic4}) to such points whenever $\boldsymbol{\theta}^\infty\in \mathrm{crit}~F$ is a strict saddle of $F$. This motivates the following   definition. 



\begin{definition}[Consensual   strict saddle  points]
	\label{SSP_def}
	Let $\Theta^\ast_{ss}=\{\boldsymbol{\theta}^\ast\in \mathrm{crit}~F\,:\,\lambda_{\text{min}}(\nabla^2 F(\boldsymbol{\theta}^\ast))<0\}$ denote the set of strict saddles of $F$. The set of \emph{consensual strict saddle points} is defined  as
	\begin{equation}\label{def_U_g}
	\mathcal{U}^\ast\triangleq \left\{\left[\begin{array}{cc}
	\mathbf{1}_n\otimes \boldsymbol{\theta}^\ast\\ -\nabla F_c(\mathbf{1}_n\otimes \boldsymbol{\theta}^\ast)
	\end{array}
	\right]\,:\, \boldsymbol{\theta}^\star\in \Theta^\ast_{ss}\right\}.  	
	\end{equation}
	
\end{definition}

Roughly speaking,  $\mathcal{U}^\ast$ represents the candidate set of ``adversarial'' limit points which   any sequence generated by \eqref{NEXT_generic4} should   escape from.   The next step is then to write (\ref{NEXT_generic4}) as a proper dynamical system whose mapping satisfies conditions in Theorem  \ref{StableSetofUnstable_fixed_point_is_zeroM} and its  set of  unstable    fixed points $\mathcal{A}_g
$ is such that $\mathcal{U}^\ast \subseteq  \mathcal{A}_g$.\smallskip

\noindent  {\textbf{Identification of  $g$ and  $\mathcal{S}$.}}
Define   $\mathbf{u}\triangleq (\mathbf{x},\mathbf{h})$, where 
$\mathbf{x}\triangleq [\mathbf{x}_1^\top,\ldots, \mathbf{x}^\top_n]^\top$,
$\mathbf{h}\triangleq [\mathbf{h}_1^\top,\ldots, \mathbf{h}^\top_n]^\top$,
and  each $\mathbf{x}_i,\mathbf{h}_i\in \mathbb{R}^{m}$;  its   value at iteration $\nu$ is denoted by $\mathbf{u}^\nu\triangleq (\mathbf{x}^\nu,\mathbf{h}^\nu)$. Consider the dynamical system \vspace{-0.2cm}
\begin{equation}
\mathbf{u}^{\nu+1}=g(\mathbf{u}^{\nu}),\quad \text{with}\quad g\left(\mathbf{u}\right)\triangleq 
\begin{bmatrix}
\mathbf{W}_R\mathbf{x}-\alpha\nabla F_c\left(\mathbf{x}\right)-\alpha\mathbf{h}
\\
\mathbf{W}_C\mathbf{h}+\left(\mathbf{W}_C-\mathbf{I}\right)\nabla F_c\left(\mathbf{x}\right)
\end{bmatrix},
\label{NEXT_generic_dyn}
\end{equation}
and $\mathbf{u}^0\in\mathbb{R}^{nm}\times \text{span}(\mathbf{W}_C-\mathbf{I})$. The fixed-point iterate  \eqref{NEXT_generic_dyn} describes the trajectory generated by the DOGT algorithm (\ref{NEXT_generic4}).  
However, the   initialization imposed by    DOGT leads to a  $g$ that  maps $\mathbb{R}^{n m}\times \text{span}(\mathbf{W}_C-\mathbf{I})$ into $\mathbb{R}^{n m}\times \mathbb{R}^{n m}$. We show next how to unify the domain and codomain of  $g$ to a subspace $\mathcal{S}\subseteq \mathbb{R}^{nm}\times \mathbb{R}^{nm}$ as in form of the mapping in Theorem \ref{StableSetofUnstable_fixed_point_is_zeroM}.

Applying \eqref{NEXT_generic4} telescopically to    the update of  the $h$-variables yields: 
$\mathbf{h}^{\nu}=\mathbf{W}_C^{\nu}\mathbf{h}^0+\left(\mathbf{W}_C-\mathbf{I}\right)\mathbf{g}_{\text{acc}}^\nu,$ for all $\nu\geq 1$, where $\mathbf{g}_{\text{acc}}^\nu\triangleq \sum_{t=0}^{\nu-1}\mathbf{W}_C^t\nabla F_c\left(\mathbf{x}^{\nu-t-1}\right)$. Denoting $\bar{\mathbf{h}}^\nu\triangleq (\mathbf{1}_n^\top \otimes \mathbf{I}_m)  \mathbf{h}^\nu$, we have 
{\begin{equation}
	\bar{\mathbf{h}}^\nu=\cdots =\bar{\mathbf{h}}^0,\quad\text{and}\quad \mathbf{h}^\nu\in\mathbf{W}_C^{\nu}\mathbf{h}^0+\text{span}\left(\mathbf{W}_C-\mathbf{I}\right)\quad \forall \nu\geq 1.
	\label{h_embedded_eq}\vspace{-0.4cm}
	\end{equation}}

{
	The initialization $\mathbf{h}^0\in\text{span}\left(\mathbf{W}_C-\mathbf{I}\right)$ in \eqref{NEXT_generic4} 
	naturally suggests the following  $(2n-1)m$-dimensional linear subspace as candidate set  $\mathcal{S}$:
	\begin{equation}\label{def_S}
	\mathcal{S}\triangleq\mathbb{R}^{nm}\times \text{span}\left(\mathbf{W}_C-\mathbf{I}\right).
	\end{equation}
	Such an  $\mathcal{S}$ also ensures that  $g:\mathcal{S}\to \mathcal{S}$. In fact,   by (\ref{h_embedded_eq}),        $\mathbf{h}^\nu\in\text{span} (\mathbf{W}_C-\mathbf{I})$, for all $\nu\geq 1$, provided that    $\mathbf{h}^0\in\text{span}\left(\mathbf{W}_C-\mathbf{I}\right)$. Therefore,   $\{g^\nu(\mathbf{u}^0)\}\subseteq\mathcal{S}$, for all  $\mathbf{u}^0\in\mathcal{S}$. 	
}

Equipped with the mapping $g$ in (\ref{NEXT_generic_dyn}) and $\mathcal{S}$ defined in (\ref{def_S}),   we check next that   the condition  in Theorem \ref{StableSetofUnstable_fixed_point_is_zeroM}  is satisfied; we then prove that $\mathcal{U}^\ast \subseteq  \mathcal{A}_g$.

\textbf{1) ${g}$ is a diffeomorphism:}  To establish this property, we add the following  extra assumption on the weight matrices $\mathbf{R}$ and $\mathbf{C}$, which  is similar to Assumption \ref{nonsingular_W} for the DGD scheme.
\begin{assumption}
	\label{nonsingular_RC}
	Matrices $\mathbf{R}\in\mathcal{M}_n(\mathbb{R})$ and  $\mathbf{C}\in\mathcal{M}_n(\mathbb{R})$ are nonsingular. 
\end{assumption}

The above condition is not particularly restrictive and it is compatible with Assumption \ref{matrix_C_R}. A rule of thumb is to choose $\mathbf{R}=(\tilde{\mathbf{R}}+\mathbf{I})/2$ and $\mathbf{C}=(\tilde{\mathbf{C}}+\mathbf{I})/2$,  with   $\tilde{\mathbf{R}}$ and $\tilde{\mathbf{C}}$  satisfying Assumption \ref{matrix_C_R}. The new matrices still satisfy   Assumption \ref{matrix_C_R} due to the following fact:
given two nonnegative matrices $\mathbf{A},\mathbf{B}\in\mathcal{M}_n(\mathbb{R})$, if the directed graph associated with matrix $\mathbf{A}$ has a spanning tree  and $\mathbf{B}\geq \rho\mathbf{A}$, for some $\rho>0$, then the directed graph associated with matrix $\mathbf{B}$ has a spanning tree as well.

We build now the differential of $g$.  {Let $\tilde{g}$ be a smooth extension of \eqref{NEXT_generic_dyn} to $\mathbb{R}^{mn}\times \mathbb{R}^{mn}$, that is  $g=\tilde{g}|_\mathcal{S}$. The differential  $\mathrm{D}\tilde{g}(\mathbf{u})$ of $\tilde{g}$ at $\mathbf{u}\in \mathcal{S}$ reads 
	\begin{equation}
	\mathrm{D}\tilde{g}(\mathbf{u})= \begin{bmatrix}
	\mathbf{W}_R-\alpha\nabla^2 F_c(\mathbf{x}) &-\alpha\mathbf{I}
	\\
	\left(\mathbf{W}_C-\mathbf{I}\right)\nabla^2 F_c(\mathbf{x}) & \mathbf{W}_C
	\end{bmatrix};
	\label{Jacobian}
	\end{equation}	
	$\mathrm{D}\tilde{g}(\mathbf{u})$ is related to the differential of $g$ by $\mathrm{D}{g}(\mathbf{u})=\mathrm{D}\tilde{g}(\mathbf{u})\mathbf{P}_{\mathcal{T}(\mathbf{u})}$   \cite{ReimannianHess_Absil}, where $\mathbf{P}_{\mathcal{T}(\mathbf{u})}$ is the orthogonal projector  onto $\mathcal{T}(\mathbf{u})$. Using $\mathcal{T}(\mathbf{u})=\mathcal{S}$, for all $\mathbf{u}\in\mathcal{S}$ (recall that $\mathcal{S}$ is a linear subspace) and  denoting by   $\mathbf{U}_h\in{\mathbb{R}^{mn\times m(n-1)}}$ an orthonormal basis of $\mathrm{span}(\mathbf{W}_C-\mathbf{I})$, $\mathrm{D}{g}(\mathbf{u})$ reads\vspace{-0.1cm}
	\begin{equation}
	\label{Reimannian_Jacob}
	\mathrm{D}g(\mathbf{u})= \begin{bmatrix}
	\mathbf{W}_R-\alpha\nabla^2 F_c(\mathbf{x}) &-\alpha\mathbf{I}
	\\
	\left(\mathbf{W}_C-\mathbf{I}\right)\nabla^2 F_c(\mathbf{x}) & \mathbf{W}_C
	\end{bmatrix} {\mathbf{U}\mathbf{U}^\top},\quad\text{with}\quad
	\mathbf{U}\triangleq \begin{bmatrix}
	\mathbf{I} &\mathbf{0}
	\\
	\mathbf{0}&\mathbf{U}_h 
	\end{bmatrix}.\vspace{-0.1cm}
	\end{equation}}
Note that $\mathbf{P}_\mathcal{S}=\mathbf{U}\mathbf{U}^\top$.
We establish  next the conditions for   $g$ to be a $\mathcal{C}^1$ diffeomorphism, as stated in Theorem \ref{StableSetofUnstable_fixed_point_is_zeroM}.

\begin{proposition}
	\label{g_is_diffeomorphism}
	Consider the mapping $g:\mathcal{S}\to \mathcal{S}$ defined in (\ref{NEXT_generic_dyn}), under Assumptions \ref{P_assumption}-(i), \ref{matrix_C_R}, and \ref{nonsingular_RC}, with $\mathcal{S}$ defined in (\ref{def_S}).   If  the step-size is chosen according to  \vspace{-0.2cm} \begin{equation}0<\alpha<\frac{\sigma_\mathrm{min}(\mathbf{C}\mathbf{R})}{L_c},\label{eq:step_diff}\vspace{-0.1cm} \end{equation} where $L_c= L_\mathrm{max}$, then $\det\left(\mathrm{D}g(\mathbf{u})\right)\neq 0$, for all $\mathbf{u}\in\mathcal{S}$.
\end{proposition}

\begin{proof}	
	Since $\mathrm{D}g(\mathbf{u}):\mathcal{S}\rightarrow\mathcal{S}$,   it is sufficient to verify that $\mathrm{D}g(\mathbf{u})$ is an invertible linear transformation for every $\mathbf{u}\in \mathcal{S}$. Using the definition of $\mathbf{U}$, this is equivalent to show  that   $\mathbf{U}^T\mathrm{D}g(\mathbf{u})\mathbf{U}$ is invertible,  for all $\mathbf{u}\in\mathcal{S}$.
	Invoking \eqref{Reimannian_Jacob},   $\mathbf{U}^\top{\mathrm{D}}g(\mathbf{u})\mathbf{U}$ reads 
	\begin{equation}
	\begin{aligned}
	\mathbf{U}^\top{\mathrm{D}}g(\mathbf{u})\mathbf{U}= &\mathbf{U}^\top 
	\mathrm{D}\tilde{g}(\mathbf{u})	\mathbf{U}
	= \begin{bmatrix}
	\mathbf{W}_R-\alpha\nabla^2 F_c(\mathbf{x}) &-\alpha\mathbf{U}_h
	\\
	\mathbf{U}_h^\top\left(\mathbf{W}_C-\mathbf{I}\right)\nabla^2 F_c(\mathbf{x}) & \mathbf{U}_h^T\mathbf{W}_C\mathbf{U}_h
	\end{bmatrix}.
	\label{Projected_Jacobian0}
	\end{aligned}
	\end{equation}
	Since $\mathbf{U}_h^\top\mathbf{W}_C\mathbf{U}_h$ is non-singular, we can use the   Schur complement  of $\mathbf{U}^\top{\mathrm{D}}g(\mathbf{u})\mathbf{U}$ with respect to $\mathbf{U}_h^\top\mathbf{W}_C\mathbf{U}_h$ and write   
	\begin{equation}
	\begin{aligned}
	\mathbf{U}^\top{\mathrm{D}}g(\mathbf{u})\mathbf{U}=\mathbf{S}_1\begin{bmatrix}
	\mathbf{W}_R-\alpha\nabla^2 F_c(\mathbf{x})+\alpha \boldsymbol{\Phi}\left(\mathbf{W}_C-\mathbf{I}\right)\nabla^2 F_c(\mathbf{x})  &\mathbf{0}
	\\
	\mathbf{0}& \mathbf{U}_h^\top\mathbf{W}_C\mathbf{U}_h
	\end{bmatrix}\mathbf{S}_2,
	\end{aligned}
	\label{Projected_Jacobian2}
	\end{equation}
	where  $\boldsymbol{\Phi}\triangleq \mathbf{U}_h\left(\mathbf{U}_h^\top\mathbf{W}_C\mathbf{U}_h\right)^{-1}\mathbf{U}_h^\top$, and  $\mathbf{S}_1$ and $\mathbf{S}_2$ are some nonsingular matrices.
	By   \eqref{Projected_Jacobian2}, it is sufficient to  show that 
	\begin{equation}
	\begin{aligned}
	\mathbf{S} \triangleq \,&\mathbf{W}_R-\alpha\nabla^2 F_c(\mathbf{x})+\alpha \boldsymbol{\Phi}\left(\mathbf{W}_C-\mathbf{I}\right)\nabla^2 F_c(\mathbf{x})
	\\
	=&\mathbf{W}_R-\alpha\mathbf{W}_C^{-1}\nabla^2 F_c(\mathbf{x})+\alpha\left(\boldsymbol{\Phi}-\mathbf{W}_C^{-1}\right)\left(\mathbf{W}_C-\mathbf{I}\right)\nabla^2 F_c(\mathbf{x})
	\end{aligned}
	\label{Eqq}
	\end{equation}
	is non-singular.
	Using  
	$\mathbf{W}_C-\mathbf{I}=\mathbf{U}_h\boldsymbol{\Delta}$, for some 
	$\boldsymbol{\Delta}\in\mathbb{R}^{m(n-1)\times mn}$ (recall that $\mathbf{U}_h$ is an orthonormal basis of $\mathrm{span}(\mathbf{W}_C-\mathbf{I})$), we can write 
	\begin{equation}
	\label{Phi_expression}
	\begin{aligned}
	\boldsymbol{\Phi} =\,&\mathbf{U}_h\left(\mathbf{U}_h^\top\mathbf{W}_C\mathbf{U}_h\right)^{-1}\mathbf{U}_h^\top
	= \mathbf{U}_h \left(\mathbf{I}+\boldsymbol{\Delta}\mathbf{U}_h\right)^{-1}\mathbf{U}_h^\top
	\\
	\overset{(a)}{=} \,& \mathbf{U}_h\mathbf{U}_h^\top-\mathbf{U}_h\boldsymbol{\Delta}\left(\mathbf{I}+\mathbf{U}_h\boldsymbol{\Delta}\right)^{-1}\mathbf{U}_h\mathbf{U}_h^\top
	\\
	= \,& \mathbf{U}_h\mathbf{U}_h^\top-\left(\mathbf{W}_C-\mathbf{I}\right)\mathbf{W}_C^{-1}\mathbf{U}_h\mathbf{U}_h^\top
	\\
	=\, &\mathbf{W}_C^{-1}\mathbf{U}_h\mathbf{U}_h^\top,
	\end{aligned}
	\end{equation}
	where (a) we used the Woodbury identity of inverse matrices. Using \eqref{Phi_expression} in \eqref{Eqq}, we obtain\vspace{-0.2cm}
	\begin{equation*}
	\begin{aligned}
	\mathbf{S}
	=\,&\mathbf{W}_R-\alpha\mathbf{W}_C^{-1}\nabla^2 F_c(\mathbf{x})-\alpha\mathbf{W}_{C}^{-1}\underset{=\mathbf{{0}}}{\underbrace{\left(\mathbf{I}-\mathbf{U}_{h}\mathbf{U}_{h}^{\top}\right)\left(\mathbf{W}_{C}-\mathbf{I}\right)}}\nabla^2 F_c(\mathbf{x}).
	\end{aligned}
	\label{Eqq2}\vspace{-0.2cm}
	\end{equation*}	
	Therefore,   if $\alpha<\frac{\sigma_{\min}(\mathbf{CR})}{L_c}$, $\mathbf{S}$ is invertible, and consequently,  so is $\mathbf{U}^\top{\mathrm{D}}g(\mathbf{u})\mathbf{U}$.
\end{proof}

\textbf{2) The consensual strict saddle points are unstable fixed points of $g$ ($\mathcal{U}^\ast \subseteq \mathcal{A}_g$): }   First of all, note that every  limit point of the sequence generated by   \eqref{NEXT_generic4}  is a fixed point of $g$ on $\mathcal{S}$; the converse might not be true. The next result  establishes the desired connection between the set $\mathcal{A}_g$ of unstable fixed points of $g$ (cf.~Definition \ref{UFP_def}) and  the set $\mathcal{U}^\ast$ of consensual strict saddle points (cf.~Definition \ref{SSP_def}).  This will let us   infer the instability of $\mathcal{U}^\ast$ from that of $\mathcal{A}_g$.

\begin{proposition}
	\label{CSS_is_in_unstable_set}
	{ Suppose that Assumptions \ref{P_assumption}-(i) and  \ref{matrix_C_R} hold along with  one of the following two conditions }
	\begin{itemize}
		\item[(i)] The weight matrices $\mathbf{R}$ and $\mathbf{C}$ are   symmetric;  
		\item[(ii)] $m=1$.
	\end{itemize}
	Then,  any consensual strict saddle point  is an unstable fixed point of $g$, i.e.,
	\begin{equation}
	\mathcal{U}^\ast\subseteq\mathcal{A}_g,
	\end{equation}
	with $\mathcal{A}_g$ and $\mathcal{U}^\ast$    defined in    \eqref{def_unstable_fixpoint} and \eqref{def_U_g}, respectively. 
\end{proposition}
\begin{proof}		
	Let $\mathbf{u}^\ast\in\mathcal{U}^\ast$; $\mathbf{u}^\ast$ is a fixed point of   $g$ defined in \eqref{NEXT_generic_dyn}. It is thus sufficient to show that    $\mathrm{D}g(\mathbf{u}^\ast)$    has an eigenvalue with magnitude   greater than one. 
	
	To do so, we begin   showing that the differential $\mathrm{D}\tilde{g}(\mathbf{u}^\ast)$ of $\tilde{g}$ at $\mathbf{u}^\ast$ has an eigenvalue greater than one. Using \eqref{Jacobian}, 
	$\mathrm{D}\tilde{g}(\mathbf{u}^\ast)$  reads 
	\begin{equation}
	\mathrm{D}\tilde{g}(\mathbf{u}^\ast)= \begin{bmatrix}
	\mathbf{W}_R-\alpha\nabla^2 F_c^\ast &-\alpha\mathbf{I}
	\\
	\left(\mathbf{W}_C-\mathbf{I}\right)\nabla^2 F_c^\star & \mathbf{W}_C
	\end{bmatrix},
	\label{Jacobian_star}
	\end{equation}
	where we defined the shorthand   $\nabla^2 F_c^\ast\triangleq \nabla^2 F_c\left(\mathbf{1}\otimes \boldsymbol{\theta}^\ast\right)$, and  $\boldsymbol{\theta}^\ast\in\Theta^\ast_{ss}$.  
	We need to prove \vspace{-0.2cm}
	\begin{equation}
	\det\left(\mathrm{D}\tilde{g}(\mathbf{u}^\ast)-\lambda_u\mathbf{I}\right)=0,\quad\text{ for some }\quad |\lambda_u|>1.
	\label{det_eq1}
	\end{equation}
	If $|\lambda_u|>1$,   $\mathbf{W}_C-\lambda_u\mathbf{I}$ is nonsingular (since $\mathrm{spradii}(\mathbf{C})=1$). Using the   Schur complement of   $\mathrm{D}\tilde{g}(\mathbf{u}^\ast)-\lambda_u\mathbf{I}$ with respect to $\mathbf{W}_C-\lambda_u\mathbf{I}$, we have  
	\begin{equation}
	\begin{aligned}
	\mathrm{D}\tilde{g}(\mathbf{u}^\ast)-\lambda_u\mathbf{I}=	\tilde{\mathbf{S}}_1\begin{bmatrix}
	\left(\mathrm{D}\tilde{g}(\mathbf{u}^\ast)-\lambda_u\mathbf{I}\right)/\left(\mathbf{W}_C-\lambda_u\mathbf{I}\right)&\mathbf{0}
	\\
	\mathbf{0} & \mathbf{W}_C-\lambda_u\mathbf{I}
	\end{bmatrix}
	\tilde{\mathbf{S}}_2,
	\end{aligned}
	\label{NEXT_generic_dyn_Jacobian_Decomposition}
	\end{equation}
	for some    $\tilde{\mathbf{S}}_1,\tilde{\mathbf{S}}_2\in\mathcal{M}_{2mn}(\mathbb{R})$, with  $\det(\tilde{\mathbf{S}}_1)=\det(\tilde{\mathbf{S}}_2)=1$. 
	Given \eqref{NEXT_generic_dyn_Jacobian_Decomposition},  \eqref{det_eq1} holds if and only if
	\begin{equation*}
	\det
	\begin{bmatrix}
	\mathbf{W}_R-\lambda_u\mathbf{I}-\alpha\nabla^2 F_c^\star+\alpha \left(\mathbf{W}_C-\lambda_u\mathbf{I}\right)^{-1} \left(\mathbf{W}_C-\mathbf{I}\right) \nabla^2 F_c^\ast  & \mathbf{0}
	\\
	\mathbf{0} & \mathbf{W}_C-\lambda_u\mathbf{I}
	\end{bmatrix}=0,
	\end{equation*}
	or equivalently \vspace{-0.2cm}
	\begin{equation}
	\det
	\left(
	\mathbf{W}_R-\lambda_u\mathbf{I}-\alpha\nabla^2 F_c^\ast+\alpha \left(\mathbf{W}_C-\lambda_u\mathbf{I}\right)^{-1} \left(\mathbf{W}_C-\mathbf{I}\right) \nabla^2 F_c^\ast\right)=0.
	\label{det_eq3}
	\end{equation}
	Multiplying both sides of  \eqref{det_eq3} by $\det(\mathbf{W}_C-\lambda_u\mathbf{I})$ yields   \vspace{-0.2cm}
	\begin{equation}
	Q(\lambda_u)\triangleq \det
	\bigg(
	\underbrace{\left(\mathbf{W}_C-\lambda_u\mathbf{I}\right)\left(\mathbf{W}_R-\lambda_u\mathbf{I}\right)+\alpha(\lambda_u-1)\nabla^2 F_c^\ast}_{\triangleq\mathbf{T}(\lambda_u)}\bigg)=0.
	\label{det_eq4}
	\end{equation}
	Trivially  $Q(\lambda_u)>0$, if $\lambda_u\gg 1$. Therefore,  to show that \eqref{det_eq1} holds, it is sufficient  to prove that there exists  some $\lambda_u>1$ such that $Q(\lambda_u)\leq 0$. Next, we prove this result  under either   condition (i) or (ii). 
	
	Suppose (i) holds; $\mathbf{R}$ and $\mathbf{C}$ are symmetric.  Define $\tilde{\boldsymbol{\upsilon}}\triangleq \mathbf{1}\otimes \boldsymbol{\upsilon}$, where $\boldsymbol{\upsilon}$ is the unitary eigenvector associated with  {a negative eigenvalue of $\nabla^2 F(\boldsymbol{\theta}^\ast)$, and let  $\lambda_\mathrm{min}(\nabla^2 F(\boldsymbol{\theta}^\ast))=-\delta$;} we can write \vspace{-0.2cm}
	\begin{equation}
	\tilde{\boldsymbol{\upsilon}}^\top\mathbf{T}(\lambda_u) \tilde{\boldsymbol{\upsilon}}=n(\lambda_u-1)\left(\lambda_u-1- \alpha\delta/n\right)<0,
	\end{equation}
	for all $1<\lambda_u<1+\alpha\delta/n$. By Rayleigh-Ritz theorem, $\mathbf{T}(\lambda_u)$ has a negative eigenvalue, implying that there exists  some real value $\bar{\lambda}_u>1$ such that $Q(\bar{\lambda}_u)=0$.

	Suppose now that  conditions (ii) holds;  $\mathbf{W}_R$ and  $\mathbf{W}_C$ reduce to  $\mathbf{R}$ and $\mathbf{C}$, respectively. Note that  $\mathbf{R}$ and $\mathbf{C}$ are now not symmetric. Let $\lambda_u=1+\epsilon$, and consider the Taylor expansion of \vspace{-0.2cm}
	\begin{equation}
	Q(1+\epsilon)=\det
	\bigg(
	\left(\mathbf{C}-\mathbf{I}\right)\left(\mathbf{R}-\mathbf{I}\right)+\epsilon\left(\alpha\nabla^2 F_c^\ast+2\mathbf{I}-\mathbf{C}-\mathbf{R}\right)+\epsilon^2\mathbf{I}\bigg),\vspace{-0.2cm}
	\end{equation} 
	around $\epsilon=0$. Define $\mathbf{M}\triangleq \left(\mathbf{C}-\mathbf{I}\right)\left(\mathbf{R}-\mathbf{I}\right)$ and $\mathbf{N}\triangleq \alpha\nabla^2 F_c^\ast+2\mathbf{I}-\mathbf{C}-\mathbf{R}$.
	It is clear that $Q(1)=0$; then, by the Jacobi's formula, we have  \vspace{-0.2cm}
	\begin{equation}
	Q(1+\epsilon)=\text{tr}\Big(\text{adj}\left(\mathbf{M}\right) \mathbf{N}\Big)\epsilon+O(\epsilon^2).\vspace{-0.2cm}
	\label{Q_function_decomposition1}
	\end{equation}
	Expanding \eqref{Q_function_decomposition1} yields \vspace{-0.2cm}
	\begin{equation}
	\begin{aligned}
	Q(1+\epsilon)= & \text{tr}\Big(\text{adj}\left(\mathbf{R}-\mathbf{I}\right)\text{adj}\left(\mathbf{C}-\mathbf{I}\right) \mathbf{N}\Big)\epsilon+O(\epsilon^2)
	\\
	=&\text{tr}\Big(\mathbf{1}\tilde{\mathbf{r}}^\top\tilde{\mathbf{c}}\mathbf{1}^\top \mathbf{N}\Big)\epsilon+O(\epsilon^2)
	=(\tilde{\mathbf{r}}^\top\tilde{\mathbf{c}})\mathbf{1}^\top \mathbf{N}\mathbf{1}\epsilon+O(\epsilon^2),
	\end{aligned}
	\label{Q_function_decomposition2}\vspace{-0.2cm}
	\end{equation}
	where $\tilde{\mathbf{r}}$ and $\tilde{\mathbf{c}}$ are the Perron vectors of $\mathbf{R}$ and $\mathbf{C}$, respectively. The second equality in \eqref{Q_function_decomposition2} is due to the following fact:  a  rank-$(n-1)$ matrix $\mathbf{A}\in\mathcal{M}_n(\mathbb{R})$ has rank-$1$ adjugate matrix $\text{adj}\left(\mathbf{A}\right)=\mathbf{a}\mathbf{b}^\top$, where $\mathbf{a}$ and $\mathbf{b}$ are   non-zero vectors belonging to the 1-dimensional null space of  $\mathbf{A}$ and $\mathbf{A}^\top$, respectively    \cite[Sec. 0.8.2]{Horn-Johnson_book}.
	We also have $\tilde{\zeta}\triangleq \tilde{\mathbf{r}}^\top\tilde{\mathbf{c}}>0$, due to Lemma \ref{R_C_norms_plus_realAnalyticity}. 
	Furthermore, since  $\boldsymbol{\theta}^\ast\in\Theta_{ss}^\ast$,    $\mathbf{1}^\top \nabla^2F_c^\ast\mathbf{1}\leq -\delta$, for some $\delta>0$, and \vspace{-0.2cm}
	\begin{equation} 
	Q(1+\epsilon)\leq -\delta\tilde{\zeta}\alpha\epsilon+O(\epsilon^2),\vspace{-0.2cm}
	\end{equation}	
	which implies the existence of a sufficiently small $\epsilon>0$ such that $Q(1+\epsilon)<0$. Consequently, there must exist some $\bar{\lambda}_u>1$ such that \eqref{det_eq1} holds. Moreover, such $\bar{\lambda}_u$ is a real eigenvalue of $\mathrm{D}\tilde{g}(\mathbf{u}^\ast)$.
	
	To summarize, we proved that there exists an eigenpair $(\bar{\lambda}_u,\mathbf{v}_u)$ of $\mathrm{D}\tilde{g}(\mathbf{u}^\ast)$, with   $\bar{\lambda}_u>1$.	 
	Next we show that $(\bar{\lambda}_u,\mathbf{v}_u)$ is also  an eigenpair of $\mathrm{D}{g}(\mathbf{u}^\ast)$. Let us partition  $\mathbf{v}_u\triangleq (\mathbf{v}^x_u,\mathbf{v}^h_u)$ such that\vspace{-0.2cm}
	\begin{equation}
	\begin{bmatrix}
	\mathbf{W}_R-\alpha \nabla^2 F_c\left(\mathbf{x}^\ast\right) &-\alpha\mathbf{I}
	\\
	\left(\mathbf{W}_C-\mathbf{I}\right)  \nabla^2 F_c\left(\mathbf{x}^\ast\right) & \mathbf{W}_C
	\end{bmatrix}
	\begin{bmatrix}
	\mathbf{v}^x_u
	\\
	\mathbf{v}^h_u
	\end{bmatrix}
	=
	\bar{\lambda}_u\begin{bmatrix}
	\mathbf{v}^x_u
	\\
	\mathbf{v}^h_u
	\end{bmatrix}.	\end{equation}
	In particular, we have 
	$\left(\mathbf{W}_C-\mathbf{I}\right)\left( \nabla^2 F_c\left(\mathbf{x}^\ast\right)\mathbf{v}^x_u+\mathbf{v}^h_u\right)=(\bar{\lambda}_u-1)\mathbf{v}^h_u$,
	which implies  $\mathbf{v}_u^h\in\text{span}(\mathbf{W}_C-\mathbf{I})$, 	since $\bar{\lambda}_u-1\neq 0$. Therefore,   $\mathbf{v}_u\in \mathcal{S}$.
	
	Now, let $\mathbf{P}_\mathcal{S}$ be the orthogonal projection matrix onto $\mathcal{S}$. Since  $\mathbf{v}_u\in \mathcal{S}$, we have
	\begin{equation}
	\label{eigenpair_equivalence}
	\mathrm{D}\tilde{g}(\mathbf{u}^\ast)\mathbf{v}_u=\bar{\lambda}_u\mathbf{v}_u \implies \mathrm{D}\tilde{g}(\mathbf{u}^\ast)\mathbf{P}_\mathcal{S}^\top\mathbf{v}_u=\bar{\lambda}_u\mathbf{v}_u \overset{(a)}{\implies}
	\mathrm{D}{g}(\mathbf{u}^\ast)\mathbf{v}_u=\bar{\lambda}_u\mathbf{v}_u ,
	\end{equation}
	where (a) is due to  $\mathrm{D}{g}(\mathbf{u}^\ast)=\mathrm{D}\tilde{g}(\mathbf{u}^\ast)\mathbf{P}_\mathcal{S}^\top$ [cf. \eqref{Reimannian_Jacob}]. Hence $(\bar{\lambda}_u,\mathbf{v}_u)$ is also an eigenpair of $\mathrm{D}{g}(\mathbf{u}^\ast)$, which completes the proof.
\end{proof}

\begin{remark}
	Note that  condition (i) in Proposition \ref{CSS_is_in_unstable_set} implies that  $\mathcal{G}_C$ and $\mathcal{G}_R$ are undirected graphs. Condition (ii) extends the network model to directed topologies  under assumption $m=1$.  {For sake of completeness, we relax condition (ii) in Appendix \ref{CSS_is_in_unstable_set_extension_Appendix} to arbitrary $m\in\mathbb{N}$, under extra (albeit  mild) assumptions  on the set of strict saddle points and the weight matrices $\mathbf{R}$ and $\mathbf{C}$.} 
\end{remark}

\subsubsection{DOGT likely converges to SoS solutions of \eqref{eq:P}}
\label{sec:SOS_main} 
Combining Theorem \ref{StableSetofUnstable_fixed_point_is_zeroM}, 
Proposition \ref{g_is_diffeomorphism}, and  Proposition \ref{CSS_is_in_unstable_set}, we can readily obtain the following second-order guarantees of the DOGT algorithms.   

{\begin{theorem}
		\label{Main_SOS_guarantee_long} Consider Problem \eqref{eq:P}, under Assumptions  \ref{P_assumption} and \ref{Net-Assump}; and let $\{\mathbf{u}^\nu\triangleq (\bx^\nu, \mathbf{h}^\nu)\}$ be the sequence generated by the DOGT Algorithm (\ref{NEXT_generic4})   under the following tuning: i) the step-size $\alpha$  satisfies \eqref{alpha_conds_assymptotic_conv} [or \eqref{alpha_conds_assymptotic_conv_final}]  and \eqref{eq:step_diff}; the weight matrices $\mathbf{C}$ and $\mathbf{R}$ are chosen according to Assumptions \ref{matrix_C_R} and \ref{nonsingular_RC}; and the initialization is set to $\mathbf{u}^0\in \mathcal{S}$, with $\mathcal{S}$ defined in (\ref{def_S}). Furthermore, suppose that either (i) or (ii) in Proposition \ref{CSS_is_in_unstable_set} holds.  Then, we have \vspace{-0.1cm}
		\begin{equation}
		\mathbb{P}_{\mathbf{u}^0}\left(\lim_{\nu\rightarrow\infty}\mathbf{u}^\nu\in\mathcal{U}^\ast\right)=0,\vspace{-0.1cm}
		\label{zero_probability_to_ssp}
		\end{equation}
		where the probability is taken over $\mathbf{u}^0\in \mathcal{S}$. 
		
		{In addition, if} $F$ is a {K\L} function, then  $\{\mathbf{x}^\nu\}$   converges almost surely to $\mathbf{1}\otimes \boldsymbol{\theta}^\infty$ at a rate determined in Theorem \ref{Rate_theorem}, where ${\boldsymbol{\theta}}^\infty$ is a SoS solution of \eqref{eq:P}.
\end{theorem}}

{Note that \eqref{zero_probability_to_ssp} implies  the desired second-order guarantees only when the sequence $\{\mathbf{u}^\nu\}$ convergences [i.e., the limit in \eqref{zero_probability_to_ssp} exists]; otherwise   \eqref{zero_probability_to_ssp} is trivially satisfied, and some limit point of  $\{\mathbf{u}^\nu\}$ can  belong  to $\mathcal{U}^\ast$ with non-zero probability. A sufficient condition for the required global convergence of $\{\mathbf{u}^\nu\}$   is that $F$ is a {K\L}  function, which is stated in the second part of the above theorem.}

{
	\begin{remark}[Comparison with \cite{Mingyi_SOAlgs_AvoidsSaddles}]\label{DOGT_vs_ProxPDA_remark}
		As already discussed in Sec. \ref{sec:intro_distributed}, the primal-dual methods in  \cite{Mingyi_SOAlgs_AvoidsSaddles}  is applicable to \eqref{eq:P}; it  is proved to almost surely converge to SoS solutions.  
		Convergence of \cite{Mingyi_SOAlgs_AvoidsSaddles} is proved under stricter conditions on the problem than DOGT, namely: i) the network must be {\it  undirected}; and  ii) the Hessian of each local $f_i$ must be  Lipschitz continuous. It does \emph{not} seem possible to extend the analysis of \cite{Mingyi_SOAlgs_AvoidsSaddles}  beyond this assumptions. 
	\end{remark}
}

{
	\section{Numerical Results} 
	\label{numerical_experiments}
	In this section we test the behavior of  DGD and DOGT   around strict saddles  on  three  classes of nonconvex   problems, namely: i) a quadratic function (cf. Sec. \ref{sec:sim-quadratic}); ii) a classification problem based on  the cross-entropy risk function using sigmoid activation functions (cf. Sec. \ref{BLR_simulations}); and iii) a two Gaussian mixture model (cf. Sec. \ref{GMM_simulations}).}

\subsection{Nonconvex quadratic optimization}\label{sec:sim-quadratic}
Consider 
\vspace{-0.3cm}
\begin{equation}
\min_{\boldsymbol{\theta}\in \mathbb{R}^{m}}\, F\left(\boldsymbol{\theta}\right)= \frac{1}{2}\sum_{i=1}^{n} \left(\boldsymbol{\theta}-\mathbf{b}_i\right)^\top\mathbf{Q}_i\left(\boldsymbol{\theta}-\mathbf{b}_i\right),\vspace{-0.2cm}
\label{eq:P_test}
\end{equation}
where $m=20$; $n=10$; $\mathbf{b}_i$'s are i.i.d Gaussian zero mean random vectors with standard deviation $10^3$; and $\mathbf{Q}_i$'s are $m\times m$ randomly generated symmetric matrices where $\sum_{i=1}^n\mathbf{Q}_i$ has $m-1$ eigenvalues $\{\lambda_i\}_{i=1}^{m-1}$ uniformly distributed over $(0,n]$, and one negative eigenvalue $\lambda_m=-n\delta$, with $\delta=0.01$. Clearly \eqref{eq:P_test} is an instance of Problem \eqref{eq:P}, with   
$F$ having a unique strict saddle point   $\boldsymbol{\theta}^{\ast}=\left(\sum_i\mathbf{Q}_i\right)^{-1}\sum_i\mathbf{Q}_i\mathbf{b_i}$. The network of $n$ agents is  modeled as  a ring; the weight matrix $\mathbf{W}\triangleq\{w_{ij}\}_{i,j=1}^{n}$, compliant to the graph topology,  is generated to be  doubly stochastic. 

To test the escaping properties of DGD and DOGT from the strict saddle of $F$,  we initialize the  algorithms  in a randomly generated   neighborhood of $\boldsymbol{\theta}^{\ast}$. More specifically,    every agent's initial point  is  $\mathbf{x}_i^0=\boldsymbol{\theta}^{\ast}+\boldsymbol{\epsilon}_{x,i}$, $i\in[n]$. In addition, for the DOGT algorithm,    we set  $\mathbf{y}^0_i=\nabla f_i(\mathbf{x}_i^0)+(w_{ii}-1)\boldsymbol{\epsilon}_{y,i}+\sum_{j\neq i}w_{ij}\boldsymbol{\epsilon}_{y,j}$, where  $\boldsymbol{\epsilon}_{x,i}$'s and $\boldsymbol{\epsilon}_{y,i}$'s are realizations of  i.i.d. Gaussian random vectors with standard deviation equal to 1.  {Both algorithms use the same  step-size $\alpha=0.99\,\sigma_\mathrm{min}(\mathbf{I}+\mathbf{W})/L_c$, with $L_c=\max_i\{|\lambda_i|\}$; this is the largest theoretical step-size  guaranteeing   convergence of  the DGD algorithm (cf. Theorem \ref{DGD_conv_cons_thm}).}

\begin{figure}[t]
	\center
	\includegraphics[scale=0.22]{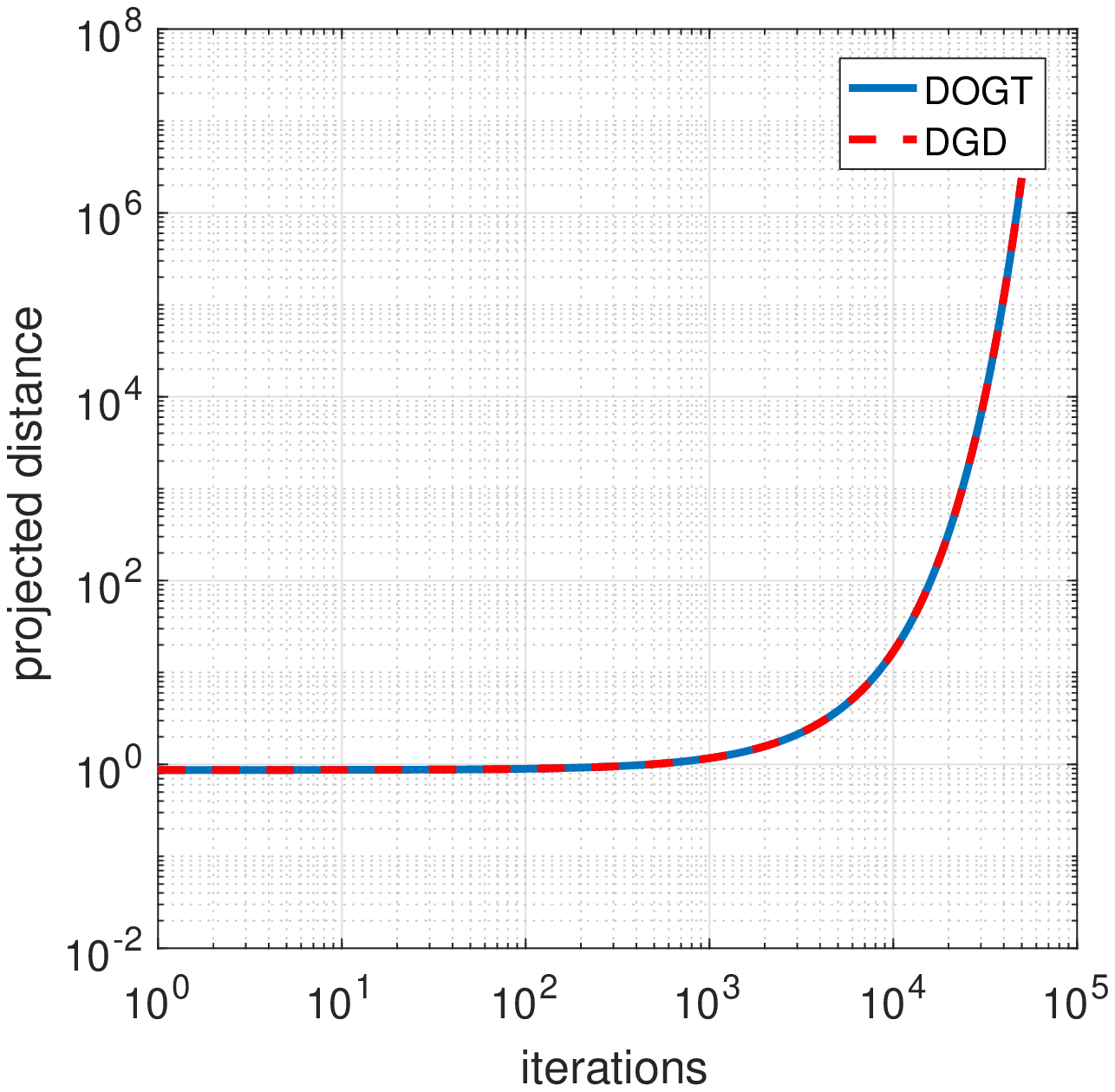}\includegraphics[scale=0.22]{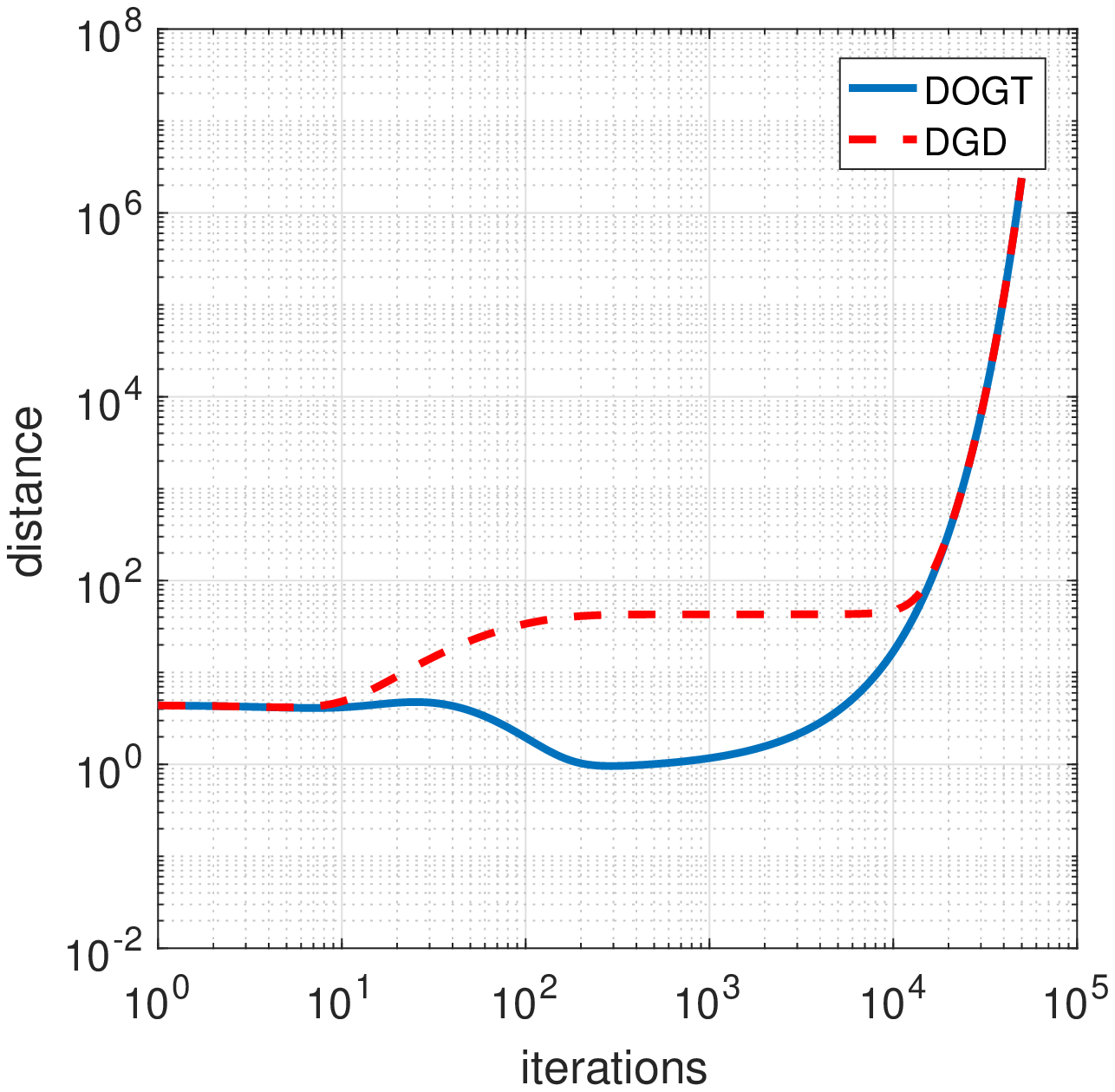}\vspace{-0.2cm}
	\caption{Escaping properties of  DGD and DOGT, applied to Problem (\ref{eq:P_test}). Left plot: distance of the average iterates from $\boldsymbol{\theta}^{\ast}$ projected onto the unstable manifold  ${E}_u$ versus the number of iterations. Right plot:    distance of the average iterates from $\boldsymbol{\theta}^{\ast}$ versus the number of iterations.  }
	\label{DGD_DOGT_escaping_Figs}\vspace{-0.3cm}
\end{figure}

In the left panel of  Fig. \ref{DGD_DOGT_escaping_Figs}, we plot the distance of the average iterates $\bar{\mathbf{x}}^\nu=(1/n)\sum_{i=1}^n\mathbf{x}_i^\nu$ from the critical point $\boldsymbol{\theta}^{\ast}$ projected on the unstable manifold  ${E}_{u}=\mathrm{span}(\mathbf{u}^u)$, where $\mathbf{u}^u$ is the   eigenvector associated with the negative eigenvalue $\lambda_m=-n\delta$. In the right panel,  we  plot   $\|\bar{\mathbf{x}}^\nu - \boldsymbol{\theta}^{\ast}\|$ versus the number of iterations. All the curves are  averaged over 50 independent initializations. Figure in the left panel shows that, as predicted by our theory,  both algorithms almost surely escapes from the  unstable subspace $E_u$, at an indistinguishable practical   rate.  
The right panel shows that   DOGT gets closer to the strict saddle; this can be justified by the fact that, differently from  DGD, DOGT  exhibits \emph{exact} convergence  to   critical points.\vspace{-0.1cm}
{ 
	\subsection{Bilinear logistic regression \cite{Dyrholm_JMLR_2007}}\label{BLR_simulations}
	Consider a classification problem with distributed training data set $\{\mathbf{s}_i,\xi_i\}_{i=1}^n$, where $\mathbf{s}_i\in\mathbb{R}^d$ is the feature vector associated with the binary class label $\xi_i\in\{0,1\}$. The bilinear logistic regression problem aims at finding the bilinear classifier  $\zeta_i(\mathbf{Q},\mathbf{w};\mathbf{s}_i)=  \mathbf{s}_i^\top\mathbf{Q}\mathbf{w}$, with $\mathbf{Q}\in\mathbb{R}^{d\times p}$ and $\mathbf{w}\in\mathbb{R}^p$ that best separates   data with distinct labels. Let $(\mathbf{s}_i,\xi_i)$ be private information for agent $i$.   Using the sigmoid activation function $\sigma(x)\triangleq 1/(1+e^{-x})$ together with the \emph{cross-entropy risk} function, the optimization problem reads
	\begin{equation}
	\min_{\mathbf{Q},\mathbf{w}}\quad  -\frac{1}{n} \sum_{i=1}^n \Big[\xi_i  \ln\big(\sigma(\mathbf{s}_i^\top\mathbf{Q}\mathbf{w})\big)+(1-\xi_i) \ln\big(1-\sigma(\mathbf{s}_i^\top\mathbf{Q}\mathbf{w})\big)\Big]+\frac{\tau}{2}\left(\norm{\mathbf{Q}}^2_F+\norm{\mathbf{w}}^2\right).
	\label{eq:P_BilinearLogisticReg}
	\end{equation}
	It is not difficult to show that $\eqref{eq:P_BilinearLogisticReg}$ is equivalent to the following instance of \eqref{eq:P}: \vspace{-0.2cm}
	\begin{equation}
	\min_{\mathbf{Q},\mathbf{w}}\quad F(\mathbf{Q},\mathbf{w})=  \sum_{i=1}^n  ~\underbrace{\frac{1}{n}\left[-\ln\big(\sigma(\tilde{\xi}_i\mathbf{s}_i^\top\mathbf{Q}\mathbf{w})\big)+\frac{\tau}{2}\left(\norm{\mathbf{Q}}^2_F+\norm{\mathbf{w}}^2\right)\right]}_{= f_i(\mathbf{Q},\mathbf{w})},
	\label{eq:P_BilinearLogisticReg_EquivalentForm}
	\vspace{-0.3cm}
	\end{equation}
	with\vspace{-0.4cm}
	\begin{equation*}
	\tilde{\xi_i}\triangleq \left\{
	\begin{array}{l l}
	-1, & \mathrm{if}\quad \xi_i=0;
	\\
	1, & \mathrm{if}\quad \xi_i=1.
	\end{array}
	\right.
	\end{equation*}
	To visualize the landscape of $F(\mathbf{Q},\mathbf{w})$ (2D plot), we consider the following setting for the free parameters.  We set  $d=p=1$; $\tau=0.2$; $n=5$; and  we generate uniformly random $\tilde{\xi}_i\in\{0,1\}$, and we draw $s_i$  from a  normal distribution with mean $\xi_i$ and variance $1$. 
	The gradient of the local loss $f_i$ reads
	\begin{equation*}
	\begin{bmatrix}
	\nabla_{Q} f_i(Q,w)
	\\
	\nabla_{w} f_i(Q,w)
	\end{bmatrix}=\frac{1}{n}
	\begin{bmatrix}
	\tau Q-\tilde{\xi}_is_iw\sigma(-\tilde{\xi}_is_iQw)
	\\
	\tau w-\tilde{\xi}_is_iQ\sigma(-\tilde{\xi}_is_iQw)
	\end{bmatrix}.
	\end{equation*}
	A surface plot of  $F({Q},{w})$ in the above setting  is plotted  in the right panel of Fig. \ref{Contour_plots}. 
	Note that such $F$ has three critical points, two of which are local minima (see the location of minima in the left or middle panel of Fig. \ref{Contour_plots}  marked by $\times$) and one strict saddle point at $(0,0)$--the Hessian at $(0,0)$,\vspace{-0.2cm}
	\begin{equation*}
	\nabla^2 F(0,0)=
	\begin{bmatrix}
	\tau & -\frac{1}{2n}\sum_i\tilde{\xi}_is_i
	\\
	-\frac{1}{2n}\sum_i\tilde{\xi}_is_i & \tau
	\end{bmatrix},\vspace{-0.15cm}
	\end{equation*}
	has an eigenvalue at $\tau-\frac{1}{2n}\sum_i\tilde{\xi}_is_i=-0.26$.

	%
	\begin{figure}[t!]
	\hspace{-.3cm}\includegraphics[scale=0.32]{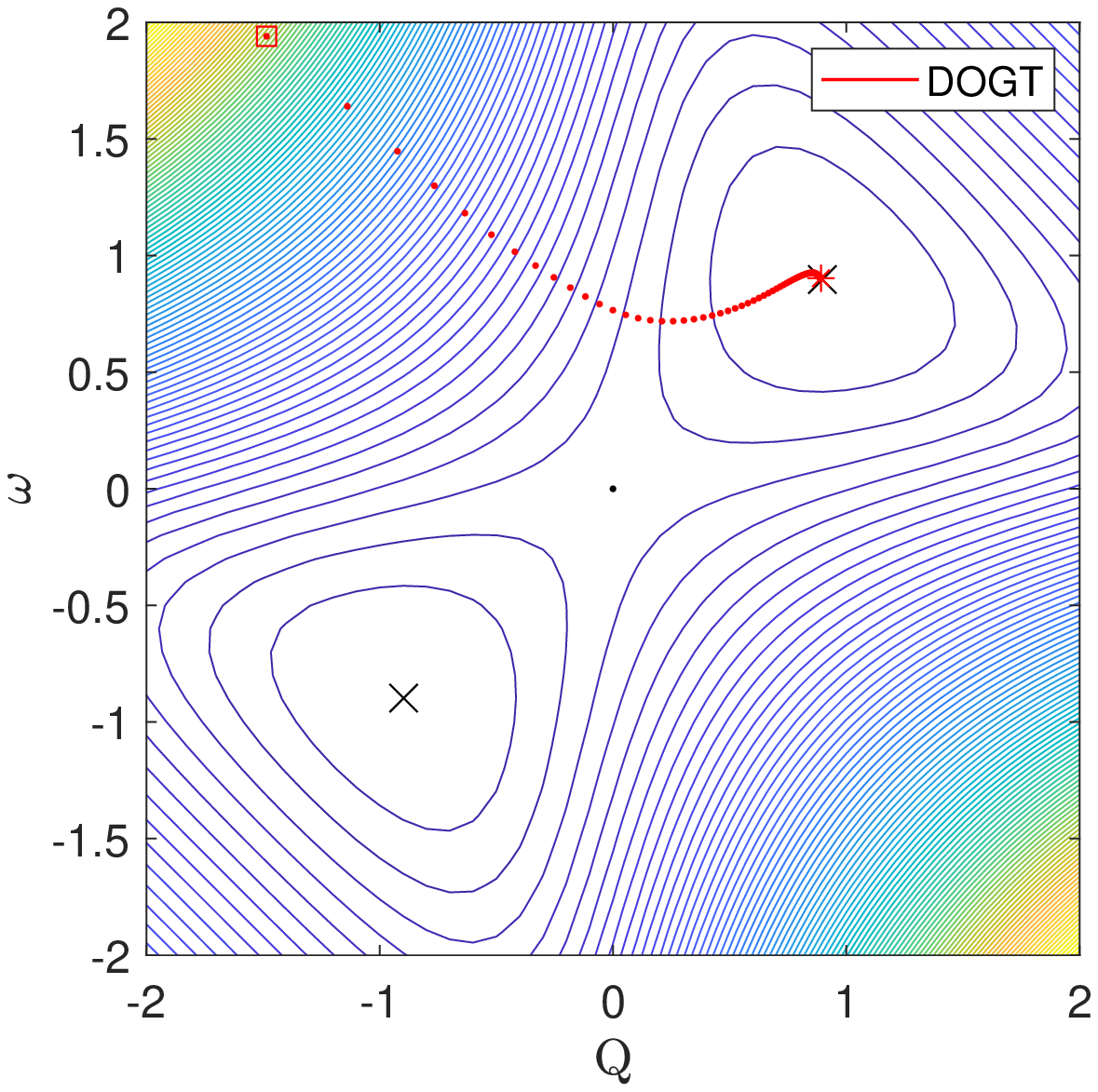}\hspace{-.2cm}\includegraphics[scale=0.32]{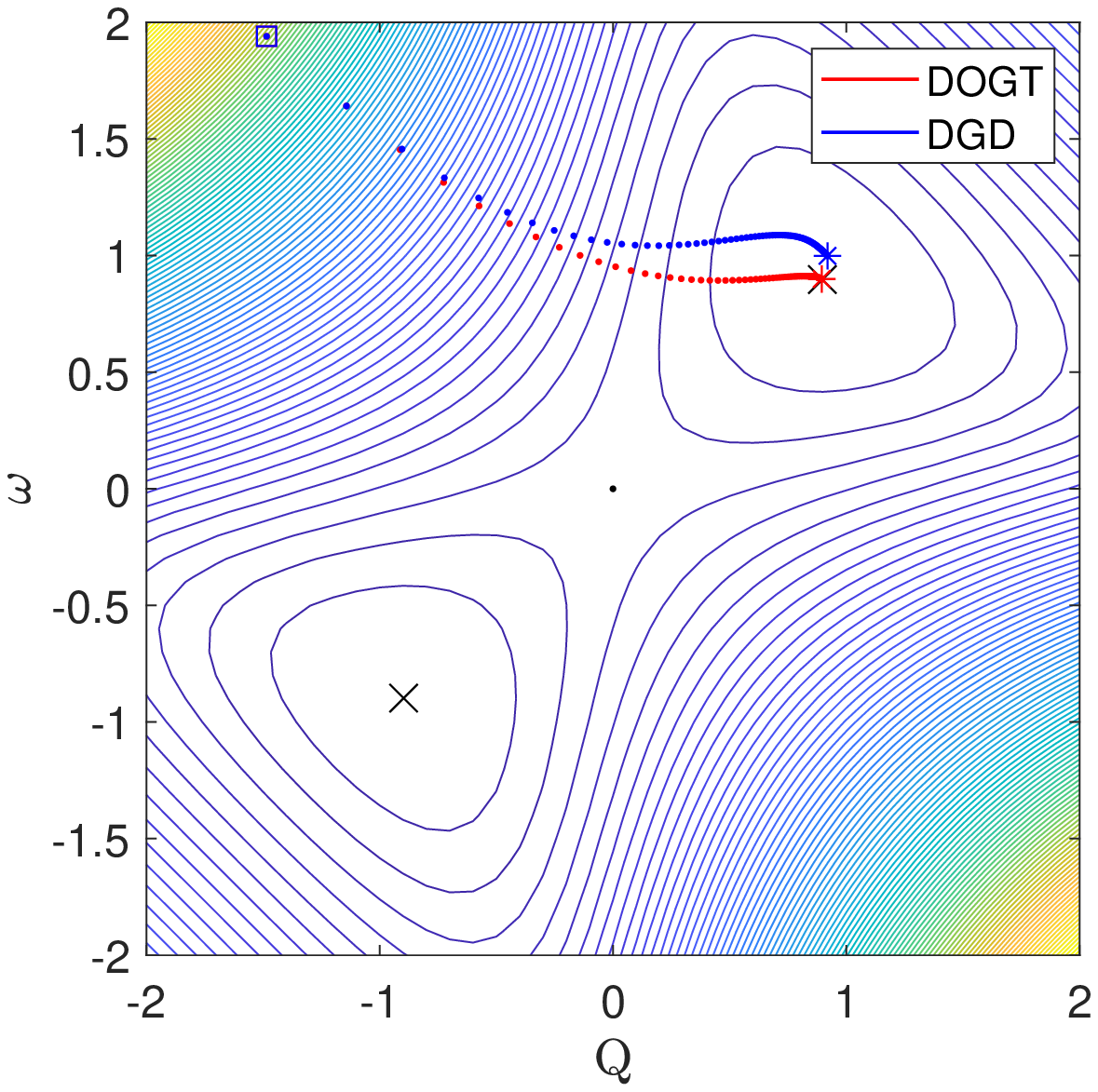}\hspace{-.15cm}\includegraphics[scale=0.43]{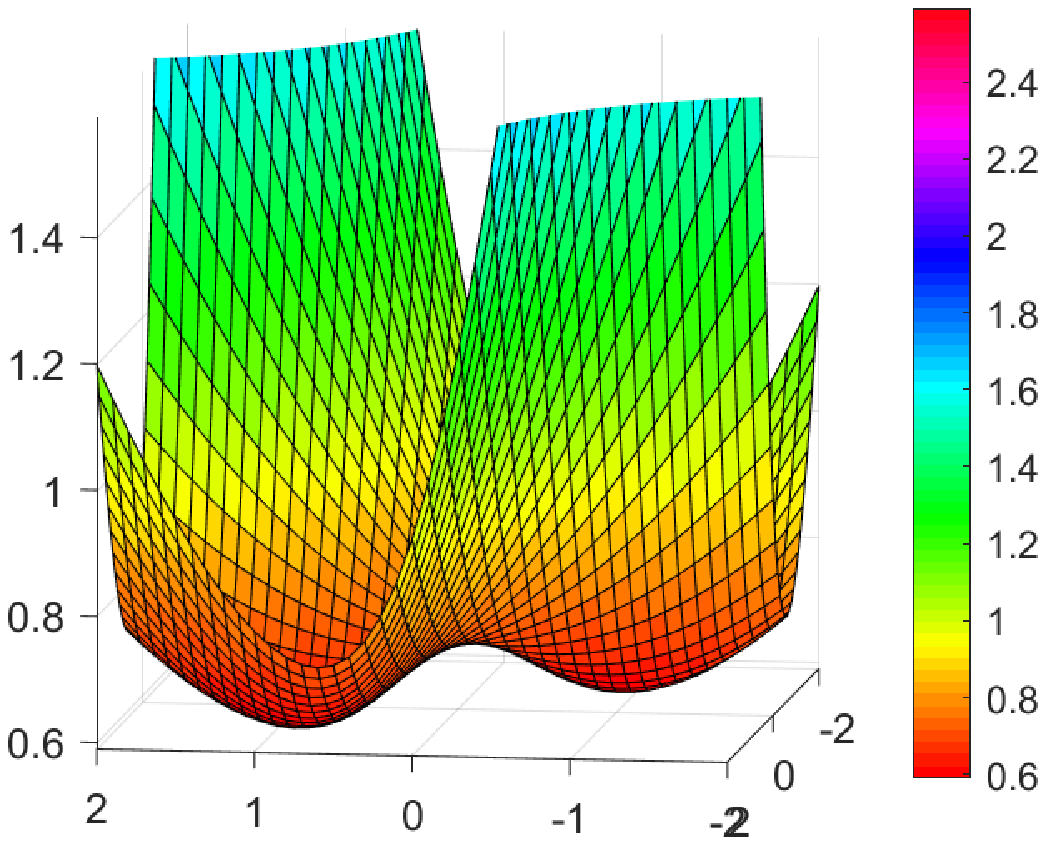}\vspace{-0.2cm}
		\caption{Escaping properties of the DGD and DOGT, applied to the bilinear logistic regression problem \eqref{eq:P_BilinearLogisticReg}. Left (resp. middle) plot: directed (resp. undirected) network;   trajectory of the average iterates on the contour  of $F$ ($(0,0)$ is the strict saddle point and $\times$ are the local minima); DGD and DOGT are  initialized at $\Box$ and terminated after $100$ iterations at $\ast$. Right plot:  plot of $F$.\vspace{-0.5cm}}
		\label{Contour_plots}
	\end{figure}

	We test DGD and DOGT   over a network of  $n=5$ agents; for DGD we considered undirected graphs whereas we run DOGT on both undirected and directed graphs. Both algorithms are initialized at the same random point and terminated after $100$ iterations;  the step-size is set to $\alpha=0.9$. We denote by $Q^\nu_i$ and $w^\nu_i$  the agent $i$'s $\nu$-the iterate of the local copies of  $Q$ and $w$, respectively. The trajectories of the average iterates $(\bar{Q}^\nu,\bar{w}^\nu)\triangleq \frac{1}{n}(\sum_iQ^\nu_i,\sum_iw^\nu_i)$ are plotted in Fig. \ref{Contour_plots}; the left panel refers to the directed graph while the  middle panel  reports the same results for the   undirected network. As expected, the DOGT algorithm converges to an exact critical point (local minimum) avoiding the strict saddle $(0,0)$ while DGD converges to a neighborhood  of the local minimum. The consensus error is  $1/n\sqrt{\sum_{i=1}^n||(Q_i^\nu,w_i^\nu)-(\bar{Q}^\nu,\bar{w}^\nu)||^2}$;  at the termination, it reads  $2.33\times 10^{-4}$ for DOGT over  the directed  network, and   $2.18\times 10^{-4}$ and $9.74\times 10^{-2}$ for DOGT and DGD, respectively, over undirected networks.

	\subsection{Gaussian mixture model}\label{GMM_simulations}
	Consider the Gaussian mixture model defined in Sec. \ref{ProblemSetting_SubSe}. The data  $\{\mathbf{z}_i\}_{i=1}^n$ where $\mathbf{z}_i\in\mathbb{R}^m$  are  realizations of the mixture model $\mathbf{z}_i\sim \frac{1}{2}\mathcal{N}(\boldsymbol{\mu}_1,\boldsymbol{\Sigma}_1)+\frac{1}{2}\mathcal{N}(\boldsymbol{\mu}_2,\boldsymbol{\Sigma}_2)$. Let each agent $i$ own $\mathbf{z}_i$. Both parameters $(\boldsymbol{\mu}_1,\boldsymbol{\mu}_2)$ and $(\boldsymbol{\Sigma}_1,\boldsymbol{\Sigma}_2)$ are unknown. The goal is  to approximate $(\boldsymbol{\mu}_1,\boldsymbol{\mu}_2)$ while $(\boldsymbol{\Sigma}_1,\boldsymbol{\Sigma}_2)$ is set to an estimate $(\tilde{\boldsymbol{\Sigma}},\tilde{\boldsymbol{\Sigma}})$. The problem reads\vspace{-0.2cm}
	\begin{equation}
	\min_{\boldsymbol{\theta}_1,\boldsymbol{\theta}_2\in \mathbb{R}^{m}}\, -\sum_{i=1}^n\log{\left(\phi_m(\mathbf{z}_i-\boldsymbol{\theta}_1)+\phi_m(\mathbf{z}_i-\boldsymbol{\theta}_2)\right)},\vspace{-0.2cm}
	\label{eq:GaussianMixtureProblem}
	\end{equation}
	where $\phi_m(\boldsymbol{\theta})$ is the $m$-dimensional normal distribution with mean $\mathbf{0}$ and covariance $\tilde{\boldsymbol{\Sigma}}$.
Consider the case of mixture of two scalar Gaussians, i.e., $m=1$. We draw $\{\mathbf{z}_i\}_{i=1}^5$ from the  the this model, with means $\mu_1=0$, $\mu_2=-5$ and variance $\sigma_1=\sigma_2=25$. The estimate of variance in problem \eqref{eq:GaussianMixtureProblem} is pessimistically set to $\tilde{\sigma}=1$. A surface plot of a random instance of above problem is depicted in right panel of Fig. \ref{GausMix_Contour_plots}. Note that this instance of problem has 2 global minima (marked by $\times$) and multiple local minima.  
We test  DGD and DOGT on the above problem over the same networks as described in Sec. \ref{BLR_simulations}.   Both algorithms are initialized at the same random point and terminated after $250$ iterations; the step-size is set to   $\alpha=0.1$.  In Fig. \ref{GausMix_Contour_plots}, we plot the trajectories of the average iterates $(\bar{\theta}_1^\nu,\bar{\theta}_2^\nu)\triangleq \frac{1}{n}(\sum_i\theta^\nu_{1,i},\sum_i\theta^\nu_{2,i})$, where   $\theta^\nu_{1,i}$ and $\theta^\nu_{2,i}$ are the agent $i$'s $\nu$-the iterate of the local copies of  $\theta_1$ and $\theta_2$, respective; the  left (resp. middle) panel refers to the undirected (resp. directed) network.  DOGT  converges to  the global minimum  while  DGD happens to converge to neighborhood of a local minima. The consensus error is measured by $(1/n)\sqrt{\sum_{i=1}^n||(\theta_{1,i}^\nu,\theta_{2,i}^\nu)-(\bar{\theta}_1^\nu,\bar{\theta}_2^\nu)||^2}$ and at the termination it is equal to $1.9\times 10^{-3}$ for DOGT on the directed graph; and   $2.8\times 10^{-3}$ and $1.135$ for DOGT and DGD, respectively over the undirected graph.
	
	
	\begin{figure}[t!]
			\hspace{-.15cm}\includegraphics[scale=0.32]{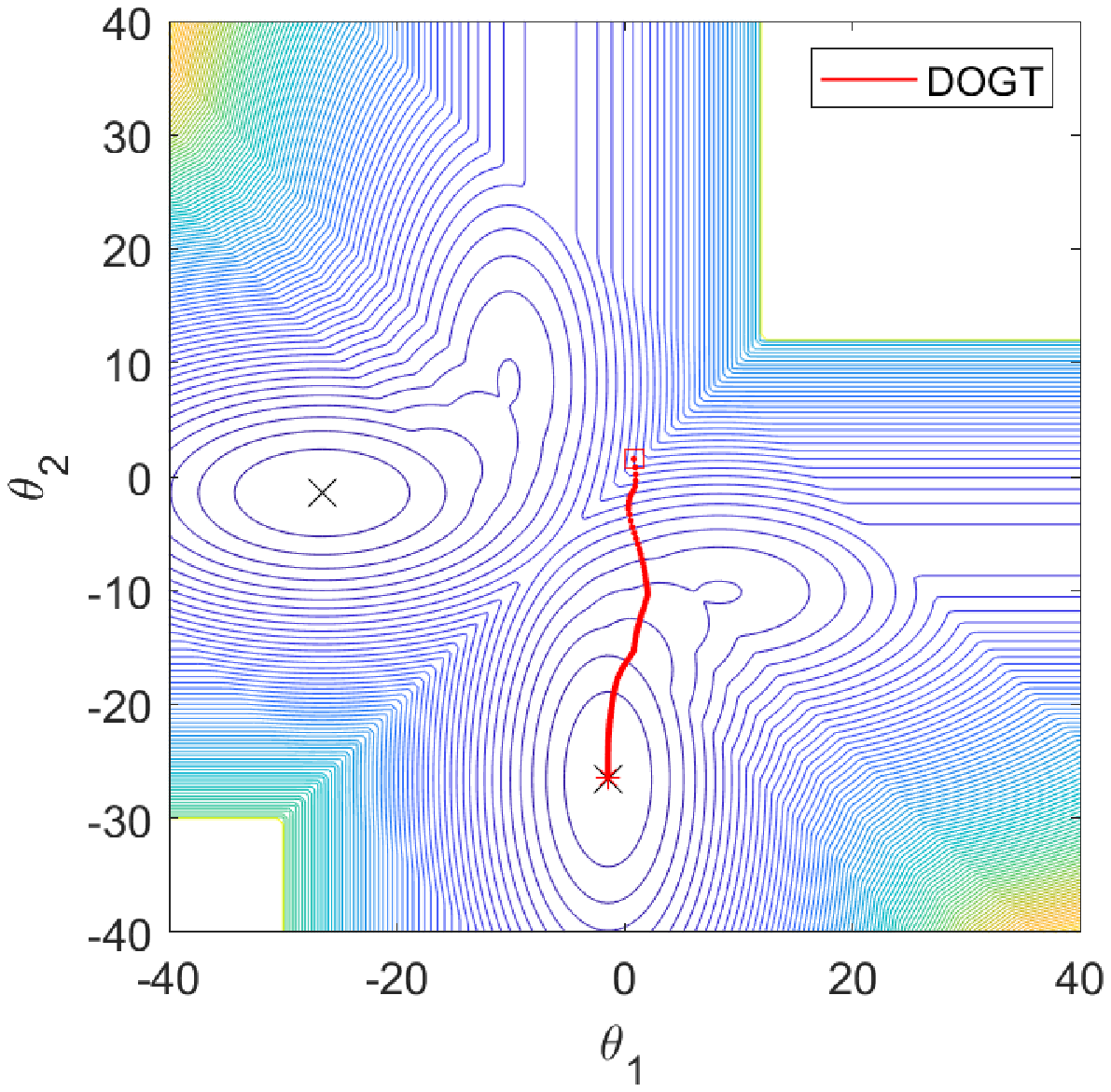}	\hspace{-.45cm}\includegraphics[scale=0.32] {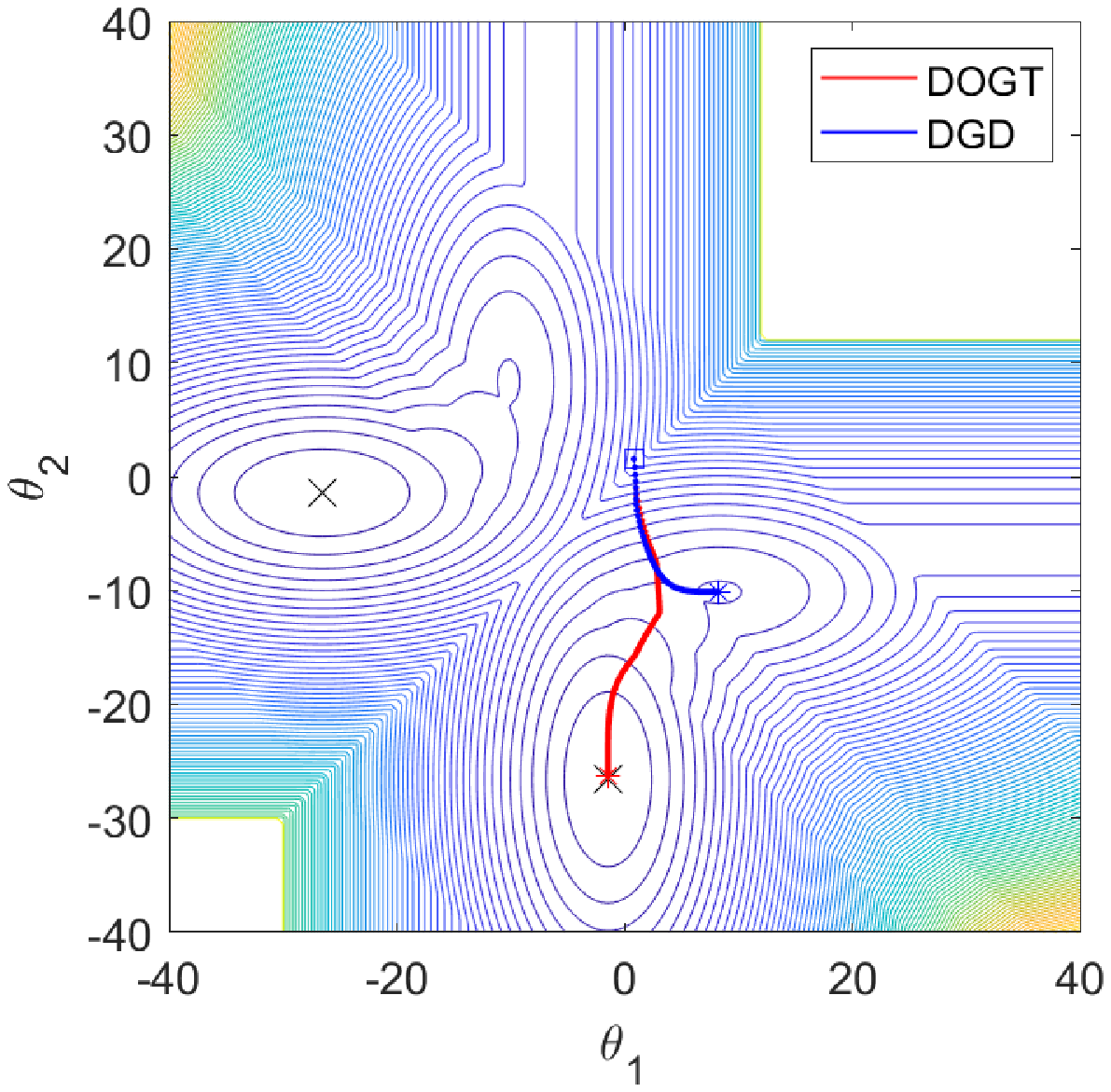}\hspace{-.32cm}\includegraphics[scale=0.38]{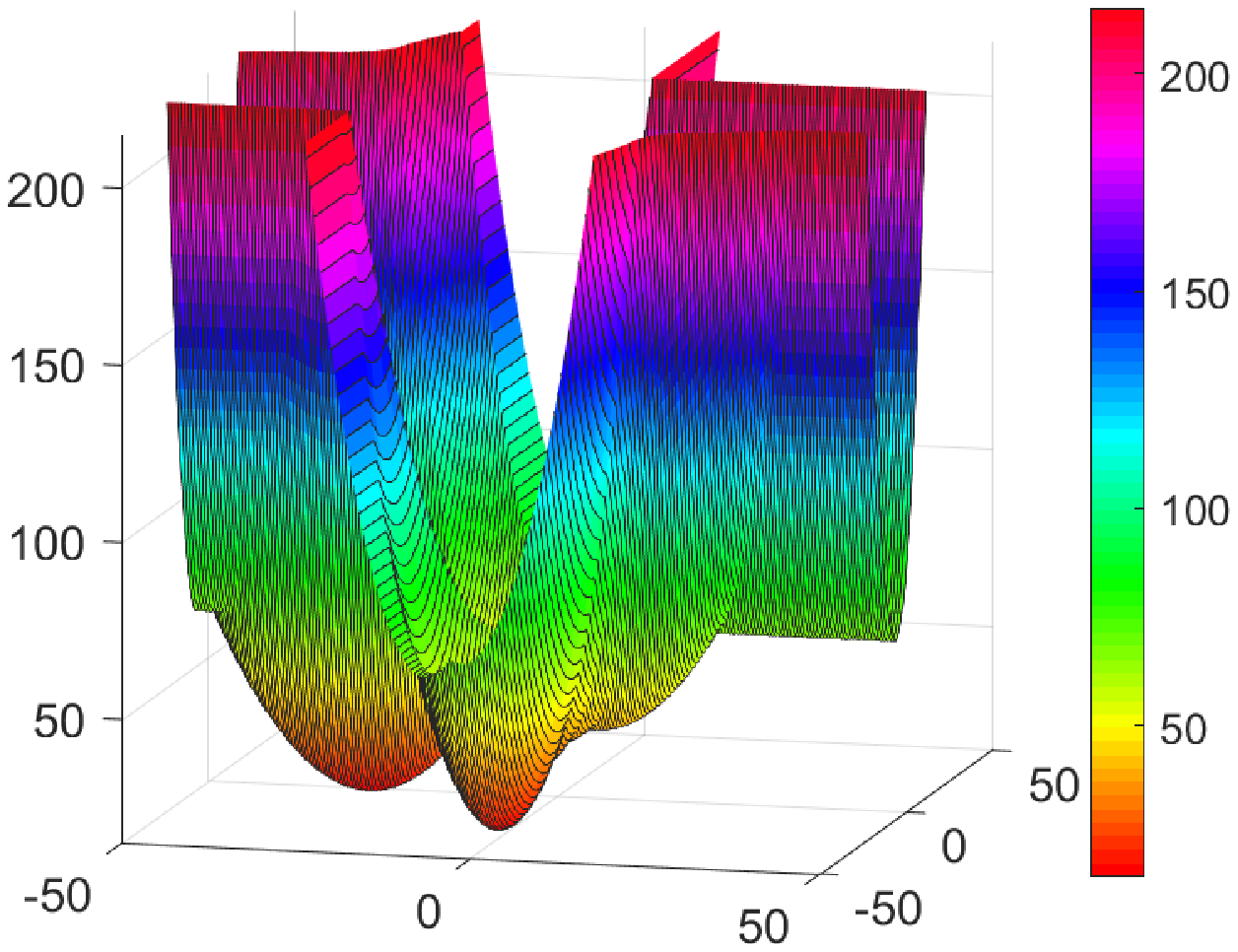}
		\vspace{-0.2cm}
		\caption{Escaping properties of the DGD and DOGT applied to the Gaussian mixture problem (\ref{eq:GaussianMixtureProblem}). Left (resp. middle) plot: directed (resp. undirected) network;   trajectory of the average iterates on the contour  of $F$ (the global minima are marked by $\times$); DGD and DOGT are  initialized at $\Box$ and terminated after $250$ iterations at $\ast$. Right plot:  plot of $F$.\vspace{-0.5cm}}
		\label{GausMix_Contour_plots}
	\end{figure}
	
	%
	%
	%
	%
	
}	

\vspace{-0.1cm}

\appendix   
\section{Appendix}\vspace{-0.1cm}

\subsection{On the problems  satisfying Assumption \ref{Strip_assumption}}\label{Strip_assumption_appendix}
{ We prove that all the functions arising from the examples in Sec.\ref{ProblemSetting_SubSe} satisfy Assumption \ref{Strip_assumption}, for sufficiently large $R$ and $R-\epsilon$. To do so, for each function, we  establish lowerbounds  implying $\langle\nabla f_i(\bt),\bt/\norm{\bt}\rangle\rightarrow\infty$ as $||\bt||\rightarrow\infty$. 
	
	\textbf{a) Distributed PCA: }
	Let us expand the objective function in \eqref{Dist_PCA} as
	\begin{equation*}
	\begin{aligned}
	F(\bt)=&\frac{1}{4}\trace{\bt\bt^\top\bt\bt^\top}-\frac{1}{2}\trace{\bt^\top \sum_{i=1}^n\mathbf{M}_i\bt}+\frac{1}{4}\trace{\sum_{i=1}^n\mathbf{M}^\top_i\sum_{i=1}^n\mathbf{M}_i}
	\\
	=&\sum_{i=1}^n\underbrace{\frac{1}{4}\left\{\frac{1}{n}\big\|\bt\bt^\top\big\|_F^2-2\trace{\bt^\top \mathbf{M}_i\bt}\right\}}_{\triangleq f_i(\bt)}+\frac{1}{4}\trace{\sum_{i=1}^n\mathbf{M}^\top_i\sum_{i=1}^n\mathbf{M}_i}.
	\end{aligned}
	\end{equation*}
	We have\vspace{-0.3cm}
	\begin{equation*}
	\begin{aligned}
	\left\langle\nabla f_i(\bt),\bt/\norm{\bt}\right\rangle=&\left\langle\frac{1}{n}\bt\bt^\top\bt-\mathbf{M}_i\bt,\bt\right\rangle/\norm{\bt}
	\\
	=&\frac{1}{n}\norm{\bt\bt^\top}_F^2/\norm{\bt}-\bt^\top\mathbf{M}_i\bt/\norm{\bt}
	\\
	\geq & \frac{1}{nK_{2,4}^4}\norm{\bt}^3-\sigma_{\max}(\mathbf{M}_i)\norm{\bt},
	\end{aligned}
	\end{equation*}
	for some $K_{2,4}>0$, where in the last inequality we used the equivalence of  $\ell_4$ and $\ell_2$ norms, i.e. $||\bt||_2\leq K_{2,4}||\bt||_4,\forall \bt\in\mathbb{R}^m$.

	\textbf{b) Phase retrieval:} 
	It is not difficult to show that for the objective function in \eqref{PR_problem}, it holds\vspace{-0.1cm}
	\begin{equation*}
	\begin{aligned}
	\left\langle \nabla f_i(\bt),\bt/\norm{\bt}\right\rangle =&\left(||\mathbf{a}_i^\top\bt||^2-y_i\right)||\mathbf{a}_i^\top\bt||^2/\norm{\bt}+\lambda\norm{\bt}
	\\
	&=\left(||\mathbf{a}_i^\top\bt||^2-y_i/2\right)^2/\norm{\bt}-\frac{y_i^2}{4\norm{\bt}}+\lambda\norm{\bt}.
	\end{aligned}
	\end{equation*}
	
	\textbf{c) Matrix sensing:} Consider the objective function in \eqref{PR_problem}; it holds
	%
	It is not difficult to show that
	\begin{equation*}
	\begin{aligned}
	\left\langle\nabla f_i(\bT),\bT/\norm{\bT}_F\right\rangle = & \left(\trace{\bT^\top\mathbf{A}_i\bT} -y_i\right)\trace{\bT^\top\mathbf{A}_i\bT}/\norm{\bT}_F+\lambda\norm{\bT}_F
	\\
	= & \trace{\bT^\top\mathbf{A}_i\bT}^2/\norm{\bT}_F-y_i\trace{\bT^\top\mathbf{A}_i\bT}/\norm{\bT}_F+\lambda\norm{\bT}_F
	\\
	= & \left(\trace{\bT^\top\mathbf{A}_i\bT} -y_i/2\right)^2/\norm{\bT}_F-\frac{y_i^2}{4\norm{\bT}_F}+\lambda\norm{\bT}_F.
	\end{aligned}
	\end{equation*}
	
	\textbf{d-f)}
	We prove the property only for the   Gaussian mixture model; similar proof applies also to the other classes of problems.  Denote $\bt=(\bt_d)_{d=1}^q$. Since $\phi_m$ is a bounded function, we have \vspace{-.4cm}
	\begin{equation*}
	\left\langle\nabla_{\theta_d}f_i(\bt_d),\bt_d\right\rangle\geq -C_d+\lambda \norm{\bt_d}^2,
	\end{equation*} 
	for some $C_d>0$. Hence, 
	$\left\langle\nabla_{\theta}f_i(\bt),\bt/\norm{\bt}\right\rangle \geq -C/\norm{\bt}+\lambda \norm{\bt},$
	with $C=\sum_dC_d$. 
	
}


\subsection{Convergence of DGD without $L-$smoothness of $f_i$'s}
\label{DGD_conv_cons_thm_extended_proof_appendix}
We sketch here how to extend the convergence results of DGD stated in Sec. \ref{sec:DGD} to the case when the gradient of the agents' loss functions is not globally Lipschitz continuous (i.e. removing Assumption \ref{P_assumption} (i)). Due to the space limitation, we prove only the counterpart of Theorem \ref{DGD_conv_cons_thm}; the other results in Sec. \ref{sec:DGD} can be extended following similar arguments.  

We begin  introducing some definitions. Under Assumptions \ref{Strip_assumption} and  \ref{DoublyStochastic_D}, define the  set $\tilde{\mathcal{Y}}\triangleq \mathcal{Y}+\mathcal{B}^{mn}_b$, with 
$\mathcal{Y}=\bar{\mathcal{L}}\cup\prod_{i=1}^n\mathcal{B}^m_R$ and 
\begin{equation}
\bar{\mathcal{L}}= \mathcal{L}_{F_c}\left(\max_{\mathbf{x}_i^0\in\mathcal{B}^m_R,   i\in [n]} \left\{\sum_{i=1}^n f_i(\mathbf{x}_i^0)\right\}+\frac{R^2}{\alpha_b}\right),
\end{equation}
where  
\begin{equation}
\begin{aligned}
& \qquad  \alpha_b=  \min_{i\in[n]}\min\{\epsilon D_{ii}/h,2D_{ii}\delta(R-\epsilon)/h^2\}>0,
\\
 h= & \max_{i\in[n],\mathbf{z}\in\mathcal{B}^m_R}||\nabla f_i(\mathbf{z})||,\quad \text{and}\quad b=\max_{\alpha\in[\alpha_b,1],\boldsymbol{\theta}\in\mathcal{Y}}||\nabla L_\alpha(\boldsymbol{\theta})||.
\end{aligned}
\end{equation}
Note that, under Assumption  \hyperlink{P_assumption_prime}{2.1'(ii)}, $\mathcal{Y}$ and $\tilde{\mathcal{Y}}$ are compact. Hence, $\nabla F_c$ is globally Lipschitz on $\tilde{\mathcal{Y}}$, and so is $\nabla L_\alpha$; we denote such  Lipschitz constants  as $\tilde{L}_{\nabla F_c}$ and $\tilde{L}_{\nabla L_\alpha}$, respectively; it is not difficult to check that  
\begin{equation}
\tilde{L}_{\nabla L_\alpha}=\tilde{L}_{\nabla F_c}+\frac{1-\sigma_{\min}(\mathbf{D})}{\alpha_b}.
\end{equation}   The following result replaces Theorem \ref{DGD_conv_cons_thm} in the above setting. 
\begin{theorem}
	\label{DGD_conv_cons_thm_extendedver2}  
	 Consider Problem \eqref{eq:P}, under  Assumptions  \hyperlink{P_assumption_prime}{2.1'(ii)}, \ref{Strip_assumption} and \ref{Net-Assump}. Let  $\{\mathbf{x}^\nu\}$ be the sequence generated by  DGD in \eqref{DGD_update} under Assumption \ref{DoublyStochastic_D}, with  $\mathbf{x}_i^0\in\mathcal{B}^m_R,  i\in[n]$ and $0<\alpha< \bar{\alpha}_{\max}\triangleq \sigma_\mathrm{min}(\mathbf{I}+\mathbf{D})/\tilde{L}_{\nabla F_c}$.   Then, same conclusions of Theorem \ref{DGD_conv_cons_thm} hold. 
\end{theorem} 

\begin{proof}
It is sufficient to show that $\{\mathbf{x}^\nu\}\subseteq \mathcal{Y}$; the rest of the proof follows similar steps as those in  \cite[lemma 2]{WYin_DGD_SIAM2016}     replacing $L_c$ with $\tilde{L}_{\nabla F_c}$. 

When $\alpha<\alpha_b$, $\{\mathbf{x}^\nu\}\subseteq \mathcal{Y}$  can be proved leveraging the same arguments used in the proof of  Lemma \ref{Y_set_found}. Therefore, in the following, we consider only the case $\alpha_b\leq \alpha <\sigma_\mathrm{min}(\mathbf{I}+\mathbf{D})/\tilde{L}_{\nabla F_c}$, with $\alpha_b< \sigma_\mathrm{min}(\mathbf{I}+\mathbf{D})/\tilde{L}_{\nabla F_c}$. 
We prove by induction. Clearly $\mathbf{x}^0\in\mathcal{L}_{L_\alpha}(L_\alpha(\mathbf{x}^0))$ and, by \eqref{L_0_Y_inclusion} (cf. Lemma \ref{Y_set_found}),
$$\mathcal{L}_{L_\alpha}(L_\alpha(\mathbf{x}^0))\subseteq \bar{\mathcal{L}}\subseteq \mathcal{Y},\quad \forall \alpha\in[\alpha_b,1].$$
 Assume $\mathcal{L}_{L_\alpha}(L_\alpha(\mathbf{x}^\nu))\subseteq \mathcal{Y}$. Since  $\mathbf{x}^\nu\in {\mathcal{Y}}$, there hold  $\bx^{\nu+1}=\bx^\nu-\alpha\nabla L_\alpha(\bx^\nu)\in \tilde{\mathcal{Y}}$ and $\theta \mathbf{x}^\nu + (1-\theta) \mathbf{x}^{\nu+1}\in \tilde{\mathcal{Y}}$, for all $\theta\in [0,1]$. Invoking the descent lemma on   $L_\alpha$ at $\mathbf{x}^{\nu+1}$  [recall that $L_\alpha$ is $\tilde{L}_{\nabla L_\alpha}$-smooth on $\tilde{\mathcal{Y}}$], we have:
\begin{equation}
\label{Descent_noLipschitz_1}
\begin{aligned}
L_{\alpha}(\bx^{\nu+1})\leq L_{\alpha}(\bx^\nu)-\alpha\left(\frac{\sigma_{\min}(\mathbf{I}+\mathbf{D})-\alpha \tilde{L}_{\nabla F_c}}{2}\right)\norm{\nabla L_\alpha(\bx^\nu)}^2\leq L_{\alpha}(\bx^\nu).
\end{aligned}
\end{equation}
Therefore,   $\mathcal{L}_{L_\alpha}(L_\alpha(\mathbf{x}^{\nu+1}))\subseteq\mathcal{L}_{L_\alpha}(L_\alpha(\mathbf{x}^\nu))\subseteq \mathcal{Y}$, which completes the induction. 
\end{proof}
 
\subsection{Proof of Theorem \ref{Rate_theorem}: Supplement}
\label{Rate_theorem_proof_appendix}
We first show that, if there exists some $\nu_0$ such that  $d^{\nu_0}=0$,    $\mathbf{z}^{\nu}=\mathbf{z}^{\nu_0}$, for all $\nu\geq\nu_0$ [see updates in \eqref{NEXT_generic}]; this means that   $\{\mathbf{z}^\nu\}$   converges in finitely many iterations. Define $\mathcal{D}\triangleq \{\nu:d^\nu\neq 0\}$ and take $\nu$ in $\mathcal{D}$. Let $\theta=0$, then the {K\L} inequality yields $||\nabla L(\mathbf{x}^\nu,\mathbf{y}^\nu)||\geq 1/c$, for all $\nu\in\mathcal{D}$. This   together with \eqref{L_descent} and Lemma \ref{nablaL_bound_lemma}, lead to
$l^{\nu+1}\leq l^{\nu}- {1}/{(Mc)^2}$,
which by Assumption \ref{P_assumption}-(ii), implies that $\mathcal{D}$ must be finite and  $\{\mathbf{z}^\nu\}$   converges in a finite number of  iterations. 

Consider \eqref{rate_bound3}. Let $\theta\in(0,1/2]$, then $(1-\theta)/\theta\geq 1$. Since $D^\nu\rightarrow 0$ as $\nu\rightarrow\infty$ [by Lemma \ref{sequence_props_a}-(ii)],   there exists a sufficiently large $\nu_0$ such that $(D^{\nu}-D^{\nu+1})^{(1-\theta)/\theta}\leq D^{\nu}-D^{\nu+1}$. By  \eqref{rate_bound3}, we have\vspace{-0.3cm}
\begin{equation*}
D^{\nu+1}\leq \frac{\tilde{M}Mc-1}{\tilde{M}Mc}D^\nu,\vspace{-0.2cm}
\end{equation*}
which proves case (ii).

Finally, let us assume  $\theta\in(1/2,1)$, then $\theta/(1-\theta)> 1$. Eq.   \eqref{rate_bound3} implies\vspace{-0.1cm}
\begin{equation*}
1\leq\frac{\bar{M}(D^{\nu}-D^{\nu+1})}{\left(D^{\nu}\right)^{\theta/(1-\theta)}}
\label{rate_bound5}\vspace{-0.1cm}
\end{equation*}
where $\bar{M}=(M\tilde{M}c)^{\theta/(1-\theta)}$. Define $h:(0,+\infty)\rightarrow\mathbb{R}$ by $h(s)\triangleq s^{-\frac{\theta}{1-\theta}}$.
Since $h$ is monotonically decreasing over $[D^{\nu+1},D^{\nu}]$, we get
\begin{equation}
\begin{aligned}
1\leq    \bar{M}(D^{\nu}-D^{\nu+1})h(D^\nu)
\leq  \bar{M} \int_{D^{\nu+1}}^{D^\nu}h(s)ds
=   \bar{M}\frac{1-\theta}{1-2\theta}\left((D^{\nu})^p-(D^{\nu+1})^p\right),
\end{aligned}
\label{rate_bound6}
\end{equation}
with $p=\frac{1-2\theta}{1-\theta}<0$. By \eqref{rate_bound6} one infers   that there exists a constant $\mu>0$ such that 
$(D^{\nu+1})^p -(D^{\nu})^p\geq \mu$.
The following chain of implications then holds: 
$(D^{\nu+1})^p \geq \mu \nu+(D^{1})^p \implies D^{\nu+1} \leq  \left(\mu \nu+(D^{1})^p\right)^{1/p}\implies D^{\nu+1} \leq   C_0\nu^{1/p},$
for some constant $C_0>0$. This proves case (iii).

\subsection{Extension of Proposition \ref{CSS_is_in_unstable_set}}
\label{CSS_is_in_unstable_set_extension_Appendix}
We relax conditions (i)-(ii) of Proposition \ref{CSS_is_in_unstable_set} under the following additional mild assumptions on the set of strict saddle points and the weight matrices $\mathbf{R}$ and $\mathbf{C}$. 

\begin{assumption}
	\label{Strict_Saddles_delta_assumption}
	There exists $\delta>0 $  such that  $\lambda_{\min} (\nabla^2 F(\boldsymbol{\theta}^\ast))\leq -\delta$, for all $\boldsymbol{\theta}^\ast\in \Theta^\ast_{ss}$ ($\Theta^\ast_{ss}$ is the set of strict saddle of $F$, cf. Definition \ref{SSP_def}).
\end{assumption}
 
\begin{assumption}\label{diag_dominant_RC}
The matrices $\mathbf{R}$ and $\mathbf{C}$ are chosen   according to
\begin{equation*}
\mathbf{R}=\frac{\widetilde{\mathbf{R}}+(t-1)\mathbf{I}}{t},\quad \mathbf{C}=\frac{\widetilde{\mathbf{C}}+(t-1)\mathbf{I}}{t},
\end{equation*}
for some $t\geq 1$, and some matrices $\widetilde{\mathbf{R}}$ and $\widetilde{\mathbf{C}}$ satisfying Assumption \ref{matrix_C_R}.
\end{assumption}

Note that   $\mathbf{R}$ and $\mathbf{C}$   satisfy  Assumption \ref{matrix_C_R} as well. The main result is given in Proposition \ref{CSS_is_in_unstable_set_extended}. Before proceeding, we recall the following result on spectral variation of non-normal matrices.
\begin{theorem}{\cite[Theorem VIII.1.1]{bhatia2013matrix}}
	\label{th:matpert}
	For arbitrary $d\times d$ matrices $\mathbf{A}$ and $\mathbf{B}$, it holds that \vspace{-0.3cm}
	\begin{equation*}
	s\left(\sigma(\mathbf{A}),\sigma(\mathbf{B})\right)\le \left(\|\mathbf{A}\|+\|\mathbf{B}\|\right)^{1-1/d} \|\mathbf{A}-\mathbf{B}\|^{1/d}\vspace{-0.2cm}
	\end{equation*}
	with\vspace{-0.3cm}
	\begin{equation*}
	s\left(\sigma(\mathbf{A}),\sigma(\mathbf{B})\right)\triangleq \max_j\min_i |\alpha_i-\beta_j|,
	\end{equation*}
	where $\alpha_1,\ldots,\alpha_d$ and  $\beta_1,\ldots,\beta_d$ are the eigenvalues of $\mathbf{A}$ and   $\mathbf{B}$, respectively. 
\end{theorem}

Following the same reasoning as in the proof of proposition, it is sufficient to show that for any $\mathbf{u}^\ast\in\mathcal{U}^\ast$, the Jacobian matrix (recall from eq. \eqref{Jacobian_star})
\begin{equation*}
\mathrm{D}\tilde{g}(\mathbf{u}^\ast)= \begin{bmatrix}
\mathbf{W}_R-\alpha\nabla^2 F_c^\ast &-\alpha\mathbf{I}
\\
\left(\mathbf{W}_C-\mathbf{I}\right)\nabla^2 F_c^\star & \mathbf{W}_C
\end{bmatrix},
\end{equation*}
has an eigenvalue with absolute value strictly greater than $1$; 
proving that such  eigenpair is also a member of $\sigma(\mathrm{D}{g}(\mathbf{u}^\ast))$ follows equivalent steps as in the proof of proposition and thus is omitted.   		
%
%
Decompose $\mathrm{D}\tilde{g}(\mathbf{u}^\ast)$ as 
\begin{equation}
\label{Jacobian_decomposition}
\mathrm{D}\tilde{g}(\mathbf{u}^\ast)= \underbrace{\begin{bmatrix}
	\mathbf{I}-\alpha\nabla^2 F_c^\ast &-\alpha\mathbf{I}
	\\
	\mathbf{0} & \mathbf{I}
	\end{bmatrix}}_{\triangleq \mathbf{Q}}
+
\underbrace{\frac{1}{t}\begin{bmatrix}
	\widetilde{\mathbf{W}}_R-\mathbf{I}&\mathbf{0}
	\\
	(\widetilde{\mathbf{W}}_C-\mathbf{I})\nabla F_c^\star  & \widetilde{\mathbf{W}}_C-\mathbf{I}
	\end{bmatrix}}_{\triangleq \mathbf{P}_t},
\end{equation}
where $\widetilde{\mathbf{W}}_R\triangleq \widetilde{\mathbf{R}}\otimes \mathbf{I}_m$ and $\widetilde{\mathbf{W}}_C\triangleq \widetilde{\mathbf{C}}\otimes \mathbf{I}_m$. Eq. \eqref{Jacobian_decomposition} reads the Jacobian matrix $\mathrm{D}\tilde{g}(\mathbf{u}^\ast)$ as a variation of $\mathbf{Q}$ by perturbation $\mathbf{P}_t$. 
For any $\mathbf{u}^\ast\in\mathcal{U}^\ast$, the spectrum of $\mathbf{Q}$ consists of $n\cdot m$ counts of $1$ along with the eigenvalues of $\mathbf{I}-\alpha\nabla^2 F_c^\ast$, which contains a real eigenvalue   $\lambda_1\geq 1+\alpha\delta/(mn)$, since $\boldsymbol{\theta}^\ast\in \Theta^\ast_{ss}$ .
Theorem \ref{th:matpert} guarantees that the spectrum variation of any perturbed arbitrary non-normal matrices is bounded by the norm of the perturbation matrix.  Thus it is sufficient to show that the perturbed $\lambda_1$, as a member of $\sigma(\mathrm{D}\tilde{g}(\mathbf{u}^\ast))$, is strictly greater than $1$.

Applying Theorem~\ref{th:matpert} gives the following sufficient conditions: denote $\tilde{d}\triangleq 2mn$,
\begin{equation}\label{eq:matcond}
\left(\|\mathbf{Q}+\mathbf{P}_t\|+\|\mathbf{Q}\|\right)^{1-1/\tilde{d}} \|\mathbf{P}_t\|^{1/\tilde{d}}  < 2\alpha\delta/\tilde{d}.
\end{equation}
By sub-additivity of the matrix norm, it is sufficient for \eqref{eq:matcond} that
\begin{equation}\label{eq:matcond2}
\left(\|\mathbf{P}_t\|+2\|\mathbf{Q}\|\right)^{1-1/\tilde{d}} \|\mathbf{P}_t\|^{1/\tilde{d}}  \leq \frac{\alpha\delta}{\tilde{d}}.
\end{equation}
Since each $\nabla f_i$ is Lipschitz continuous (cf. Assumption \ref{P_assumption}), there exist constants $C_Q>0$ and $C_P>0$ such that $\max_{\mathbf{u}^\ast\in\mathcal{U}^\ast} \|\mathbf{Q}\|\leq C_Q$ and  $\max_{\mathbf{u}^\ast\in\mathcal{U}^\ast} \|\mathbf{P}_t\|\leq C_P/t$. It is not difficult to show that a sufficient  condition for \eqref{eq:matcond2} is
\begin{equation}
\label{t_condition}
t\geq \frac{(C_P+2C_Q)^{\tilde{d}-1}C_P}{(\alpha\delta/\tilde{d})^{\tilde{d}}},\qquad  \tilde{d}=2mn.
\end{equation}

%
%
%

\begin{proposition}
	\label{CSS_is_in_unstable_set_extended}
	Let Assumptions \ref{matrix_C_R} and \ref{Strict_Saddles_delta_assumption} hold, and matrices $\mathbf{R}$ and $\mathbf{C}$ be chosen according to Assumption \ref{diag_dominant_RC}, with $t$  satisfying \eqref{t_condition}.
	Then,  any consensual strict saddle point  is an unstable fixed point of $g$, i.e.,
	$\mathcal{U}^\ast\subseteq\mathcal{A}_g$,
	with $\mathcal{A}_g$ and $\mathcal{U}^\ast$    defined in    \eqref{def_unstable_fixpoint} and \eqref{def_U_g}, respectively. 
\end{proposition}


Note that above proposition  ensures $\mathcal{U}^\ast\subseteq\mathcal{A}_g$ under \eqref{t_condition} and given  step-size $\alpha$. Convergence of the sequence is proved under \eqref{alpha_conds_assymptotic_conv_final} and \eqref{eq:step_diff} for the step-size  $\alpha$. However,  \eqref{alpha_conds_assymptotic_conv_final}  may not hold for some large $t$ (there can be instances where the set of step-size satisfying conditions \eqref{t_condition} and \eqref{alpha_conds_assymptotic_conv_final} is empty). Hence, when $m>1$, the statement in Theorem \ref{Main_SOS_guarantee_long} is conditioned to the convergence of the algorithm.

\vspace{1cm}

\bibliographystyle{siamplain}
\bibliography{references}

\end{document}